%% file: cotti_arXiv_final.tex
\documentclass[12pt,reqno]{amsart}
\makeatletter
\@namedef{subjclassname@2020}{\textup{2020} Mathematics Subject Classification}
\makeatother
\usepackage{amsmath}
\usepackage{amsthm}
\usepackage{amsfonts}
\usepackage{amssymb}
\usepackage{mathrsfs}
\usepackage{graphicx}
\usepackage[T1]{fontenc}
\usepackage[OT2,T1]{fontenc}
\usepackage[latin1]{inputenc} 
\usepackage{lmodern}
\usepackage[all,cmtip]{xy}
\usepackage{cancel}
\usepackage{bbm}
\usepackage[vcentermath]{youngtab}
\usepackage{ stmaryrd }
\usepackage{calc}
\usepackage{empheq}
\usepackage{multicol} 
\usepackage{xcolor}
\usepackage{setspace}
\usepackage{epigraph}
\usepackage{mathdots}
\usepackage{braket}
\usepackage{longtable}
\usepackage{float}

\usepackage[pdftex,     
          plainpages=false,   
           breaklinks=true,    
           colorlinks=true,
           pdftitle=My Document
           pdfauthor=My Good Self 
           colorlinks=true,
	    urlcolor=blue,
	    citecolor=red,
	    linkcolor=blue
          ]{hyperref}

\usepackage[a4paper,top=4cm,bottom=4cm,left=2.8cm,right=2.8cm]{geometry}

\newtheorem{thm}{Theorem}[section]
\newtheorem{cor}[thm]{Corollary}
\newtheorem{lem}[thm]{Lemma}
\newtheorem{prop}[thm]{Proposition}
\newtheorem{conj}[thm]{Conjecture}

\numberwithin{equation}{section}

\theoremstyle{definition}
\newtheorem{defn}[thm]{Definition}
\newtheorem{rem}[thm]{Remark}
\newtheorem*{example}{Example}

\newtheorem*{conv}{Convention}

\newlanguage\fakelanguage
\newcommand\cyr{\fontencoding{OT2}\fontfamily{wncyr}\selectfont
   \language\fakelanguage}
\DeclareTextFontCommand{\textcyr}{\cyr}

\numberwithin{equation}{section}

\usepackage{chngcntr}
\counterwithin{figure}{section}
\counterwithin{table}{section}

\DeclareMathOperator{\Hom}{Hom}

\usepackage{calligra}

\makeatletter
\newsavebox{\@brx}
\newcommand{\llangle}[1][]{\savebox{\@brx}{\(\m@th{#1\langle}\)}%
  \mathopen{\copy\@brx\kern-0.5\wd\@brx\usebox{\@brx}}}
\newcommand{\rrangle}[1][]{\savebox{\@brx}{\(\m@th{#1\rangle}\)}%
  \mathclose{\copy\@brx\kern-0.5\wd\@brx\usebox{\@brx}}}
\makeatother

\usepackage{framed,color}
\definecolor{shadecolor}{rgb}{0.67, 0.9, 0.93} 
\definecolor{ambra}{rgb}{1.0, 0.75, 0.0}
\definecolor{ametista}{rgb}{0.6, 0.4, 0.8}
\definecolor{auburn}{rgb}{0.43, 0.21, 0.1}
\definecolor{ballblue}{rgb}{0.13, 0.67, 0.8}
\definecolor{cadmiumgreen}{rgb}{0.0, 0.42, 0.24}
\definecolor{candypink}{rgb}{0.89, 0.44, 0.48}
\definecolor{caribbeangreen}{rgb}{0.0, 0.8, 0.6}

\newcommand{\clr}{\textcolor[rgb]{1.00,0.00,0.00}}

\newcommand{\bsh}{\begin{shaded}}
\newcommand{\esh}{\end{shaded}}

\newcommand{\bsy}{\boldsymbol}
\newcommand{\beq}{\begin{equation}}
\newcommand{\eneq}{\end{equation}}
\newcommand{\puqed}{\pushQED{\qed}}
\newcommand{\poqed}{\popQED}

\newcommand{\edc}{\widehat\nabla}
\newcommand{\derz}[1]{\frac{\der#1}{\der z}}

\newcommand{\rqh}{|_{\substack{{\bf Q}=1\\ \hbar=1}}}

\def\R{{\mathbb R}}
\def\N{{\mathbb N}}
\def\C{{\mathbb C}}
\def\Z{{\mathbb Z}}
\def\Q{{\mathbb Q}}
\def\P{{\mathbb P}}
\def\Fb{{\mathbb F}}
\def\L{{\mathbb L}}

\let\mc\mathcal

\let\al\alpha
\let\bt\beta
\let\dl\delta

\let\eps\varepsilon
\let\gm\gamma
\let\Gm\Gamma
\let\ka\kappa
\let\la\lambda
\let\La\Lambda

\let\phi\varphi
\let\si\sigma
\let\Si\Sigma
\let\Sig\varSigma

\let\ups\upsilon
\let\Ups\Upsilon
\let\om\omega
\let\Om\Omega

\let\der\partial

\let\geq\geqslant

\let\leq\leqslant

\let\bm\boldsymbol

\makeatletter
\def\moverlay{\mathpalette\mov@rlay}
\def\mov@rlay#1#2{\leavevmode\vtop{%
   \baselineskip\z@skip \lineskiplimit-\maxdimen
   \ialign{\hfil$\m@th#1##$\hfil\cr#2\crcr}}}
\newcommand{\charfusion}[3][\mathord]{
    #1{\ifx#1\mathop\vphantom{#2}\fi
        \mathpalette\mov@rlay{#2\cr#3}
      }
    \ifx#1\mathop\expandafter\displaylimits\fi}
\makeatother

\newcommand{\cupdot}{\charfusion[\mathbin]{\cup}{\cdot}}

\usepackage{scalerel,stackengine}  
\stackMath
\newcommand\reallywidehat[1]{%
\savestack{\tmpbox}{\stretchto{%
  \scaleto{%
    \scalerel*[\widthof{\ensuremath{#1}}]{\kern.1pt\mathchar"0362\kern.1pt}%
    {\rule{0ex}{\textheight}}
  }{\textheight}%
}{2.4ex}}%
\stackon[-6.9pt]{#1}{\tmpbox}%
}
\begin{document}

\title[Cyclic Stratum, Borel-Laplace $(\bm\alpha,\bm\beta)$-multitransforms]{Cyclic stratum of Frobenius manifolds, Borel-Laplace $(\bm \alpha,\bm \beta)$-multitransforms, and integral representations of solutions of Quantum Differential Equations
}
\author[Giordano Cotti]{Giordano Cotti}
\address{Faculdade de Ci\^encias da Universidade de Lisboa, Grupo de F\'isica Matem\'atica, Campo Grande Edif\'icio C6, 1749-016 Lisboa, Portugal}
\email{gcotti@fc.ul.pt, gcotti@sissa.it}
\maketitle

\begin{abstract}
In the first part of this paper, we introduce the notion of {\it cyclic stratum} of a Frobenius manifold $M$. This is the set of points of the extended manifold $\C^*\times M$ at which the unit vector field is a cyclic vector for the isomonodromic system defined by the flatness condition of the extended deformed connection. The study of the geometry of the complement of the cyclic stratum is addressed. We show that at points of the cyclic stratum, the isomonodromic system attached to $M$ can be reduced to a scalar differential equation, called the {\it master differential equation} of $M$. In the case of Frobenius manifolds coming from Gromov-Witten theory, namely quantum cohomologies of smooth projective varieties, such a construction reproduces the notion of quantum differential equation.

In the second part of the paper, we introduce two multilinear transforms, called {\it Borel-Laplace $(\bm \al,\bm\bt)$-multitransforms}, on spaces of Ribenboim formal power series with exponents and coefficients in an arbitrary finite dimensional $\C$-algebra $A$. When $A$ is specialized to the cohomology of smooth projective varieties, the integral forms of the Borel-Laplace $(\bm \al,\bm\bt)$-multitransforms are used in order to rephrase the Quantum Lefschetz Theorem. This leads to explicit Mellin-Barnes integral representations of solutions of the quantum differential equations for a wide class of smooth projective varieties, including Fano complete intersections in projective spaces.

In the third and final part of the paper, as an application, we show how to use the new analytic tools, introduced in the previous parts, in order to study the quantum differential equations of Hirzebruch surfaces. For Hirzebruch surfaces diffeomorphic to $\P^1\times \P^1$, this analysis reduces to the simpler quantum differential equation of $\P^1$. For Hirzebruch surfaces diffeomorphic to the blow-up of $\P^2$ in one point, the quantum differential equation is integrated via Laplace $(1,2;\frac{1}{2},\frac{1}{3})$-multitransforms of solutions of the quantum differential equations of $\P^1$ and $\P^2$, respectively. This leads to explicit integral representations for the Stokes bases of solutions of the quantum differential equations, and finally to the proof of Dubrovin Conjecture for all Hirzebruch surfaces.
\end{abstract}

\newpage
\tableofcontents

\section{Introduction}

\noindent1.1$\,\,\,${\bf Background.} In the last decades, we have been witnessing a growing and fruitful interaction between theoretical physics and various branches of geometry, leading to new developments in both disciplines. {\it Enumerative geometry} - an old subject and an active field in the 19th century- has been revolutionized by new ideas from the physics of string theory. After the categorical axiomatization of physical theories of quantum fields \cite{atiyahtop,atiyahintro,segal}, the emergence of new mathematical objects was noticed. 
In such an inspiring context, rich structures known as Frobenius manifolds naturally arise, together with the construction of several invariants of symplectic and algebraic varieties.

The notion of {\it Frobenius manifolds} was introduced by B.\,Dubrovin, who first recognized its emergence in the study of classification of two dimensional topological field theories \cite{dubronapoli,dubro1,dubro0}. A Frobenius manifold consists\footnote{Precise definitions will be given in the main body of the paper.} of a complex manifold $M$ whose tangent spaces admit an associative, commutative, and unital algebra structure $(T_pM,\circ_p)$, holomorphically depending on the point $p\in M$. The structure is further enriched with a non-degenerate symmetric bilinear form $\eta$, whose Levi-Civita connection is flat, and which is compatible with the product, that is
\[\eta(Y\circ W,Z)=\eta(Y,W\circ Z),
\]for any local vector fields $Y,W,Z$ on $M$. This condition makes $(T_pM,\circ_p,\eta_p)_{p\in M}$ a family of Frobenius algebras. Pretty soon, it was understood that  Frobenius manifolds are a unifying notion in mathematics. These structures play a central role in mirror symmetry, theory of unfolding spaces of singularities, and enumerative geometry \cite{manin,hertling,sabbah}. Remarkably enough, results proved for classes of Frobenius manifolds emerging in a certain mathematical theory turn out to be valid in general. This {\it universality} of Frobenius manifolds usually leads to unexpected connections between the aforementioned mathematical theories \cite{dubro3}.

{\it Quantum cohomology}, introduced by E.\,Witten \cite{witten} and C.\,Vafa \cite{vafa} in their study of topological nonlinear sigma model, is one of the most interesting example of Frobenius manifold, associated with any complex smooth projective variety $X$, or a more general compact symplectic manifold  \cite{dubro1,kon,manin}. 
From the physical point of view, the space $X$ is the target of two-dimensional fields, and the Frobenius algebras that arise are a highly non-linear deformation of the classical cohomological ring $H^\bullet(X,\mathbb C)$. 
If the classical cohomology ring of a variety encodes information about the intersections of its subvarieties, the non-functorial construction of quantum cohomology is an instrument to understand how they are related by {\it rational} (or, in the general symplectic case, {\it pseudo-holomorphic}) curves. This information
is codified in the {\it Gromov--Witten invariants} \cite{gromov,witten,witten91}, used to define the quantum perturbation of the product.
Gromov-Witten invariants count curves on $X$: for each $\beta\in H_2(X,\Z)/$torsion, and cycles $Z_1,\dots, Z_n\subseteq X$ in general position, the Gromov-Witten invariant\footnote{Here $PD(\al)$ denotes the Poincar\'e dual class of $\al$.}
\[\langle PD(Z_1), \dots, PD(Z_n)\rangle_{g,n,\bt}^X\in\Q
\]heuristically equals the number of curves $C\subseteq X$, of genus $g$, with homology class $[C]=\bt$, and intersecting all the cycles $Z_i$. Consider the generating function 
\[F_0^X(\bm \gm)=\sum_{n=0}^\infty\sum_{\bt}\frac{1}{n!}\langle\underbrace{\bm\gm,\dots,\bm\gm}_{n\text{ times}}\rangle^X_{0,n,\bt},\qquad \bm\gm\in H^\bullet(X,\C),
\]of genus 0 Gromov-Witten invariants of $X$, and assume that this sum is convergent on a non-empty domain $\Om\subseteq H^\bullet(X,\C)$. The quantum cohomology $QH^\bullet(X)$ is the Frobenius manifold structure on $\Om$, the flat metric $\eta$ being given by the Poincar\'e pairing 
\[\eta(Y,W):=\int_XY\cup W,
\]for any local vector fields\footnote{The tangent space $T_p\Om$ is canonically identified with $H^\bullet(X,\C)$ for any $p\in\Om$. Thus the $\cup$-product $Y\cup W$ of local vector fields is well-defined.} $Y,W$ on $\Om$, 
and the product $Y\circ W$ of vector fields being defined by the identity
\[\eta(Y\circ W,Z)=(YWZ)\,F_0^X,
\]for arbitrary flat local vector fields $Y,W,Z$ on $\Om$.

\noindent1.2$\,\,\,${\bf The Main Problem.} At the core of the analytic theory of Frobenius manifolds, there is the local identification of semisimple\footnote{A point $p\in M$ is semisimple if the Frobenius algebra $(T_pM,\circ_p,\eta_p)$ is with no nilpotents.} points $p\in M$ with the parameters of isomonodromic deformations of ordinary differential equations with rational coefficients. Such an identification - one of the main points of the theory of Dubrovin- was originally established in \cite{dubro1,dubro0,dubro2}, and subsequently extended in \cite{CDG0,CDG,CG1,CG2}. 

In this paper, we mainly consider the example of analytic Frobenius manifolds given by the quantum cohomology $QH^\bullet(X)$ of a complex smooth projective variety $X$, see  \cite{dubro1,kon,manin}. In such a case, points $p\in QH^\bullet(X)$ are parameters of isomonodromic deformations of a linear system of differential equations of the form
\beq
\label{de}
\frac{\der}{\der z}\zeta(z,p)=\left({\bm{\mc U}}(p)+\frac{1}{z}{\bm\mu}(p)\right)\zeta(z,p).
\eneq
Here $\zeta$ is a $z$-dependent vector field of $QH^\bullet(X)$, whereas $\bm{\mc U}$ and $\bm\mu$ are $(1,1)$-tensors on $QH^\bullet(X)$: the first\footnote{Precise definitions will be given in the main body of the paper.} is the operator of quantum multiplication by the Euler vector field - a distinguished vector field on $QH^\bullet(X)$ which equals the first Chern class $c_1(X)$ along the locus of small quantum cohomology - the second, called \emph{grading operator}, keeps track of the non-vanishing degrees of $H^\bullet(X,\C)$. 

Equation \eqref{de} is a rich object associated with the variety $X$: it encapsulates information not only about its {\it Gromov-Witten theory}, but also (conjecturally) about its {\it topology}, its {\it algebraic geometry}, and their mutual relations. The study of the monodromy of solutions of \eqref{de} is the way to disclose such an amount of information, see \cite{dubro0,gamma1,CDG1}. In this paper we address the following 

{{\bf Main Problem:} \it to find integral representations of solutions of \eqref{de} for Fano complete intersections in Fano varieties}. 

\noindent We split the Main Problem in two parts:
\begin{enumerate}
\item to reduce the system of differential equations \eqref{de} to a distinguished scalar linear differential equation, the \emph{master differential equation;}
\item to find integral representations of solutions of master differential equations.
\end{enumerate}
The study of these questions leads us to introduce some relevant notions, both in the analytic theory of Frobenius manifolds and in the theory of integral transforms. The first three ingredients are the notions of \emph{cyclic stratum}, \emph{master differential equations} and \emph{master functions} of a Frobenius manifold. The second new analytical tool is a pair of integral multilinear transforms of functions, that we call \emph{Borel-Laplace $(\bm \al,\bm \bt)$-multitansforms}. We are going to briefly outline these objects.

\noindent 1.3$\,\,\,${\bf Master functions and master differential equations.}  The rich geometry of a Frobenius manifold $M$ is (almost) completely encoded in 
integrability conditions of the \emph{extended deformed connection} or \emph{first structural connection} of $M$ \cite{dubro1,dubro2,manin}. This is a flat meromorphic connection $\edc$ defined on the pullback $\pi^*TM$ of the tangent bundle of $M$ on the extended manifold $\widehat M:=\C^*\times M$, by the natural projection $\pi\colon \widehat M\to M$. Equation \eqref{de} is equivalent to the equation
\beq\label{de2}
\edc_{\frac{\der}{\der z}}\xi=0,\quad \xi\in\Gamma(\pi^*T^*M),
\eneq
the one-form $\xi$ and the vector field $\zeta$ being identified via a flat metric $\eta$ on $M$. We call \emph{master function} at $p\in M$ any function\footnote{Here $\widetilde {\C^*}$ denotes the universal cover of $\C^*$.} $\Phi_{\xi}\in\mc O(\widetilde {\C^*})$ of the form 
\[\Phi_{\xi}(z)=z^{-\frac{d}{2}}\langle \xi(z,p),e(p)\rangle,
\]where $\xi$ is as in \eqref{de2}, and $d$ is the \emph{charge} of the Frobenius manifold $M$. 

In the first part of the paper, we address the problem of reducing the system of differential equations \eqref{de2}  to a scalar differential equation, whose coefficients depend on the point $p\in M$. This is a well known problem in the theory of ordinary differential equations, equivalent to the choice of a \emph{cyclic vector} \cite[Lemma II.1.3]{delignede}. 
On Frobenius manifold, however, we have a \emph{natural} candidate, namely the unit vector field $e\in\Gamma(TM)$.

In Section \ref{22.05.17-1} we introduce the \emph{cyclic stratum} $\widehat M^{\rm cyc}\subseteq\widehat M$ defined as the set of points $(z,p)$ at which the iterated covariant derivatives
\beq\label{cyc}
e,\quad\edc_\derz{}e,\quad\edc^2_\derz{}e,\dots,\quad \edc^{n-1}_\derz{}e,\quad n:=\dim_\C M,
\eneq
define a basis of the fiber $\pi^*TM|_{(z,p)}$. 
The complement of $\widehat M^{\rm cyc}$ in $\P^1\times M$ admits a natural stratification, whose study is addressed in Section \ref{strat}. A particular role is played  by the $\mc A_\La$-stratum of $M$, defined as the set of points $p\in M$ such that $$\C^*\times\{p\}\subseteq \widehat M\setminus \widehat M^{\rm cyc}.$$
Introducing the cyclic coframe $\om_0\dots,\om_{n-1}\in\Gamma(\pi^*T^*M)$ as the dual frame of \eqref{cyc}, the system of differential equations \eqref{de2}, specialized at points $p\in M\setminus \mc A_\La$, reduces to a scalar differential equation - the {\it master differential equation} - in the function $\langle\xi,e\rangle$. 
Hence, at points $p\in M\setminus\mc A_\La$, we obtain a one-to-one correspondence
\[\left\{\begin{aligned}
&\text{Solutions $\xi$ of the system}\\
&\text{ \eqref{de2} specialized at $p$}\\
\end{aligned}\right\}\Longleftrightarrow
\left\{
\text{Master functions $\Phi_\xi$'s at $p$}
\right\}.
\]See Theorems \ref{thmde} and \ref{niso}. Thus, if integral representations for a basis of master functions are found, the Main Problem is \emph{solved} at points in $M\setminus \mc A_\La$.

Some motivational comments for introducing these new tools are in order. The notions of master functions and master differential equations define analogs, for an arbitrary Frobenius manifold, of well-known objects in Gromov-Witten and quantum cohomology theories. Namely, in the case of quantum cohomology the components of Givental's $J$-function (w.r.t. an \emph{arbitrary} cohomology basis) define a generating set of master functions. Moreover, the master differential equation is (up to re-scaling of the unknown function) a quantum differential equation as defined e.g. in \cite[Section 10.3]{cox}, see Section \ref{sec5}. In our opinion the concepts of cyclic stratum, master functions, and master differential equations may represent relevant notions in the analytic theory of Frobenius manifolds. For example, any contingent relations with the geometry of distinguished subsets of Frobenius manifolds (e.g.\,\,bifurcation diagram, Maxwell stratum, caustic) deserve further investigations. In that regard, 
it would be interesting to study relations with results of \cite{CDG0,CDG}, concerning the isomonodromic description of Frobenius manifolds at  semisimple coalescing points. This point will be addressed in a future publication.

\noindent 1.4$\,\,\,${\bf Borel-Laplace multitransforms.} In Section \ref{BL}, we introduce a pair of multilinear transforms in both a formal and an analytical setting.

For $h\in\N^*$, and a given $h$-tuple $\bm\kappa\in (\C^*)^h$, we introduce a ring $\mathscr F_{\bm \kappa}(A)$ of Ribenboim generalized power series \cite{riben2,riben1} with both coefficients and exponents in a finite dimensional, commutative, associative, and unitary $\C$-algebra $A$. The numbers $\kappa_i$'s play a role of ``weights'' for the exponents of the power series. In such a formal setting, given $\bm\al,\bm\bt\in(\C^*)^h$, we introduce the {\it Borel-Laplace $(\bm\al,\bm\bt)$-multitransforms} as two $A$-multilinear maps rescaling the weights
\begin{align*}&\mathscr B_{\bsy\alpha,\bsy\beta}\colon \bigotimes_{j=1}^h\mathscr F_{\kappa_j}(A)\to \mathscr F_{\bm\al^{-1}\cdot\bm\bt^{-1}\cdot\bm\kappa}(A),&&\bm\al^{-1}\cdot\bm\bt^{-1}\cdot\bm\kappa:=\left(\frac{\kappa_1}{\al_1\bt_1},\dots,\frac{\kappa_h}{\al_h\bt_h}\right),\\
&\mathscr L_{\bsy\alpha,\bsy\beta}\colon \bigotimes_{j=1}^h\mathscr F_{\kappa_j}(A)\to \mathscr F_{\bm\al\cdot\bm\bt\cdot\bm\kappa}(A),&&\bm\al\cdot\bm\bt\cdot\bm\kappa:=(\al_1\bt_1\kappa_1,\dots,\al_h\bt_h\kappa_h).
\end{align*}
See Sections \ref{fka} and \ref{formal} for precise definitions.

In the analytical setting, given $h$ functions $\Phi_1,\dots,\Phi_h\colon \widetilde{\C^*}\to A$, we define their Borel-Laplace $(\bm\al,\bm\bt)$-multitransforms by
\begin{align*}\mathscr B_{\bm\al,\bm\bt}[\Phi_1,\dots,\Phi_h](z):=\frac{1}{2\pi i}\int_\gamma\prod_{j=1}^h\Phi_j\left(z^\frac{1}{\al_j\bt_j}\la^{-\bt_j}\right)e^\la\frac{d\la}{\la},\\
\mathscr L_{\bsy{\alpha},\bsy{\beta}}\left[\Phi_1,\dots,\Phi_h\right](z):=\int_0^\infty\prod_{i=1}^h\Phi_i(z^{\alpha_i\beta_i}\lambda^{\beta_i})e^{-\lambda}d\lambda,
\end{align*}
provided that the integrals exist. The contour $\gamma$ is a Hankel-type contour beginning from $-\infty$, circling the origin once in the positive direction, and returning to $-\infty$ (see Figure \ref{gammahankel}).

\noindent 1.5$\,\,\,${\bf Main results.} Consider a Fano smooth projective variety $X$, and let $\iota\colon Y\to X$ be a Fano subvariety defined as the zero locus of a regular section of a vector bundle $E\to X$. 
The classical cohomology groups $H^k(Y,\C)$ can be (partially) recovered by the cohomology groups $H^k(X,\C)$ by Lefschetz Hyperplane Theorem. Quantum Lefschetz Theorem  is a quantum improvement of the classical result: it describes how to reconstruct the Gromov-Witten theory of $Y$ starting from the Gromov-Witten theory of $X$ \cite{lee,coatgiv,qLcoates}.

In this paper, by using Quantum Lefschetz Theorem, we give explicit integral representations of master functions of $Y$ in terms of Laplace $(\bm\al,\bm\bt)$-multitransforms of master functions of the ambient space $X$ under the following assumptions on $X$ and $E$:
\begin{enumerate}
\item[\bf Case 1.] We assume that $E$ is a direct sum of fractional powers of the determinant bundle $\det TX$ of $X$;
\item[\bf Case 2.] We assume that $X=X_1\times\dots\times X_h$ is a product of Fano varieties $X_i$'s, and that $E$ is the external tensor product of fractional powers of the determinant bundles $\det TX_i$.
\end{enumerate}

Our first main result concerns Case 1. Our Theorem \ref{TH1} asserts that any master function of $Y$, at points $\iota^*\delta\in H^2(Y,\C)$ of its small quantum cohomology, can be expressed in terms of iterated Laplace $(\al,\bt)$-transforms (simple transforms of a single function) of master functions of $X$ at the point $\delta\in H^2(X,\C)$. More precisely, if $E=\bigoplus_{j=1}^rL^{\otimes d_j}$, and $\det TX=L^\ell$ for an ample line bundle $L$, then any master function of $Y$ at $\iota^*\delta$ is a $\C$-linear combination of integrals of the form
\begin{align*}e^{-c_\delta z}\mathscr L_{\frac{\ell-\sum_{i=1}^s d_i}{d_s},\frac{d_s}{\ell-\sum_{i=1}^{s-1} d_i}}\circ\dots\circ\mathscr L_{\frac{\ell-d_1-d_2}{d_2},\frac{d_2}{\ell-d_1}}\circ\mathscr L_{\frac{\ell-d_1}{d_1},\frac{d_1}{\ell}}[\Phi]\\
=
e^{-c_\delta z}\int_0^\infty\dots\int_0^\infty\Phi\left(z^{\frac{\ell-\sum_{j=1}^rd_j}{\ell}}\prod_{i=1}^r\zeta_i^\frac{d_i}{\ell}\right)e^{-\sum_{i=1}^r\zeta_i}d\zeta_1\dots d\zeta_r,
\end{align*}where $\Phi$ is a master function of $X$ at $\delta$, and $c_\delta\in\C$ is a complex number depending on $\delta$.

Our second main result concerns Case 2. In particular, Theorem \ref{TH2} asserts that any master function of $Y$, at points $\iota^*\delta\in H^2(Y,\C)$ of the small quantum locus, can be expressed in terms of Laplace $(\bm \al,\bm\bt)$-multitransforms of master functions of $X_j$ at the point $\delta_j\in H^2(X,\C)$, where
\[\delta=\sum_{j=1}^h1\otimes\dots\otimes\delta_j\otimes \dots\otimes 1.
\] More precisely, if $E=\boxtimes_{j=1}^hL_j^{\otimes d_j}$ and $\det TX_j=L_j^{\ell_j}$ for ample line bundles $L_j$,  
then any master function of $Y$ at $\iota^*\delta$ is a $\C$-linear combination of integrals of the form
\[e^{-c_\delta z}\mathscr L_{\bm\al,\bm\bt}[\Phi_1,\dots,\Phi_h](z)=e^{-c_\delta z}\int_0^\infty\prod_{j=1}^h\Phi_j\left(z^{\frac{\ell_j-d_j}{\ell_j}}\la^\frac{d_j}{\ell_j}\right)e^{-\la}d\la,
\]
where $(\bm\al,\bm\bt)=(\frac{\ell_1-d_1}{d_1},\dots,\frac{\ell_h-d_h}{d_h};\frac{d_1}{\ell_1},\dots,\frac{d_h}{\ell_h})$, $\Phi_j$ is a master function of $X_j$ at $\delta_j$, and $c_\delta\in\C$ is a complex number depending on $\delta$.

Assumptions of Cases 1 and 2 are clearly satisfied when the varieties $X$ and $X_j$'s have Picard rank one. Therefore, Theorems \ref{TH1} and \ref{TH2} can be applied to all Fano complete intersections in $\P^n$ and Fano hypersurfaces in products of projective spaces, in order to obtain explicit Mellin-Barnes integral representations of master functions. In particular, if $Y\subseteq \P^{n-1}$ is a Fano complete intersection defined by homogeneous polynomials of degrees $d_1,\dots, d_h$, our Theorem \ref{TH1b} asserts that any master function of $Y$ at $0\in H^\bullet(Y,\C)$ is a linear combination of {\it one}-dimensional Mellin-Barnes integrals 
\[
G_j(z)=\frac{e^{-cz}}{2\pi\sqrt{-1}}\int_\gamma\Gamma(s)^n\prod_{k=1}^h\Gamma\left(1-d_ks\right)z^{-(n-\sum_{k=1}^hd_k)s}\phi_j(s)ds,\quad j=0,\dots, n-1,
\]where $c\in\Q$, $\gamma$ is a parabola (of the form ${\rm Re}\,s=-\rho_1({\rm Im}\,s)^2+\rho_2$, for suitable $\rho_1,\rho_2\in\R_+$) encircling the poles of the factor $\Gamma(s)^n$ and separating them from the poles of the factors $\Gamma\left(1-d_ks\right)$, and the function $\phi_j(s)$ are defined by
 \begin{empheq}[left={\phi_j(s):=}\empheqlbrace]{align*} 
 \exp\left({2\pi\sqrt{-1}js}\right),\quad n\text{ even,}\\
 \exp\left({2\pi\sqrt{-1}js}+\pi\sqrt{-1}s\right),\quad n \text{ odd}.
 \end{empheq}
In the case of a Fano hypersurface $Y\subseteq \P^{n_1-1}\times \dots\times \P^{n_h-1}$ defined by a homogeneous polynomial of multi-degree $(d_1,\dots, d_h)$, then our Theorem \ref{TH2b} asserts that any master function of $Y$ at $0\in H^\bullet(Y,\C)$ is a linear combination of the $h$-dimensional Mellin-Barnes integrals
\[
H_{\bm j}(z):=\frac{e^{-cz}}{(2\pi\sqrt{-1})^h}\int_{\bigtimes\gamma_{i}}\left[\prod_{i=1}^h\Gamma(s_i)^{n_i}\phi_{j_i}^{i}(s_i)\right]\Gamma\left(1-\sum_{i=1}^hs_i\right)z^{-\sum_{i=1}^hd_is_i}ds_1\dots ds_h,
\]where $c\in\Q$, $\gamma_i$ are parabolas (of the form ${\rm Re}\,s_i=-\rho_{1,i}({\rm Im}\,s_i)^2+\rho_{2,i}$, for suitable $\rho_{1,i},\rho_{2,i}\in\R_+$) encircling the poles of the factors $\Gamma(s_i)^{n_i}$, and the functions $\phi_{j_i}^{i}(s_i)$ are defined by
\begin{empheq}[left={\phi_{j_i}^{i}(s_i):=}\empheqlbrace]{align*}
\exp\left({2\pi\sqrt{-1}j_is_i}\right),\quad n_i\text{ even},\\
\exp\left({2\pi\sqrt{-1}j_is_i}+\pi\sqrt{-1}s_i\right),\quad n_i\text{ odd},
\end{empheq}
for any $h$-tuple $\bm j=(j_1,\dots, j_h)$ with $0\leq j_h\leq n_i-1$.

Some comments are in order. Given a Fano variety $X$, Mirror Symmetry provides other kinds of integral representations of solutions of equation \eqref{de2}.\footnote{More precisely, for the equations $\edc_{\frac{\der}{\der t^\al}}\xi=0$, where $t^1,\dots, t^n$ are coordinates on $QH^\bullet(X)$, and not w.r.t. the spectral parameter $z$.}  These are complex oscillating integrals associated with the Landau-Ginzburg models mirror to $X$, see \cite{givehomgeo,give2,givemirtor,eguchi,konens,horivafa}. In these representations the cycles of integration are multi-dimensional\footnote{Notice, for example, that already in the case of $\P^n$ these oscillating  integrals are over $n$-dimensional cycles. On the other hand, one-dimensional Mellin-Barnes integral representations of solutions of the equation \eqref{de} associated with $\P^n$ were obtained in \cite{guzzetti1}. Their asymptotics in sectors of $\widetilde{\C^*}$ is easier to study.}.  This fact typically makes more difficult the study of the asymptotic expansions of solutions, and of the determination of the corresponding validity sectors in $\widetilde{\C^*}$.  Furthermore, let us recall another technical issue which may be faced:  Landau-Ginzburg models may not have enough critical points, and suitable compactification procedures have to be applied in order to recover the right number, see \cite{konny,maximvassily,pechrie}. This could represent a delicate point for the computation of the Stokes bases of solutions of equation \eqref{de}, whose exponential growth is ruled by the critical values of the Landau-Ginzburg potential. 

We believe that one-dimensional Mellin-Barnes integrals of Theorem \ref{TH1b} represent a more advantageous representation of the solutions to the purpose of asymptotic analysis. Moreover, even for multi-dimensional Mellin-Barnes integrals of Theorem \ref{TH2b} the study of their asymptotics is tame: it is equivalent to the study of the asymptotics of one-dimensional generalized Fax\'en integrals 
\[I(\la;c_1,\dots, c_r):=\int_0^\infty\exp\left[-\la\left(x^\mu+\sum_{k=1}^rc_kx^{m_k}\right)\right]dx,\]\[ \text{with }\mu>m_1>m_2>.\dots >m_r>0,
\]which have saddle points whose exponential contributions dominate algebraic terms in the asymptotic expansion. See \cite[Chapter 7]{pkmellinbook}, \cite[Section 5]{pkmellin} for a detailed asymptotic analysis, and also \cite{bur,bak,wright} for some special cases. This will be exemplified in Section \ref{aLap}.

\noindent 1.6$\,\,\,${\bf Dubrovin Conjecture for Hirzebruch surfaces.} Equation \eqref{de} has two singularities: a Fuchsian singularity at $z=0$ and  an irregular singularity at $z=\infty$ of Poincar\'e rank 1. The monodromy of its solutions is quantified by a finite set of matrices:  
\begin{itemize}
\item a monodromy matrix $M_0$, quantifying the monodromy of solutions of \eqref{de} at $z=0$,
\item a Stokes matrix $S$, describing the Stokes phenomenon at $z=\infty$, 
\item and a central connection matrix $C$ gluing the monodromy data $M_0$ and $S$ at the two singularities.
\end{itemize}
Remarkably, the monodromy data  define a sort of ``system of coordinates'' in the space of solutions of WDVV equations: from the knowledge of their numerical values, the whole Frobenius manifold structure can be reconstructed via a Riemann-Hilbert problem \cite{dubro1,dubro2,guzzetti2}.

In \cite{dubro0}, B.\,Dubrovin formulated an intriguing conjecture concerning the geometrical meaning of the numerical values of the monodromy data of quantum cohomologies of Fano varieties. In the \emph{qualitative} part of the conjecture,  for a given Fano variety $X$, the semisimplicity condition of $QH^\bullet(X)$ is claimed to be equivalent to the existence of full exceptional collections in the derived category $\mc D^b(X)$ of coherent sheaves on $X$. Moreover, in the refined \emph{quantitative} part of the conjecture, formulated in \cite[Conjecture 5.2]{CDG1}, the Stokes and central connection matrices $(S_p,C_p)$ computed at any point $p\in QH^\bullet(X)$ are claimed to be determined by characteristic classes of $X$ and of objects of a full exceptional collection $\frak E_p$ in  $\mc D^b(X)$. 

In particular, the central connection matrix $C_p$ is claimed to equal the matrix associated with the morphism
\begin{align}\label{Drussa}
\textnormal{\textcyr{D}}_X^-\colon K_0(X)_\mathbb C\longrightarrow& H^\bullet(X,\C)\\ 
\nonumber F\xmapsto{\quad\quad}& \frac{(\sqrt{-1})^{\overline{d}}}{(2\pi)^{\frac{d}{2}}}\widehat\Gamma^-_X\exp(-\pi\sqrt{-1}c_1(X)){\rm Ch}(F),
\end{align}
where $d=\dim_\C X$,  $\overline d$ is its residue class modulo 2, $\widehat\Gamma^-_X$ is the characteristic class of $X$ defined by
\begin{align*}\widehat\Gamma^-_X:=\prod_{j=1}^{\dim_\C X}\Gamma(1-\delta_j),\quad
\Gamma(1-t)=\exp\left(\gamma t+\sum_{n=2}^\infty\frac{\zeta(n)}{n}t^n\right),\quad \delta_j\text{ Chern roots of $TX$,}
\end{align*}
and ${\rm Ch}(F)$ is the graded Chern character defined on vector bundles by the formula ${\rm Ch}(V):=\sum_{j=1}^{{\rm rk} V}\exp(2\pi\sqrt{-1}\eps_j)$, $\eps_j$'s being the Chern roots of $V$. The matrix of \textcyr{D}$^-_X$ is computed w.r.t. the exceptional basis $[\frak E_p]$ of $K_0(X)_{\C}$, defined by the $K$-theoretical classes of objects of $\frak E_p$, and an arbitrary\footnote{The choice of a basis of $H^\bullet(X,\C)$ in \eqref{Drussa} corresponds to the choice of a system of flat coordinates on $QH^\bullet(X)$ w.r.t. which the monodromy data $(M_0,S,C)$ are computed.} basis of $H^\bullet(X,\C)$.  
Furthermore, if the central connection matrix $C_p$ is related to the morphism $\textcyr{D}^-_X$ as explained above, then the Stokes matrix $S_p$ \emph{automatically} equals the inverse of the Gram matrix of the Grothendieck-Euler-Poincar\'e $\chi$-pairing on $K_0(X)$ w.r.t. the exceptional basis $[\frak E_p]$, see \cite[Corollary 5.8]{CDG1}.

It is important to stress that the monodromy data $(M_0,S,C)$ are defined up to several choices: the choice of a system of flat coordinates on the Frobenius manifold $QH^\bullet(X)$, choices of normalizations (at both $z=0$ and $z=\infty$) of solutions of equation \eqref{de}, and the choice of an ``admissible ray'' in $\C^*$. Remarkably, all these operations have a geometrical counterpart in derived categories, see \cite[Theorem 5.9]{CDG1}. Deserving special mention is $\Gamma$- conjecture II of \cite{gamma1}: it consists of an equivalent conjectural statement about the central connection matrix, though w.r.t. a choice of a solution in ``Levelt form'' at $z=0$ not \emph{natural} from the point of view of the theory of Frobenius manifolds. See \cite[Section 5.6]{CDG1} for details.

The explicit computation of the monodromy data of quantum cohomologies is typically a rather delicate operation. To the best knowledge of the author, the only cases in which the computation of the complete set of monodromy data $(S,C)$ of equation \eqref{de} has been carried out in all the details (including the determination of the corresponding full exceptional collections) are the cases of projective spaces \cite{dubro2,guzzetti1} and of complex Grassmannians \cite{gamma1,CDG1}. We believe that the main results of the current paper, namely the integral representations described in Theorems \ref{TH1}, \ref{TH2}, \ref{TH1b}, and \ref{TH2b}, will represent a fundamental tool for the development of this study \cite{inprep}.

As an application, in Sections \ref{Df2k} and \ref{Df2k1}, we show how to use the Laplace $(\bm\al,\bm\bt)$-multitransform, and the main results described above, in order to prove the quantitative part of Dubrovin Conjecture for Hirzebruch surfaces \cite{hirzsup}. These are surfaces $\Fb_k$, with $k\in \Z$, defined as the total space of the projective bundle $\P(\mc O\oplus\mc O(-k))$ on $\P^1$. The interest of this example is highlighted by the fact that 
\begin{itemize}
\item only two Hirzebruch surfaces are Fano varieties (namely $\Fb_0$ and $\Fb_{1}$),
\item all others Hirzebruch surfaces are deformation equivalent to either $\Fb_0$ or $\Fb_1$.
\end{itemize}
Results of A.\,Bayer already suggested the unnecessity of the Fano assumption for the validity of the qualitative part of Dubrovin Conjecture, see \cite{bay4}. Moreover, X.\,Hu proved that, in a smooth family of complete varieties, the existence of full exceptional collection on a fiber preserves for the fibers in a neighborhood, see \cite{xu}.  See also \cite[Corollary B]{steve} for an analogue result for arbitrary semiorthogonal decompositions. To the best of our knowledge, the study of the monodromy of the isomonodromic systems \eqref{de} associated with Hirzebruch surfaces, developed in Sections \ref{Df2k} and \ref{Df2k1}, represents the first example in literature which addresses also the quantitative part of Dubrovin Conjecture, in both the non-Fano case and the case of deformations of the complex structures.

The case of Hirzebruch surfaces $\Fb_{2k}$ (resp. $\Fb_{2k+1}$) can be reduced to the single case of $\Fb_{0}=\P^1\times \P^1$ (resp. $\Fb_1={\rm Bl}_{pt}\P^2$). The monodromy data of $QH^\bullet(\Fb_0)$ can be easily reconstructed from the monodromy data of $QH^\bullet(\P^1)$, see Theorem \ref{DubconjF2k}. In the case of $QH^\bullet(\Fb_1)$, the computation is more delicate, and reduces to the study of the quantum differential equation 
\begin{multline*}
(283 z-24)\vartheta^4\Phi+\left(283 z^2-590 z+24\right)\vartheta^3\Phi+ \left(-2264 z^2+192 z+3\right)\vartheta^2\Phi\\
\nonumber-4 z^2 \left(2547 z^2+350 z-104\right)\vartheta\Phi+z^2 \left(-3113 z^3-9924 z^2+1476 z+192\right)\Phi=0,
\end{multline*}
where $\vartheta:=z\frac{d}{dz}$.  In Section \ref{LapF1}, we show that the solutions of  this equation can be expressed as linear combinations of integrals of the form
\[e^{-z}\mathscr L_{(1,2;\frac{1}{2},\frac{1}{3})}[\Phi_1,\Phi_2;z]=e^{-z}\int_0^\infty\Phi_1\left(z^\frac{1}{2}\la^\frac{1}{2}\right)\Phi_2\left(z^\frac{2}{3}\la^\frac{1}{3}\right)e^{-\la}d\la,
\]where $\Phi_1$ and $\Phi_2$ are solutions of quantum differential equations of $\P^1$ and $\P^2$ respectively, that is
\[\vartheta^2\Phi_1=4z^2\Phi_1,\quad \vartheta^3\Phi_2=27z^3\Phi_2.
\] This allows the study of the asymptotics of solutions in sectors of $\widetilde{\C^*}$, to reconstruct the Stokes bases of solutions of the quantum differential equation of $\Fb_1$, and finally to the computation of both Stokes and central connection matrices, see Theorem \ref{CSF1}.

From these results, the quantitative part of the Dubrovin Conjecture is proved for all Hirzebruch surfaces $\Fb_{k}$, by making explicit the exceptional collections in $\mc D^b(\Fb_k)$ which arise from the monodromy data, see Theorems \ref{DubconjF2k} and \ref{DubrconjF2k1}. 

\noindent 1.7$\,\,\,${\bf Plan of the paper.} The paper is organized as follows. In Section \ref{22.05.17-1}, we introduce the notion of \emph{cyclic stratum} in the general context of Frobenius manifolds theory. A first study of the geometry of the cyclic stratum, and its complement in the extended manifold $\C^*\times M$, is addressed.

In Section \ref{sec3}, we recall basic definitions in Gromov-Witten theory, including the definition of the Frobenius manifold structure on the quantum cohomology of a smooth projective variety.
In Section \ref{sec4}, we recall the definitions of topological-enumerative solution of the isomonodromic system \eqref{de}, and also of its monodromy data. We also recall the main properties and natural transformations of the complete set of monodromy data.

In Section \ref{sec5}, we recall the definition of Givental's $J$-function, and we explain how it is related to the space of master functions, see Theorem \ref{topJ} and Corollary \ref{masJ}. We recall the formulation of the Quantum Lefschetz Theorem,  and we obtain an upper bound for the dimension of the space of master functions of a Fano hypersurface of a smooth projective variety $X$, see Theorem \ref{inmasY}.

In Section \ref{BL}, we recall the notion of generalized power series in the sense of P.\,Riben\-boim, and we introduce the ring $\mathscr F_{\bm \kappa}(A)$ of generalized power series with coefficients and exponents in a finite-dimensional $\C$-algebra. We introduce the notions of Borel-Laplace $(\bm\al,\bm\bt)$-multitransforms, in both formal and analytic setting, and we prove the compatibility of the two definitions, see Theorem \ref{BLAF}.

In Section \ref{sec7}, we explain how the $J$-function can be identified (in several ways) with elements of rings of Ribenboim generalized power series. We prove the main results of this paper, Theorems \ref{TH1}, \ref{TH2}, \ref{TH1b}, \ref{TH2b}.

In Section \ref{sec8}, we recall the notions of exceptional collections in derived categories of coherent sheaves, exceptional bases in $K$-theory, their mutations and helices.  We then describe the refined statement of Dubrovin Conjecture, as formulated in \cite{CDG1}. 

In Section \ref{sec9}, we describe the classical and quantum cohomology rings of Hirzebruch surfaces. 

In Section \ref{Df2k}, we explicitly compute the monodromy data of the quantum cohomologies $QH^\bullet(\Fb_{2k})$, and we prove Dubrovin conjecture for Hirzebruch surfaces $\Fb_{2k}$. 

In Section \ref{Df2k1}, we address the study of the quantum differential equations of Hirzebruch surfaces $\Fb_{2k+1}$. We show how to use the Laplace $(1,2;\frac{1}{2},\frac{1}{3})$-multitransform in order to give integral representations of solutions, how to reconstruct Stokes fundamental solutions, and hence how to compute the monodromy data. This leads to a proof of Dubrovin Conjecture for Hirzebruch surfaces $\Fb_{2k+1}$.

\noindent{\bf Acknowledgements.} The author thanks C.\,Bartocci, A.\,Brini, U.\,Bruzzo, G.\,Carlet, B.\,Dubrovin, D. Guzzetti, C.\,Hertling, A.\,Its, P.\,Lorenzoni, D.\,Masoero, M.\,Mazzocco, A.\,T.\,Ricolfi, V.\,Roubtsov, C.\,Sab\-bah, M.\,Smirnov, A.\,Tacchella, A.\,Varchenko, D.\,Yang for very useful discussions, and also an anonymous referee for corrections and suggestions improving the exposition of the paper. The author is thankful to Max-Planck Institute f\"ur Mathematik in Bonn, Germany, where this project was started, for providing excellent working conditions. This research was supported by MPIM (Bonn, Germany), the EPSRC Research Grant EP/P021913/2, and the FCT Project PTDC/MAT-PUR/ 30234/2017 ``Irregular connections on algebraic curves and Quantum Field Theory''.

\section{Cyclic stratum of Frobenius manifolds}\label{22.05.17-1}

\subsection{Frobenius manifolds} Given a complex manifold $M$, we denote by $TM$ (resp. $T^*M$) its holomorphic tangent (resp. cotangent) bundle. If $E$ is a holomorphic vector bundle on $M$, we denote by $\bigodot^kE$ its $k$-th symmetrized tensor power, and by $\Gm(E)$ the vector space of global holomorphic sections of $E$.

\begin{defn}A \emph{Frobenius manifold} structure on a complex manifold $M$ of dimension $n$ is defined by giving
\begin{enumerate}
\item[(FM1)] a symmetric $\mathcal O(M)$-bilinear form $\eta\in\Gamma\left(\bigodot^2T^*M\right)$, called \emph{metric}\footnote{In what follows, we denote by $(-)^\flat$ and $(-)^\sharp$ the \emph{musical isomorphisms} induced by the metric $\eta$. These are the isomorphisms between vector spaces of mixed tensors.  If $v\in\Gm(TM)$, the 1-form $v^\flat\in\Gm(T^*M)$ is defined by $v^\flat(w)=\eta(w,v)$, where $w\in\Gm(TM)$. Conversely, if $\xi\in\Gm(T^*M)$, the vector field $\xi^\sharp\in\Gm(TM)$ is uniquely defined by the identity $\xi(w)=\eta(w,\xi^\sharp)$, where $w\in\Gm(TM)$. Thus, $(-)^\flat\colon \Gm(TM)\to\Gm(T^*M)$ and $(-)^\sharp\colon\Gm(T^*M)\to \Gm(TM)$ are mutually inverse. In components, these operations are also known as ``lowering'' and ``raising'' of indices, respectively. These operations naturally extend to mixed tensors. For example, given a (1,2)-tensor $c\in\Gm(TM\otimes T^*M\otimes T^*M)$, the tensor $c^\flat$ is the (0,3)-tensor defined by $c^\flat(v_1,v_2,v_3)=\eta(v_1,c(v_2,v_3))$, where $v_1,v_2,v_3\in\Gm(TM)$.}, whose corresponding Levi-Civita connection $\nabla$ is flat;
\item[(FM2)] a $(1,2)$-tensor $c\in\Gamma\left(TM\otimes\bigodot^2T^*M\right)$ such that
\begin{enumerate}
\item the induced multiplication of vector fields $X\circ Y:=c(-,X,Y)$, for $X,Y\in\Gamma(TM)$, is \emph{associative}, 
\item $c^\flat\in\Gamma\left(\bigodot^3T^*M\right)$,
\item $\nabla c^\flat\in\Gamma\left(\bigodot^4T^*M\right)$;
\end{enumerate}
\item[(FM3)] a vector field $e\in\Gamma(TM)$, called the \emph{unity vector field}, such that
\begin{enumerate}
\item the bundle morphism $c(-,e,-)\colon TM\to TM$ is the identity morphism,
\item $\nabla e=0$;
\end{enumerate}
\item[(FM4)] a vector field $E\in\Gamma(TM)$, called the \emph{Euler vector field}, such that
\begin{enumerate}
\item $\frak L_Ec=c$,
\item $\frak L_E\eta=(2-d)\cdot \eta$, where $d\in\mathbb C$ is called the \emph{charge} of the Frobenius manifold.
\end{enumerate}
\end{enumerate}
\end{defn}

At any point $p\in M$ the triple $(T_pM,\eta_p, \circ_p)$ is a complex \emph{Frobenius algebra}, namely an associative commutative algebra with unity whose product is compatible with the metric, in the sense that
\begin{equation}\label{09.05.17-1}\eta_p(a\circ_p b, c)=\eta_p(a, b\circ_p c),\quad \text{for all }a,b,c\in T_pM,
\end{equation}
by axioms (FM2-a),(FM2-b),(FM3-a). Moreover, there exists an open neighborhood $\Omega\subseteq M$ of $p$ and a function $F\colon\Omega\to \mathbb C$ such that 
\begin{align}\label{09.12.18-1}c^\flat=&\nabla^{3}F,\\
\label{09.12.18-2} \eta=&\nabla_e\nabla^{2}F. 
\end{align}This follows from the axiom (FM2-b). Any such a function $F$ will be called \emph{potential} of $M$.

\begin{rem}
The Euler vector field $E$ is an affine vector field, i.e.
\[\nabla^{2}E=0.
\]This follows\footnote{For a generic vector field $X$ on a pseudo-riemannian manifold $(M,g)$, a simple computation (invoking the first Bianchi identities) shows that
\[\nabla_\beta\nabla_\alpha X_\lambda=\sum_\mu R_{\lambda\alpha\beta\mu}X^\mu+\frac{1}{2}\left(\nabla_\beta K_{\alpha\lambda}+\nabla_\alpha K_{\beta\lambda}-\nabla_\lambda K_{\alpha\beta}\right),
\]where
\[K_{\alpha\beta}=(\frak L_Xg)_{\alpha\beta}=\nabla_\alpha X_\beta+\nabla_\beta X_\alpha.
\] If $X$ is Killing conformal, and $\frak L_Xg=\omega g$ for a function $\omega$, then 
\[\nabla_\beta\nabla_\alpha X_\lambda=\sum_\mu R_{\lambda\alpha\beta\mu}X^\mu+\frac{1}{2}\left(g_{\alpha\lambda}\partial_\beta\omega +g_{\beta\lambda}\partial_\alpha\omega -g_{\alpha\beta}\partial_\lambda\omega \right).
\]In our case $R=0$ and $\omega$ is a constant function.} from the axioms (FM1) and (FM4-b).
\end{rem}

\begin{conv}In this paper, we assume that the flat endomorphism $X\mapsto \nabla_XE$ of $TM$ is \emph{diagonalizable}. By introducing $\nabla$-flat coordinates $\bm t=(t^\alpha)_{\alpha=1}^n$ on $M$, w.r.t. which the metric $\eta$ is constant and the connection $\nabla$ coincides with partial derivatives,  we have that
\[E=\sum_{\alpha=1}^n\left((1-q_\alpha)t^\alpha+r_\alpha\right)\frac{\partial}{\partial t^\alpha},\quad q_\alpha,r_\alpha\in\mathbb C.
\]
Following \cite{dubro1, dubro0, dubro2}, we choose flat coordinates $\bm t$ so that $\frac{\partial}{\partial t^1}\equiv e$ and $r_\alpha\neq 0$ only if $q_\alpha=1$. This can always be done, up to an affine change of coordinates.
\end{conv}

\begin{rem}
The associativity of the algebra is equivalent to the following conditions for $F$, called WDVV-equations:
\[\sum_{\gm,\dl=1}^n{\partial_\alpha\partial_\beta\partial_\gamma}F~\eta^{\gamma\delta}{\partial_\delta\partial_\epsilon\partial_\nu}F=\sum_{\gm,\dl=1}^n{\partial_\nu\partial_\beta\partial_\gamma}F~\eta^{\gamma\delta}{\partial_\delta\partial_\epsilon\partial_\alpha}F,
\]while axiom (FM4) is equivalent to
\[\eta_{\alpha\beta}=\partial_1\partial_\alpha\partial_\beta F,\quad \frak L_EF=(3-d)F+Q(\bm t),
\]with $Q(\bm t)$ a quadratic expression in $t_\alpha$'s. Conversely, given a solution of the WDVV equations, satisfying the quasi-homogeneity conditions above, a structure of Frobenius manifold is naturally defined on an open subset of the space of parameters $t^\alpha$'s.
\end{rem}

\begin{defn}
We call \emph{grading operator} of $M$ to be the tensor ${\bm\mu}\in \Gamma(TM\otimes T^*M)$ defined by
\[{\bm\mu}(Y):=\frac{2-d}{2}Y-\nabla_YE,\quad Y\in\Gamma(TM).
\]
In what follows we will also denote by $\bm{\mathcal U}$ the (1,1)-tensor defined by $\circ$-multiplication by the Euler vector field, i.e.
\[{\bm{\mathcal U}}(Y):=E\circ Y,\quad Y\in\Gamma(TM).
\]We denote by $\mu$ and $\mc U$ the matrices of components of the tensors $\bm\mu$, and $\bm{\mc U}$ respectively, w.r.t.\,\,the system $\bm t$ of $\nabla$-flat coordinates.
\end{defn}

\subsection{Semisimple points and bifurcation set} 
\begin{defn}
A point $p\in M$ is \emph{semisimple} if and only if the corresponding Frobenius algebra $(T_pM,*_p,\eta_p, \frac{\partial}{\partial t^1}|_p)$ is without nilpotents. Denote by $M_{ss}$ the open dense subset of $M$ of semisimple points.
\end{defn}

In this paper, only generically semisimple Frobenius manifolds are considered. In other words, we will always assume $M_{ss}\neq \emptyset$.

On $M_{ss}$ there are $n$ well-defined idempotent vector fields $\pi_1,\dots,\pi_n\in\Gamma(TM_{ss})$, satisfying
\beq
\pi_i*\pi_j=\delta_{ij}\pi_i,\quad\eta(\pi_i,\pi_j)=\delta_{ij}\eta(\pi_i,\pi_i),\quad i,j=1,\dots,n.
\eneq

\begin{thm}[\cite{dubronapoli,dubro1,dubro2}]\label{canu}
The idempotent vector fields pairwise commute: $[\pi_i,\pi_j]=0$ for $i,j=1,\dots, n$. Hence, there exist holomorphic local coordinates $(u_1,\dots, u_n)$ on $M_{ss}$ such that $\frac{\partial}{\partial u_i}=\pi_i$ for $i=1,\dots, n$. 
\end{thm}

\begin{defn}
The coordinates $(u_1,\dots, u_n)$ of Theorem \ref{canu} are called \emph{canonical coordinates}. 
\end{defn}

\begin{prop}[\cite{dubro1,dubro2}]
Canonical coordinates are uniquely defined up to ordering and shifts by constants. The eigenvalues of the tensor $\bm{\mathcal U}$ define a system of canonical coordinates in a neighborhood of any semisimple point of $M_{ss}$.
\end{prop}

\begin{defn}
Given a Frobenius manifold $M$, we call \emph{bifurcation set} of $M$ the set $\mc B_M$ of points $p\in M$ at which the spectrum of the operator $\mc U(p)$ is not simple, i.e. $u_i(p)=u_j(p)$ for some $i\neq j$. 

Following the terminology of \cite{CG2,CDG,CDG1}, the points of $\mc B_M$ which are semisimple are called \emph{semisimple coalescing points}. We define the\footnote{The name is taken from singularity theory: for Frobenius structures defined on the universal space of unfoldings of singularities the two notions coincide, see \cite{AGZV,arnold1,singularity1}.} \emph{Maxwell stratum} of $M$ to be the closure of the set of semisimple coalescing points, i.e. $\mc M_M:=\overline{M_{ss}\cap \mc B_M}$. 

The \emph{caustic} of $M$, is the set-theoretic difference $\mc K_M:=\mc B_M\setminus M_{ss}$.
\end{defn}

\begin{lem}
We have $\mc B_M=\mc M_M\cup\mc K_M$.\qed
\end{lem}

\begin{defn}
We call \emph{orthonormalized idempotent frame} a frame $(f_i)_{i=1}^n$ of $TM_{ss}$  defined by
\beq\label{fvects}
f_i:=\eta(\pi_i,\pi_i)^{-\frac{1}{2}}\pi_i,\quad i=1,\dots,n,
\eneq 
for arbitrary choices of signs of the square roots. The $\Psi$-matrix is the matrix $(\Psi_{i\alpha})_{i,\alpha=1}^n$ of change of tangent frames,  defined by
 \beq
 \frac{\partial}{\partial t^\alpha}=\sum_{i=1}^n\Psi_{i\alpha}f_i,\quad \alpha=1,\dots,n.
 \eneq
\end{defn}

\begin{rem}
In the orthonormalized idempotent frame, the operator $\bm{\mathcal U}$ is represented by a diagonal matrix, and the operator $\bm\mu$ by an antisymmetric matrix:
\beq
U:=\operatorname{diag}(u_1,\dots,u_n),\quad\Psi\mathcal U\Psi^{-1}=U,
\eneq
\beq
V:=\Psi\mu\Psi^{-1},\quad V^T+V=0.
\eneq
\end{rem}

\subsection{Extended deformed connection}
Given a Frobenius manifold $M$, let us introduce the extended manifold $\widehat M:=\mathbb C^*\times M$, and let us consider the pull-back $\pi^*TM$ of the tangent bundle of $M$ along the obvious projection $\pi\colon\widehat M\to M$. We will denote the natural lifts on $\widehat M$ of the tensors $\eta,c,e,E,\bm\mu,\bm{\mathcal U}$ by the same symbols. Moreover, we also denote by $\nabla$ the pull-backed Levi-Civita connection: it is the connection on the vector bundle $\pi^*TM$, uniquely defined by the further requirement that
\[\nabla_{\frac{\partial}{\partial z}} Y=0,\quad \text{for all }Y\in\pi^{-1}\mathscr T_M,
\]where $z$ denotes the natural coordinate on $\mathbb C^*$, and $\mathscr T_M$ denotes the tangent sheaf of $M$. We are going now to define a second connection $\widehat\nabla$ on $\pi^*TM$ which is a deformation of $\nabla$.

\begin{defn}
We define the \emph{extended deformed connection} $\widehat \nabla$ as the connection on $\pi^*TM$ given by
\begin{align*}
\widehat\nabla_XY=&\nabla_XY+z X\circ Y,\\
\widehat\nabla_{\frac{\partial}{\partial z}}Y=&\nabla_{\frac{\partial}{\partial z}}Y+\bm{\mathcal U}(Y)-\frac{1}{z}\bm\mu(Y),
\end{align*}
for all $X,Y\in\Gamma(\pi^*TM)$.
\end{defn}

\begin{thm}[\cite{dubro2}]
The extended deformed connection $\widehat\nabla$ if flat. More precisely, its flatness is equivalent to the totality of the following conditions:
\begin{enumerate}
\item $\nabla c^\flat\in \Gamma( \odot^4T^*M)$,
\item the product on each tangent space of $M$ is associative,
\item $\nabla^2E=0$,
\item $\frak L_Ec=c$.\qed
\end{enumerate}
\end{thm}

The connection $\widehat\nabla$ induces a flat connection on $\pi^*T^*M$, denoted by the same symbol.

\subsection{Cyclic stratum, and cyclic (co)frame}

\begin{defn}Given a Frobenius manifold $M$, we define infinitely many sections $e_j\in\Gamma(\pi^*TM)$ as
\[e_k:=\widehat\nabla_{\frac{\partial}{\partial z}}^ke,\quad k\in\mathbb N.
\]
We will call the \emph{cyclic stratum} ${\widehat M}^{\rm cyc}$ to be the maximal open subset $U$ of $\widehat M$ such that the bundle $\pi^*TM|_U$ is trivial and the collection of sections $(e_k|_{U})_{k=0}^{n-1}$ defines a basis of each fiber. On ${\widehat M}^{\rm cyc}$ we will also introduce the dual coframe $(\omega_j)_{j=0}^{n-1}$, by imposing
\beq
\label{eom}
\langle\omega_j,e_k\rangle=\delta_{jk}.
\eneq
The frame $(e_k)_{k=0}^{n-1}$ will be called \emph{cyclic frame}, and its dual $(\om_j)_{j=0}^{n-1}$ \emph{cyclic coframe}.
\end{defn}

\begin{defn}
Define the matrix-valued function $\La=(\La_{i\al}(z,p))$, holomorphic on ${\widehat M}^{\rm cyc}$, by the equation
\beq\label{Lam}
\frac{\der}{\der t^\al}=\sum_{i=0}^{n-1}\La_{i\al}e_i,\quad \al=1,\dots,n.
\eneq
\end{defn}

\begin{rem}
The $\La$-matrix should be thought as an analogue of the $\Psi$-matrix. The former relates the flat coordinate frame $(\frac{\der}{\der t^\al})_{\al=1}^n$ to the cyclic frame $(e_i)_{i=0}^{n-1}$. The latter relates the flat coordinate frame $(\frac{\der}{\der t^\al})_{\al=1}^n$ to the normalized idempotent frame $(f_i)_{i=1}^n$. 
\end{rem}

\begin{lem}\label{derom}
For $j=1,\dots,n-1$, we have $\widehat\nabla_{\frac{\partial}{\partial z}}\om_j=-\om_{j-1}$.
\end{lem}
\proof
From \eqref{eom}, for any $k=0,\dots,n-2$, we have 
\begin{align*}\langle\widehat\nabla_{\frac{\partial}{\partial z}}\om_j,e_k\rangle+\langle\om_j,e_{k+1}\rangle=0\quad\Longrightarrow &\quad\langle\widehat\nabla_{\frac{\partial}{\partial z}}\om_j,e_k\rangle=-\delta_{j,k+1}\\
\puqed
\quad\Longrightarrow &\quad \widehat\nabla_{\frac{\partial}{\partial z}}\om_j=-\om_{j-1}. \qedhere
\poqed
\end{align*}

\begin{prop}\label{structcycvec}
The vector fields $e_k$, with $k\in\mathbb N$, have the following form
\[e_k=\sum_{j=0}^k\frac{1}{z^j}p^k_j(E),
\]where the vector fields $p^k_j(E)$ do not depend on $z$ and satisfy the difference equations
\begin{align*}
p^{k+1}_0(E)=&E\circ p^k_0(E),\\
p^{k+1}_j(E)=&E\circ p^k_j(E)-\bm\mu(p^k_{j-1}(E))+(1-j)p^k_{j-1}(E),\quad j=1,\dots, k\\
p^{k+1}_{k+1}(E)=&-\mu(p^k_{k}(E))-kp^k_{k}(E),
\end{align*}
with the only initial datum $p^{0}_j(E)=\delta_{0j}\cdot e$.\qed
\end{prop}

\subsection{Properties of the function $\det \La$}
The holomorphic function $\det \La\colon \widehat M^{\rm cyc}\to \C^*$ extends meromorphically to a function on $\P^1\times M$.
\begin{thm}\label{strLa}
The function $\det\La$ is a meromorphic function on $\P^1\times M$  of the form
\[\det\La(z,p)=\frac{z^{\binom{n-1}{2}}}{z^{\binom{n-1}{2}}A_0(p)+\dots+A_{\binom{n-1}{2}}(p)},
\]
where $A_0,\dots, A_{\binom{n-1}{2}}$ are holomorphic functions on $M$. Moreover, if $n>2$ and if the eigenvalues of the grading operator $\bm\mu$ are not pairwise distinct, then the function $A_{\binom{n-1}{2}}$ is identically zero. 
\end{thm}
We need a preliminary result. 

\begin{lem}\label{polwk}
For $k\in\{0,\dots, n-1\}$, the polyvector field $e_0\wedge\dots\wedge e_k\in\Gamma(\bigwedge^{k+1}\pi^*TM)$ admits a pole at $\{0\}\times M$ of order at most $\binom{k}{2}$. More precisely, we have
\[e_0\wedge\dots\wedge e_k=w_0+\frac{1}{z}w_1+\dots+\frac{1}{z^{\binom{k}{2}}}w_{\binom{k}{2}},\qquad w_j\in\Gm\left(\bigwedge\nolimits^{k+1}\pi^*TM\right),
\]with
\[w_{\binom{k}{2}}=(-1)^{\binom{k}{2}}\,e\wedge E\wedge \bm\mu(E)\wedge\bm\mu^2(E)\wedge\dots\wedge \bm\mu^{k-1}(E).
\]
\end{lem}
\proof By induction on $k$. For the base cases $k=0$ and $k=1$, we have $e_0=e$ and $e_0\wedge e_1=e\wedge E$, respectively. So, for $k=0,1$ the claim holds true.

Assume that $e_0\wedge\dots\wedge e_{k-1}$ is of the form
\[e_0\wedge\dots\wedge e_{k-1}=w_0+\frac{1}{z}w_1+\dots+\frac{1}{z^{\binom{k-1}{2}}}w_{\binom{k-1}{2}},
\]with 
\[w_{\binom{k-1}{2}}=(-1)^{\binom{k-1}{2}}\,e\wedge E\wedge \bm\mu(E)\wedge\bm\mu^2(E)\wedge\dots\wedge \bm\mu^{k-2}(E).
\]We have 
\[e_0\wedge\dots\wedge e_k=\left(\sum_{j=0}^{\binom{k-1}{2}}z^{-j}w_j\right)\wedge\left(\sum_{\ell=0}^kz^{-\ell}p^k_\ell(E)\right).
\]We claim that the coefficient $w_{\binom{k-1}{2}}\wedge p^k_k(E)$ of $z^{-\binom{k-1}{2}-k}$ vanishes. Indeed, $p^k_k(E)$ is proportional to $e$: we have
\[p^k_k(E)=\frac{d}{2}\left(\frac{d}{2}-1\right)\dots\left(\frac{d}{2}-k+1\right)e,\qquad k\geq 0,
\]as it can be easily seen by induction (the key property is $\mu(e)=-\frac{d}{2}e$, together with the last difference equation of Proposition \ref{structcycvec}). Consequently, we have $w_{\binom{k-1}{2}}\wedge p^k_k(E)=c\cdot (e\wedge\dots\wedge e)=0$.

Hence, the (possibly non-vanishing) most polar term of $e_0\wedge\dots \wedge e_k$ equals
\begin{align*}z^{-\binom{k-1}{2}-k+1}\cdot w_{\binom{k-1}{2}}\wedge p^k_{k-1}(E)&=z^{-\binom{k}{2}}\cdot w_{\binom{k-1}{2}}\wedge ((-1)^{k-1}\mu^{k-1}(E))\\
&=z^{-\binom{k}{2}}(-1)^{\binom{k}{2}}\,e\wedge E\wedge \bm\mu(E)\wedge\dots\wedge \bm\mu^{k-1}(E).
\end{align*}For the first equality we have used the difference equation for $p^k_{k-1}(E)$ of Proposition \ref{structcycvec}. 
\endproof

\proof[Proof of Theorem \ref{strLa}]
The polyvector field $e_0\wedge\dots\wedge e_{n-1}$ has the form
\beq\label{wp}
e_0\wedge\dots\wedge e_{n-1}=w_0(p)+\frac{1}{z}w_1(p)+\dots\frac{1}{z^{\binom{n-1}{2}}}w_{\binom{n-1}{2}}(p),
\eneq
where $w_0,w_1,\dots,w_{\binom{n-1}{2}}$ are holomorphic $n$-vector fields on $M$, by Lemma \ref{polwk}. Introduce holomorphic functions $A_0(p),\dots, A_{\binom{n-1}{2}}(p)$, such that $$w_j(p)=A_j(p)\cdot \frac{\der}{\der t^1}\wedge\dots\wedge\frac{\der}{\der t^n}.$$
From the identity 
\[\frac{\der}{\der t^1}\wedge\dots\wedge\frac{\der}{\der t^n}=\det\La\cdot e_0\wedge\dots\wedge e_{n-1},
\]we deduce
\[
1=\det\La(z, p)\left(A_0(p)+\frac{1}{z}A_1(p)+\dots\frac{1}{z^{\binom{n-1}{2}}}A_{\binom{n-1}{2}}(p)\right).
\] The last statement on $A_{\binom{n-1}{2}}$ follows from the explicit formula for $w_{\binom{n-1}{2}}$ given in Lemma \ref{polwk}.
\endproof

\begin{thm}\label{asydetLa}We have
\[A_0(p)=\frac{\prod_{i<j}(u_j(p)-u_i(p))}{{\rm Jac}(p)},\quad {\rm Jac}(p):=\det\left.\left(\frac{\der u_i}{\der t^\al}\right)\right|_p.
\]
\end{thm}

\proof The polyvector field $w_0$ in equation \eqref{wp} is
\[w_0=\bigwedge_{j=0}^{n-1}p^{j}_0(E).
\]By Proposition \ref{structcycvec}, we have that $$p^{j}_0(E)=E^{\circ j},\quad j\in\mathbb N,$$
and using the idempotent vielbein $(\frac{\partial}{\partial u_i})_{i=1}^n$, we can write $w_0$ as follows
\begin{align*}
w_0=&\begin{vmatrix}
1&\dots&1\\
u_1&\dots&u_n\\
u_1^2&\dots&u_n^2\\
&\vdots&\\
u_1^{n-1}&\dots&u_n^{n-1}
\end{vmatrix}\frac{\partial }{\partial u_1}\wedge\dots\wedge\frac{\partial }{\partial u_n}
\label{bz}
=\left(\prod_{i<j}(u_j-u_i)\right)\frac{\partial }{\partial u_1}\wedge\dots\wedge\frac{\partial }{\partial u_n}\\
\puqed
=&\left(\prod_{i<j}(u_j-u_i)\right)\cdot \frac{1}{{\rm Jac}}\cdot  \frac{\der}{\der t^1}\wedge\dots\wedge\frac{\der}{\der t^n}.\qedhere
\poqed
\end{align*}

\begin{rem}
We also have
\[
\frac{\der}{\der t^1}\wedge\dots\wedge\frac{\der}{\der t^n}=\det\Psi f_1\wedge\dots\wedge f_n=\frac{\det \Psi}{\prod_{i=1}^n\eta(\frac{\der}{\der u_i},\frac{\der}{\der u_i})^\frac{1}{2}}\frac{\der}{\der u_1}\wedge\dots\wedge\frac{\der}{\der u_n},
\]
so that 
\[{\rm Jac}(p)=\left.\frac{\det \Psi}{\prod_{i=1}^n\eta(\frac{\der}{\der u_i},\frac{\der}{\der u_i})^\frac{1}{2}}\right|_p=\left.\frac{(\det \eta)^\frac{1}{2}}{\prod_{i=1}^n\eta(\frac{\der}{\der u_i},\frac{\der}{\der u_i})^\frac{1}{2}}\right|_p.
\]
The last equality follows from $\Psi^T\Psi=\eta$.
\end{rem}

\subsection{Geometry of the complement of the cyclic stratum in $\P^1\times M$}\label{strat} Consider the tuple of functions $(A_0,\dots, A_{\binom{n-1}{2}})$, and extend it to the sequence $(A_k)_{k\in\N}$ by setting $A_k=0$ for $k>\binom{n-1}{2}$. Set $$\bar n:=\min\{j\in\N\colon A_{h}(p)=0\quad\forall p\in M,\,\forall h>j\}.$$
We necessarily have $0\leq \bar n\leq \binom{n-1}{2}$. By Theorem \ref{strLa}, we have $\bar n<\binom{n-1}{2}$ if $\bm \mu$ has not simple spectrum. 
The function $\det\La$ takes the form
\[\det \La=\frac{z^{\bar n}}{z^{\bar n} A_0(p)+z^{\bar n-1}A_1(p)\dots+A_{\bar n}(p)}.
\]
Define the subsets $\mc P_\La, M_0,M_\infty\subseteq \P^1\times M$ and $\mc A_\La,\mc I_\La^\infty,\mc I_\La^0\subseteq M$ by
\[\mc P_\La:=\left\{(z,p)\in \widehat M\colon\quad  z^{\bar n}A_0(p)+\dots+A_{\bar n}(p)=0\right\},
\]
\[M_0:=\{0\}\times M,\quad M_\infty:=\{\infty\}\times M,
\]
\[\mc A_\La:=\left\{p\in M\colon\quad A_0(p)=\dots=A_{\bar n}(p)=0\right\},
\]
\[\mc I_\La^\infty:=\left\{p\in M\colon \quad A_0(p)=0\right\},
\]
\[\mc I_\La^{0}:=\left\{p\in M\colon\quad A_{\bar n}(p)=0\right\}.
\]

\begin{lem}\label{ala}
We have the obvious inclusions 
\[\puqed
\C^*\times \mc A_\La\subseteq \mc P_\La,\quad \mc A_\La\subseteq \mc I^0_\La\cap \mc I^\infty_\La. \qedhere
\poqed
\]
\end{lem}

\begin{table}
\[\begin{array}{|c|c|}
\hline
&\\
\text{Poles of } \det\La& \mc P_\La\cupdot\ (\{\infty\}\times \mc I^\infty_\La)\\
&\\
\hline
&\\
\text{Zeros of } \det\La& M_0\setminus (\{0\}\times \mc I^0_\La)\\
&\\
\hline
&\\
\text{Indeterminacy locus of } \det\La& \{0\}\times \mc I^0_\La\\
&\\
\hline
\end{array}
\]
\caption{In this table, we summarize the location of poles, zeros and indeterminacy locus for the meromorphic function $\det\La$ on $\P^1\times M$.}
\end{table}

The set $\mc P_\La$ is an analytic subspace of $\P^1\times M$ of codimension 1 along which the function $\det\La$ admits a pole. The function $\det \La$ admits poles along a further analytic subspace, namely $\{\infty\}\times \mc I^\infty_\La$.

The set $\mc P_\La$ is the complement $\widehat M\setminus \widehat M^{\rm cyc}$ of the cyclic stratum.
The complement of $\widehat M^{\rm cyc}$ in $\P^1\times M$ is the disjoint union 
\[\mc P_\La\ \cupdot\ M_0\ \cupdot\ M_\infty. 
\] The geometry of $\mc P_\La$ is rather complicated: in general it admits several irreducible components. For example, $\mc A_\La$ itself does, and consequently also $\C^*\times \mc A_\La$.  The projection $\pi\colon \widehat M\to M$, if restricted to $\mc P_\La\setminus (\C^*\times \mc A_\La)$, defines a ramified covering of degree $\bar n$. 

The set $\{0\}\times\mc I^0_\La$ is an analytic subspace of $\P^1\times M$ of codimension 2 and it is the \emph{indeterminacy locus} of the function $\det\La$.

Each of the sets $\mc I^\infty_\La, \mc I^0_\La, \mc A_\La$ seems to be strictly related to other distinguished subsets of the Frobenius manifold $M$, namely its bifurcation set $\mc B_M$, and its two  components, the {Maxwell stratum} $\mc M_M$ and the {caustic} $\mc K_M$. We limit to the following observation.

\begin{thm}\label{ib}
We have $\mc I^\infty_\La\subseteq \mc B_M$.
\end{thm}
\proof Let $p\notin \mc B_M$. 
On the complement of $\mc B_M$, the eigenvalues $(u_1,\dots, u_n)$ define a holomorphic system of coordinates. Hence, ${\rm Jac}(p)\neq 0$. Moreover, by definition we have $\prod_{i<j}(u_j(p)-u_i(p))\neq 0$. Hence, $p\notin \mc I^\infty_\La$ by Theorem \ref{asydetLa}.
\endproof

\begin{figure}[ht!]
\centering
\def\svgscale{.7}
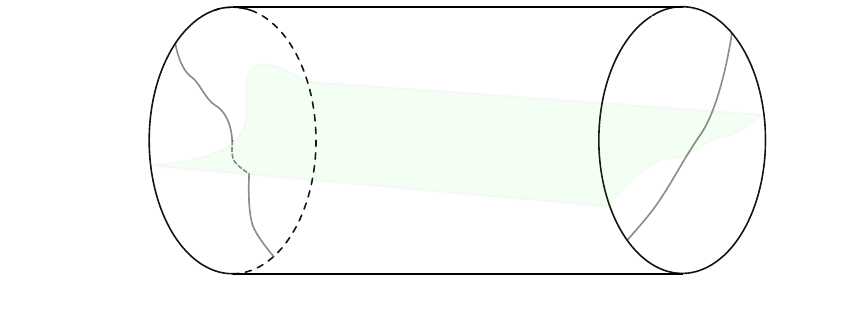
\caption{Configuration of the sets $\mc P_\La$, $\{\infty\}\times \mc I^\infty_\La$, and $\{0\}\times\mc I^0_\La$ in $\P^1\times M$. }
\label{cyclic}
\end{figure}

In order to obtain more precise results on contingent relations between the sets $\mc I^\infty_\La, \mc I^0_\La, \mc A_\La$ and $\mc B_M, \mc M_M,\mc K_M$ a more detailed study of the polyvector fields $p^k_j(E)$ of Proposition \ref{structcycvec} is needed. 
We plan to address this problem in a future project. We conclude this section with three low dimensional examples.

\begin{example}
For 2-dimensional Frobenius manifolds, we have $\mc I^\infty_\La= \mc B_M$. In this case, indeed, we have
\[e_0=e,\quad e_1= E+\frac{d}{2z}e\quad\Rightarrow\quad e_0\wedge e_1=e\wedge E.
\]The bivector $e\wedge E$ vanishes if and only if $u_1=u_2$.
\end{example}

\begin{example}
Consider the $A_3$-Frobenius manifold, that is the space $M\cong \C^3$ of polynomials $f(x,\bm a)=x^4+a_2x^2+a_1x+a_0$, where $\bm a=(a_0,a_1,a_2)\in \C^3$ are natural coordinate. Fix $\bm a_o\in M$, and define the Kodaira-Spencer isomorphism $\kappa\colon T_{\bm a_o}M\to \C[x]/\langle \der_x f(x,\bm a_o)\rangle$, by identifying $\der_{a_i}$ with the class of $\der_{a_i} f(x,\bm a_o)$. This allows to pull-back the product of the Jacobi-Milnor algebra $\C[x]/\langle \der_x f(x,\bm a_o)\rangle$ on $T_{\bm a_o}M$. Consider the Grothendieck residue metric
\[\eta_{\bm a}\left(\frac{\der}{\der a_i},\frac{\der}{\der a_j}\right):=\frac{1}{2\pi i}\int_{\Gm_{\bm a}}\left.\frac{\frac{\der f}{\der a_i}\frac{\der f}{\der a_j}}{\frac{\der f}{\der x}}\right|_{(u,\bm a)}du,
\]
where $\Gm_{\bm a}$ is a circle, positively oriented, bounding a disc containing all the roots of $\frac{\der f}{\der x}(u,\bm a)$. One can show that the coordinates $\bm t=(t_1,t_2,t_3)$ given by
\[t_1=a_0-\frac{1}{8}a_2^2,\quad t_2=a_1,\quad t_3=a_2,
\]are flat for the metric $\eta$. In $\bm t$-coordinates, the Euler vector field is given by
\[E=t_1\frac{\der}{\der t_1}+\frac{3t_2}{4}\frac{\der}{\der t_2}+\frac{t_3}{2}\frac{\der}{\der t_3}.
\]
The Maxwell stratum is the set $\{t_2=0\}$, and the caustic is the set $\{8t_3^3+27t_2^2=0\}$.
We have the following formulae for the $\La$-matrix and for $\det \La$:
\begin{align*}
&\La(z,\bm t)=\\
&\left(
\begin{array}{ccc}
 1 & \frac{z^2 t_3^5-21 z^2 t_2^2 t_3^2-64 z^2 t_1^2 t_3-12 t_3-18 z t_2^2-72 z^2 t_1 t_2^2}{2 z t_2 \left(8 z t_3^3-6 t_3+27 z t_2^2\right)} & \frac{-3 z^2 t_3^4-16 z t_3^2-64 z^2 t_1 t_3^2+63 z^2 t_2^2 t_3+192 z^2 t_1^2+48 z t_1+48}{4 z \left(8 z t_3^3-6 t_3+27 z t_2^2\right)} \\
 0 & \frac{4 \left(9 z t_2^2+16 z t_1 t_3\right)}{t_2 \left(8 z t_3^3-6 t_3+27 z t_2^2\right)} & -\frac{4 \left(-4 z t_3^2+24 z t_1+3\right)}{8 z t_3^3-6 t_3+27 z t_2^2} \\
 0 & -\frac{32 z t_3}{t_2 \left(8 z t_3^3-6 t_3+27 z t_2^2\right)} & \frac{48 z}{8 z t_3^3-6 t_3+27 z t_2^2} \\
\end{array}
\right),
\end{align*}
\[
\det\La(z,\bm t)=\frac{64 z}{(8 t_2t_3^3 +27 t_2^3) z-6 t_2t_3}.
\]
We have 
\begin{itemize}
\item $\mc I^\infty_\La=\mc B_M$, 
\item $\mc I^0_\La=\mc M_M\cup\{t_3=0\}$,
\item $\mc A_\La=\mc M_M$.
\end{itemize}
\end{example}

\begin{example}
The $A_2\times A_2$--Frobenius manifold is the Frobenius structure $M$ on $\C^4$, with flat coordinates $(\bm t,\bm s)=(t_0,t_1,s_0,s_1)$, defined by the WDVV--potential 
\[F(\bm t,\bm s)=\frac{1}{2}\left(t_0^2t_1+s_0^2s_1\right)+\frac{1}{24}(t_1^4+s_1^4).
\]In these coordinates, the unit vector field is $e=\frac{\der}{\der t_0}+\frac{\der}{\der s_0}$, and the flat metric $\eta$ has components 
\[\eta=\begin{pmatrix}
0&1&0&0\\
1&0&0&0\\
0&0&0&1\\
0&0&1&0
\end{pmatrix}.
\] The Euler field equals
\[E=t_0\frac{\der}{\der t_0}+\frac{2}{3}t_1\frac{\der}{\der t_1}+s_0\frac{\der}{\der s_0}+\frac{2}{3}s_1\frac{\der}{\der s_1}.
\]The bifurcation diagram $\mc B_{M}$ equals
\[\mc B_{M}=\mc M_{M}\cup \mc K_{M},
\]where the Maxwell stratum is
\[\mc M_{M}=\left\{-8 t_1^3 \left(9 (s_0-t_0)^2+4 s_1^3\right)+\left(4 s_1^3-9 (s_0-t_0)^2\right)^2+16 t_1^6=0\right\},
\]and the caustic is
\[\mc K_{M}=\{t_1=0\}\cup\{s_1=0\}.
\]After some computations, one finds that 
\begin{multline*}\det \La(z,\bm t,\bm s)=\\ \frac{729 z^2}{4 s_1 t_1 \left(z^2 \left(-8 t_1^3 \left(9 (s_0-t_0)^2+4 s_1^3\right)+\left(4 s_1^3-9 (s_0-t_0)^2\right)^2+16 t_1^6\right)+45 (s_0-t_0)^2\right)}.
\end{multline*}
We have
\begin{itemize}
\item $\mc I^\infty_\La=\mc B_{M}$,
\item $\mc I^0_\La=\mc K_{M}\cup\{s_0=t_0\}$,
\item $\mc A_\La=\mc K_{M}\cup \{s_0=t_0,\,s_1^3=t_1^3\}$.
\end{itemize}
\end{example}

\subsection{Master differential equation and master functions}\label{secmas}

Let $\xi\in\Gamma(\pi^*T^*M)$ be a $\edc$-flat section. Consider the corresponding vector field $\zeta\in\Gamma(\pi^*TM)$ via musical isomorphism, i.e. such that $\xi(v)=\eta(\zeta,v)$ for all $v\in\Gamma(\pi^*TM)$.

The vector field $\zeta$ satisfies the following system\footnote{We consider the joint system \eqref{eq1}, \eqref{qde} in matrix notations ($\zeta$ is a column vector whose entries are the components $\zeta^\alpha(z, \bm t)$ w.r.t. $\frac{\partial}{\partial t^\alpha}$). Bases of solutions are arranged in invertible $n\times n$-matrices, called \emph{fundamental systems of solutions}.} of equations
 \begin{align}
 \label{eq1}
 \frac{\partial}{\partial t^\alpha}\zeta&=z\mathcal C_\alpha\zeta,\quad \alpha=1,\dots, n,\\
\label{qde}
  \frac{\partial}{\partial z}\zeta&=\left(\mathcal U+\frac{1}{z}\mu\right)\zeta.
 \end{align}
Here $\mathcal C_\alpha$ is the $(1,1)$-tensor defined by $(\mathcal C_\alpha)^\beta_\gamma:=c^\beta_{\alpha\gamma}$. 
\vskip 3mm
Multiply by $\eta$ (on the left) the l.h.s. and r.h.s. of \eqref{eq1}, \eqref{qde}: we obtain the equivalent system of differential equations
 \begin{align}
 \label{eq1du}
 \frac{\partial}{\partial t^\alpha}\xi&=z\mathcal C^T_\alpha\xi,\quad \alpha=1,\dots, n,\\
\label{qdedu}
  \frac{\partial}{\partial z}\xi&=\left(\mathcal U^T-\frac{1}{z}\mu\right)\xi,
 \end{align}
where $\xi$ is a column vector whose entries are the components $\xi_\al(z,\bm t)$ w.r.t. $dt^\al$. At points $(z,p)\in \widehat M^{\rm cyc}$, let us introduce the column vector $\bar\xi$ by
\beq
\label{barxixi}
\bar\xi=(\La^{-1})^T\xi,
\eneq
where $\La$ is defined as in \eqref{Lam}. The entries of $\bar\xi$ are the components $\bar\xi_j$ w.r.t. the cyclic coframe $\om_j$. The vector $\bar \xi$ satisfies the system
\begin{align}
\frac{\der\bar\xi}{\der t^\al}&= \left(z(\La^{-1})^T\mc C_\al \La^T+\frac{\der(\La^{-1})^T}{\der t^\al}\La^T\right)\bar\xi,\\
\label{qdediag}
\derz{\bar\xi}&=\left((\La^{-1})^T\mc U^T \La^T-\frac{1}{z}(\La^{-1})^T\mu\La^T+\frac{\der(\La^{-1})^T}{\der t^\al}\La^T\right)\bar\xi.
\end{align}

\begin{prop}
Let $\xi\in\Gamma(\pi^*T^*M)$ be a $\edc$-flat section, and let $(\bar\xi_j(z,p))_{j=0}^{n-1}$ be its components w.r.t. the cyclic co-frame, i.e. $\xi=\sum_j\bar\xi_j\om_j$. We have
\beq
\frac{\partial \bar\xi_j}{\partial z}=\bar\xi_{j+1},\quad j=0,\dots, n-2.
\eneq
\end{prop}

\proof
We have
\begin{align*}
0=\edc_{\derz{}}\xi&=\sum_j\derz{\bar\xi_j}\om_j+\sum_j\bar\xi_j\edc_{\derz{}}\om_j\\
&=\sum_j\derz{\bar\xi_j}\om_j-\sum_j\bar\xi_j\om_{j-1},
\end{align*}
by Lemma \ref{derom}. The claim follows. 
\endproof

\begin{cor}
The system of differential equations \eqref{qdediag} is the companion system of a \emph{scalar} differential equation in $\bar\xi_0$. \qed
\end{cor}

\begin{rem}
Notice that $\xi_1=\bar\xi_0$. Indeed, we have $e_0=e=\frac{\der}{\der t^1}$, so that $\La_{i1}=\dl_{i1}$. The claim then follows from \eqref{barxixi}. 
\end{rem}

\begin{thm}\label{thmde}
Consider the system of differential equations \eqref{qdedu}, specialized at a point $p\in M\setminus \mc A_\La$. The system can be reduced to a single scalar ordinary differential equation of order $n$ in the unknown function $\xi_1$. The scalar differential equation admits at most $\binom{n-1}{2}$ apparent singularities.
\end{thm}
\proof
If $p\in M\setminus \mc A_\La$, then there exist $\bar n$ complex numbers $z_1,\dots, z_{\bar n}$, not necessarily distinct, such that $(z_i,p)\notin \widehat{M}^{\rm cyc}$. The numbers $z_i$ are the zeroes of the denominator of the function $\det\La(z,p)$.
\endproof

The scalar differential equation to which the system \eqref{qdedu} can be reduced will be called the \emph{master differential equation} of $M$.

\begin{defn}
Fix $p\in M$. Consider the system of differential equations \eqref{qdedu} specialized at $p$, and set $\mc X_p$ be the $\C$-vector space of its solutions. Let $\nu_p\colon \mc X_p\to\mc O(\widetilde{\C^*})$ be the morphism defined by
\[\xi\mapsto \Phi_\xi(z),\quad \Phi_\xi(z):=z^{-\frac{d}{2}}\langle\xi (z,p),e(p)\rangle,
\]where $d$ is the charge of the Frobenius manifold. Set $\mc S_p(M):={\rm im} (\nu_p)$. Elements of $\mc S_p(M)$ will be called \emph{master functions} at $p$.
\end{defn}

\begin{thm}\label{niso}
At points $p\in M\setminus\mc A_\La$ the morphism $\nu_p$ is injective.
\end{thm}
\proof
Given $\Phi_\xi\in\mc S_p(M)$, the function $\xi_1(z)=z^{\frac{d}{2}}\Phi_\xi(z)$ is a solution of the master differential equation at $p$. By Theorem \ref{thmde}, the solution $\xi(z)$ can be reconstructed from the component $\xi_1(z)$ only.
\endproof

\section{Gromov-Witten theory}\label{sec3}

\subsection{Notations and conventions}\label{notX} Let $X$ be a smooth projective variety over $\mathbb C$. In order not to introduce superstructures, in what follows we assume that $H^{\rm odd}(X,\mathbb C)=0$. Denote by $b_k(X)$ the $k$-th Betti number of $X$.
\vskip 3mm
Attached to $X$ there is an infinite dimensional $\mathbb C$-vector space $\mathcal P_X$, called the \emph{big phase space}, defined as the infinite product of countable many copies of the classical cohomology space of $X$, that is
\[\mathcal P_X:=\prod_{n\in\mathbb N}H^\bullet(X,\mathbb C).
\]
Let us fix a homogeneous basis $(T_0,\dots, T_N)$ of $H^\bullet(X, \mathbb C)$ such that
\begin{itemize}
\item $T_0=1$, 
\item and $T_1,\dots, T_r$ is a nef integral basis of $H^2(X,\Z)$.
\end{itemize}
In particular, we have $b_2(X)=r$. Set $\bm t=(t^0,\dots, t^N)$ the dual coordinates of $H^\bullet(X,\C)$.
\vskip3mm
Denote by $(\tau_pT_0,\dots, \tau_pT_N)$ the corresponding basis of the $p$-th copy of $H^\bullet(X,\mathbb C)$ in $\mathcal P_X$. The element $\tau_pT_\alpha$ will be called a \emph{descendant} of $T_\alpha$ with level $p$. The coordinate of a point ${\bsy\gamma}\in\mathcal P_X$ w.r.t. the basis $(\tau_pT_\alpha)_{\alpha,p}$ will be denoted by $\bm t^\bullet=(t^{\alpha,p})_{\alpha,p}$.  
Instead of denoting by ${\bsy\gamma}=(t^{\alpha,p}\tau_{p}T_\alpha)_{\alpha,p}$ a generic element of $\mathcal P_X$ we will usually write this as a formal series
\[{\bsy\gamma}=\sum_{\al=1}^m\sum_{p=0}^\infty t^{\alpha,p}\tau_pT_\alpha.
\] 
We identify $H^\bullet(X,\C)$ with the $0$-th factor of $\mc P_X$, called the \emph{small phase space}. This allow us to identify $t^\alpha\equiv t^{\alpha,0}$ for $\al=0,\dots, N$.

We denote by $\eta\colon H^\bullet(X,\mathbb C)\times H^\bullet(X,\mathbb C)\to H^\bullet(X,\mathbb C)$ the Poincaré pairing defined by
\[\eta(u,v):=\int_Xu\cup v,
\]and we set $\eta_{\alpha\beta}:=\eta(T_\alpha, T_\beta)$ for $\alpha,\beta=0,\dots,N$. The numbers $\eta_{\al\bt}$ will be collected in the Gram\footnote{We denote the metric tensor and its Gram matrix by the same symbol $\eta$. This is a standard abuse of notation.} matrix $\eta=(\eta_{\al\bt})_{\al,\bt=0}^N$, with inverse matrix $\eta^{-1}=\left(\eta^{\al\bt}\right)_{\al,\bt=0}^N$. We also introduce the dual basis $(T^0,\dots, T^N)$ of $H^\bullet(X,\mathbb C)$, by setting
\[T^\al:=\sum_{\la=0}^NT_\la\eta^{\la\al},\quad \al=0,\dots, N.
\]
Define the \emph{Novikov ring} $\La_X$ as the ring of formal sums
\[\sum_{\bt\in H_2(X,\Z)}a_\bt {\bf Q}^\bt,\quad a_\bt\in\Q,
\]such that
\[{\rm card}\left\{\bt\colon a_\bt\neq 0\text{ and }\int_\bt\om<C\right\}<\infty,\quad \text{for any }C\in\R,
\]where $\om$ is the K\"ahler form of $X$.

\subsection{Descendant Gromov-Witten invariants}
For any given $g,n\in\mathbb N$ and $\beta\in H_2(X,\mathbb Z)$, let us denote by $\overline{\mathcal M}_{g,n}(X,\beta)$ the Kontsevich-Manin moduli stack of genus $g$, $n$-pointed stable maps of degree $\beta$ with target $X$: it parametrizes isomorphism classes of pairs $((C,\bsy x),f)$ where 
\begin{itemize}
\item $C$ is a genus $g$ nodal connected projective curve,
\item $\bsy x=(x_1,\dots, x_n)$ is an $n$-tuple of pairwise distinct points of the smooth locus of $C$,
\item $f\colon C\to X$ is a morphism with $f_*[C]=\beta$,
\item a morphism between two pairs $((C,\bsy x), f),((C',\bsy x'), f')$ is a morphism $\sigma\colon C\to C'$ such that $\sigma(x_i)=x_i'$ for all $i$, and making commutative the diagram 
\[\xymatrix{C\ar[rr]^{\sigma}\ar[dr]_{f}&&C'\ar[dl]^{f'}\\&X&}
\]
\item the group of automorphisms of $((C,\bsy x),f)$ is finite.
\end{itemize}
The moduli space $\overline{\mathcal M}_{g,n}(X,\beta)$ is a proper Deligne-Mumford stack of virtual dimension
\[{\rm vir\ dim}_{\mathbb C}\overline{\mathcal M}_{g,n}(X,\beta):=(1-g)(\dim_{\mathbb C}X-3)+\int_\beta c_1(X)+n.
\]Let us denote by $\mathcal L_i$, with $i=1,\dots, n$, the $i$-th \emph{tautological line bundle} on $\overline{\mathcal M}_{g,n}(X,\beta)$ whose fiber at the point $[((C,\bsy x),f)]\in\overline{\mathcal M}_{g,n}(X,\beta)$ is the cotangent space $T^*_{x_i}C$. Set $\psi_j:=c_1(\mc L_j)$ for $j=1,\dots,n$.

We have naturally defined \emph{evaluation morphisms}$${\rm ev}_i\colon \overline{\mathcal M}_{g,n}(X,\beta)\to X,\quad [((C,\bsy x), f)]\mapsto f(x_i)$$ for $i=1,\dots,n$.

\begin{defn}Let $d_1,\dots, d_n$ be non-negative integers. The \emph{genus $g$ descendant Gromov-Witten invariants} (or \emph{genus $g$ gravitational correlators}) are the rational numbers defined by the integrals
\[\langle\tau_{d_1}\alpha_1, \dots, \tau_{d_n}\alpha_n\rangle_{g,n,\beta}^{X}:=\int_{[\overline{\mathcal M}_{g,n}(X,\beta)]^{\rm virt}}\prod_{j=1}^n\psi_j^{d_j}\cup{\rm ev}^*_j(\alpha_j),
\]where $\alpha_1,\dots,\alpha_n\in H^\bullet(X,\mathbb C)$, and the class $$[\overline{\mathcal M}_{g,n}(X,\beta)]^{\rm virt}\in {\rm CH}_{D}(\overline{\mathcal M}_{g,n}(X,\beta)),\quad D={\rm vir}\dim_{\mathbb C}\overline{\mathcal M}_{g,n}(X,\beta),$$ denotes the virtual fundamental class of $\overline{\mathcal M}_{g,n}(X,\beta)$.
\end{defn}

\begin{defn}
The \emph{genus $g$ total descendant potential} of $X$ is the generating function $\mathcal F_g^X\in\La_X[\![\bm t^\bullet]\!]$ of descendant $GW$-invariants of $X$ defined by
\begin{align*}\mathcal F^X_g(\bsy{t}^\bullet,{\bf Q}):=&\sum_{n=0}^\infty\sum_{\beta\in {\rm Eff}(X)}\frac{\bf Q^\beta}{n!}\langle{\bsy\gamma},\dots,\bm\gamma\rangle_{g,n,\beta}^{X}\\
=&\sum_{n=0}^\infty\sum_\bt\sum_{\al_1,\dots,\al_n=0}^N\sum_{p_1,\dots,p_n=0}^\infty\frac{t^{\al_1,p_1}\dots t^{\al_n,p_n}}{n!}\langle\tau_{p_1}T_{\al_1},\dots,\tau_{p_n}T_{\al_n}\rangle^X_{g,n,\bt}\bf Q^\bt.
\end{align*}
Setting $t^{\alpha,0}=t^\alpha$ and $t^{\alpha,p}=0$ for $p>0$, we obtain the \emph{genus $g$ Gromov-Witten potential} of $X$
\beq\label{gwpot}
F^X_g(\bm t,{\bf Q}):=\sum_{n=0}^\infty\sum_{\beta}\sum_{\alpha_1,\dots,\alpha_n=0}^N\frac{t^{\alpha_1}\dots t^{\alpha_n}}{n!}\langle T_{\alpha_1},\dots, T_{\alpha_n}\rangle_{g,n,\beta}^X\bf Q^\bt.
\eneq
It will also be convenient to introduce the  \emph{genus g correlation functions} defined by the derivatives
\beq\label{gcorfun}
\llangle\tau_{d_1}T_{\alpha_1},\dots,\tau_{d_n}T_{\alpha_n}\rrangle_g:=\left.\frac{\partial}{\partial t^{\alpha_1,d_1}}\dots\frac{\partial}{\partial t^{\alpha_n,d_n}}\mathcal F^X_g(\bm t^\bullet,\bf Q)\right|_{\substack{t^{\al,p}=0\text{ for }p>0,\\ t^{\al,0}=t^\al}}.
\eneq
\end{defn}
\subsection{Quantum cohomology}
Let $\bt_1,\dots,\bt_r\in H_2(X,\Z)$ be the homology classes dual to $T_1,\dots, T_r$. By the Divisor axiom, the genus $0$ Gromov-Witten potential $F^X_0(\bm t, {\bf Q})$ can be seen as an element of the ring $\mathbb C[\![t^0,{\bf Q}^{\bt_1}e^{t^1},\dots,{\bf Q}^{\bt_r}e^{t^r},t^{r+1},\dots, t^N]\!]$. In what follows we will be interested in the cases when   $F^X_0$ is a convergent  series 
expansion 
\begin{equation}
\label{26gennaio2020-3}
F^X_0
\in
\mathbb C\left\{t^0,{\bf Q}^{\bt_1}e^{t^1},\dots,{\bf Q}^{\bt_r}e^{t^r},t^{r+1},\dots, t^N\right\}.
\end{equation}
 Without loss of generality we can put ${\bf Q}=1$. Under the  assumption (\ref{26gennaio2020-3}),  $F^X_0(\bm t)$ defines an analytic function in an open neighbourhood $\Omega\subseteq H^\bullet(X,\mathbb C)$ of the point
\begin{equation}
\label{classical1}
t^i=0,\quad i=0,r+1,\dots, N;
\quad\quad
\quad
\operatorname{Re} t^i\to -\infty,\quad i=1,3,\dots, r.
\end{equation}
The function $F^X_0$ is a solution of WDVV equations \cite{kon,manin,tian,voisin}, and thus it defines an analytic Frobenius manifold structure on $\Omega$. Using the canonical identifications of tangent spaces $T_p\Omega\cong H^\bullet(X;\mathbb C)\colon \partial_{t^\alpha}\mapsto T_\alpha$, the unit vector field is $e=\partial_{t^0}\equiv 1$, and the Euler vector field is
\[
E:=c_1(X)+\sum_{\alpha=0}^N\left(1-\frac{1}{2}\deg T_\alpha\right)t^\alpha T_\alpha,
\]
which satisfies
\[\frak L_EF^X_0=(3-\dim_\C X)F_0^X.
\]The Frobenius manifold structure on $\Om$ can be extended by analytic continuation.
The resulting maximal Frobenius structure is called \emph{quantum cohomology of }$X$, denoted $QH^\bullet (X)$.

In the recent paper \cite{cotti2021degenerate}, a useful convergence criterion for formal power series solutions of WDVV equations is given. In the case of quantum cohomologies of Fano varieties, we have the following result.

Assume that $X$ is Fano, and consider the finite dimensional $\C$-algebra $(H^\bullet(X,\C),\circ_0)$, where the product $\circ_0$ is defined by
\[T_\al\circ_0 T_\gm=\sum_{\la=0}^Nc^\la_{\al\gm}T_\la,\qquad c^\la_{\al\gm}:=\sum_{\eps=0}^N\sum_{\bt\in {\rm Eff}(X)}\langle T_\al,T_\gm, T_\eps\rangle^X_{0,3,\bt}\eta^{\eps\la},\quad \al,\gm=0,\dots, N.
\]Notice that the sums defining the structure constants $c^\la_{\al\bt}$ are finite, due to the Fano assumption. 

\begin{thm}[\cite{cotti2021degenerate}]
If $(H^\bullet(X,\C),\circ_0)$ is semisimple, then the Gromov--Witten potential $F^X_0(\bm t,{\bf Q})$ is convergent. That is, condition \eqref{26gennaio2020-3} holds. \qed
\end{thm}

For a further convergence result, beyond the Fano case, see \cite[Sec.\,6]{cotti2021degenerate}. See also \cite{cotti2021riemann}, where the convergence criteria of \cite{cotti2021degenerate} have been generalized to solutions of the more general ``oriented associativity equations'' \cite{manin2005f}.

\section{Monodromy data of quantum cohomology}\label{sec4}

\subsection{Topological-enumerative solution}For $\beta=0,\dots, N$ and $k\in\mathbb N$, introduce the functions
\begin{gather}
\theta_{\beta,k}(\bm t):=\left.\llangle \tau_kT_\bt,1\rrangle_0\right|_{\bf Q=1},\\
\theta_\beta(z,\bm t):=\sum_{k=0}^\infty\theta_{\beta,k}(\bm t)z^k.
\end{gather}
Define the matrix $\Theta(z,\bm t)$ by
\beq
\Theta(z,\bm t)^\alpha_\beta:=\sum_{\la=0}^N\eta^{\alpha\lambda}\frac{\partial\theta_\beta(z,\bm t)}{\partial t^\lambda},\quad \alpha,\beta=0,\dots, N.
\eneq
Denote by $R$ the matrix associated with the morphism $$c_1(X)\cup\colon H^\bullet(X,\C)\to H^\bullet(X,\C),\quad x\mapsto c_1(X)\cup x,$$
w.r.t.\,\,the basis $(T_0,\dots, T_N)$.

Consider the joint system \eqref{eq1}, \eqref{qde} attached to the Frobenius manifold $QH^\bullet(X)$.
\begin{thm}[\cite{dubro2,CDG}]
The matrix $Z_{\rm top}(z,\bm t):=\Theta(z,\bm t)z^\mu z^{R}$ is a fundamental system of solutions of the joint system \eqref{eq1}, \eqref{qde}.\qed
\end{thm}

The fundamental system of solutions $Z_{\rm top}(z,\bm t)$ is called \emph{topological-enumerative solution} of the joint system \eqref{eq1}, \eqref{qde}.

Let $M_0(\bm t)$ be the monodromy matrix defined by
\beq
Z_{\rm top}(e^{2\pi\sqrt{-1}}z,\bm t)=Z_{\rm top}(z,\bm t)M_0(\bm t),\quad z\in\widetilde{\C^*}.
\eneq

\begin{lem}\label{M0}
We have
\beq
M_0(\bm t)=\exp(2\pi\sqrt{-1}\mu)\exp(2\pi\sqrt{-1}R).
\eneq
In particular, $M_0$ does not depend on $\bm t$.\qed
\end{lem}

\subsection{Stokes rays and $\ell$-chamber decomposition}
\begin{defn}
We call \emph{Stokes rays} at a point $p\in \Om$ the oriented rays $R_{ij}(p)$ in $\mathbb C$  defined by
\beq
R_{ij}(p):=\left\{-\sqrt{-1}(\overline{u_i(p)}-\overline{u_j(p)})\rho\colon \rho\in\mathbb R_{+}\right\},
\eneq
where $(u_1(p),\dots, u_n(p))$ is the spectrum of the operator $\mathcal U(p)$ (with a fixed arbitrary order).
\end{defn}
Fix an oriented ray $\ell$ in the universal cover $\widetilde{\C^*}$.  
\begin{defn}
We say that $\ell$ is \emph{admissible} at $p\in\Om$ if the projection of the  ray $\ell$ on $\mathbb C^*$ does not coincide with any Stokes ray $R_{ij}(p)$.
\end{defn}
\begin{defn}
Define the open subset $O_\ell$ of points $p\in \Om$ by the following conditions:
\begin{enumerate}
\item the eigenvalues $u_i(p)$ are pairwise distinct,
\item $\ell$ is admissible at $p$.
\end{enumerate}
We call \emph{$\ell$-chamber} of $\Om$ any connected component of $O_\ell$. 
\end{defn}

\subsection{Stokes fundamental solutions at $z=\infty$}
Fix an oriented ray $\ell\equiv \left\{\arg z=\phi\right\}$ in $\widetilde{\C^*}$. For $m\in\mathbb Z$, define the sectors in $\widetilde{\C^*}$
\begin{gather}
\Pi_{L,m}(\phi):=\left\{z\in\widetilde{\C^*}\colon \phi+2\pi m<\arg z<\phi+\pi+2\pi m\right\},\\
\Pi_{R,m}(\phi):=\left\{z\in\widetilde{\C^*}\colon \phi-\pi+2\pi m < \arg z< \phi+2\pi m\right\}.
\end{gather}
Denote by $\mc B_{X}$ the bifurcation diagram of the quantum cohomology of $X$.
\begin{thm}[\cite{dubro1,dubro2}]
There exists a unique formal solution $Z_{\rm form}(z,\bm t)$ of the joint system \eqref{eq1}, \eqref{qde} of the form
\begin{align}
Z_{\rm form}(z,\bm t)&=\Psi(\bm t)^{-1}G(z,\bm t)\exp(z U(\bm t)),\\
G(z,\bm t)&= I+\sum_{k=1}^\infty\frac{1}{z^k}G_k(\bm t),
\end{align}
where the matrices $G_k(\bm t)$ are holomorphic on $\Om\setminus \mc B_X$.
\end{thm}

\begin{thm}[\cite{dubro1,dubro2}]\label{thfs}
Let $m\in\mathbb Z$. There exist unique fundamental systems of solutions $Z_{L,m}(z,\bm t)$, $Z_{R,m}(z,\bm t)$ of the joint system \eqref{eq1}, \eqref{qde} with asymptotic expansion
\begin{align}
Z_{L,m}(z,\bm t)&\sim Z_{\rm form}(z,\bm t),\quad |z|\to\infty,\quad z\in \Pi_{L,m}(\phi),\\
Z_{R,m}(z,\bm t)&\sim Z_{\rm form}(z,\bm t),\quad |z|\to\infty,\quad z\in \Pi_{R,m}(\phi),
\end{align}
respectively.
\end{thm}

\begin{defn}
The solutions $Z_{L,m}(z,\bm t)$ and $Z_{R,m}(z,\bm t)$ are called \emph{Stokes fundamental solutions} of the joint system \eqref{eq1}, \eqref{qde} on the sectors $\Pi_{L,m}(\phi)$ and $\Pi_{R,m}(\phi)$ respectively.
\end{defn}

\subsection{Monodromy data}
Let $\ell\equiv \left\{\arg z=\phi\right\}$ be an oriented ray in $\widetilde{\C^*}$ and consider the corresponding Stokes fundamental systems of solutions $Z_{L,m}(z,\bm t)$, $Z_{R,m}(z,\bm t)$, for $m\in\mathbb Z$.

\begin{defn}
We define the \emph{Stokes} and \emph{central connection} matrices $S^{(m)}(p)$, $C^{(m)}(p)$, with $m\in\mathbb Z$, at the point $p\in O_\ell$ by the identities 
\begin{align}
Z_{L,m}(z,\bm t(p))=Z_{R,m}(z,\bm t(p)) S^{(m)}(p),\quad z\in\widetilde{\C^*}\\
 Z_{R,m}(z,\bm t(p))=Z_{\rm top}(z,\bm t(p) )C^{(m)}(p), \quad z\in\widetilde{\C^*}.
\end{align}
Set $S(p):=S^{(0)}(p)$ and $C(p):=C^{(0)}(p)$.
\end{defn}

\begin{defn}
The \emph{monodromy data} at the point $p\in O_\ell$ are defined as the $4$-tuple of matrices $(\mu,R,S(p),C(p))$, where
\begin{itemize}
\item $\mu$ is the matrix associated to the grading operator,
\item $R$ is the matrix associated to the operator $c_1(X)\cup\colon H^\bullet(X,\C)\to H^\bullet(X,\C)$,
\item $S(p), C(p)$ are the Stokes and central connection matrices at $p$, respectively.
\end{itemize}
\end{defn}

\begin{defn}
Fix a point $p\in O_\ell$ with canonical coordinates $(u_i(p))_{i=1}^n$. Define the oriented rays $L_j(p,\phi)$, $j=1,\dots,n$, in the complex plane by the equations
\beq
L_j(p,\phi):=\left\{u_j(p)+\rho e^{\sqrt{-1}(\frac{\pi}{2}-\phi)}\colon \rho\in\mathbb R_+\right\}.
\eneq
The ray $L_j(p,\phi)$ is oriented from $u_j(p)$ to $\infty$. We say that $(u_i(p))_{i=1}^n$ are in $\ell$-\emph{lexicographical order} if $L_j(p,\phi)$ is on the left of $L_k(p,\phi)$ for $1\leq j<k\leq n$. 
\end{defn}
In what follows, it is assumed that the $\ell$-lexicographical order of canonical coordinates is fixed at all points of $\ell$-chambers.
\begin{lem}[\cite{CDG1,dubro2}]
If the canonical coordinates $(u_i(p))_{i=1}^n$ are in $\ell$-lexicogra\-phical order at $p\in O_\ell$, then the Stokes matrices $S^{(m)}(p)$, $m\in\mathbb Z$, are upper triangular with $1$'s along the diagonal.
\end{lem}

By Lemma \ref{M0}, the matrices $\mu$ and $R$ determine the monodromy of solutions of the $qDE$,
\beq
M_0:=\exp(2\pi\sqrt{-1}\mu)\exp(2\pi\sqrt{-1}R).
\eneq Moreover, $\mu$ and $R$ do not depend on the point $p$. The following theorem furnishes a refinement of this property.

\begin{thm}[\cite{CDG1,dubro1,dubro2}]\label{thmd}
The monodromy data $(\mu,R, S,C)$ are constant in each $\ell$-chamber. Moreover, they satisfy the following identities:
\begin{gather}
\label{const1}
CS^TS^{-1}C^{-1}=M_0,\\
\label{const2}
S=C^{-1}\exp(-\pi\sqrt{-1}R)\exp(-\pi\sqrt{-1}\mu)\eta^{-1}(C^T)^{-1},\\
\label{const3}
S^T=C^{-1}\exp(\pi\sqrt{-1}R)\exp(\pi\sqrt{-1}\mu)\eta^{-1}(C^T)^{-1}.
\end{gather}
\end{thm}

\begin{thm}[\cite{CDG1}]
The Stokes and central connection matrices $S_m,C_m$, with $m\in\mathbb Z$, can be reconstructed from the monodromy data $(\mu, R,S,C)$:
\beq\label{allsc}
S^{(m)}=S,\quad C^{(m)}=M_0^{-m}C,\quad m\in\mathbb Z.
\eneq
\end{thm}

\begin{rem}\label{coal}
Points of $O_\ell$ are semisimple. The results of \cite{CDG0,CDG1,CG1,CG2} imply that the monodromy data $(\mu,R,S,C)$ are well defined also at points $p\in\Om_{ss}\cap\mc B_{\Om}$, and that Theorem \ref{thmd} still holds true.
\end{rem} 
 \begin{rem}\label{RH}
 From the knowledge of the monodromy data $(\mu,R,S,C)$ the Gromov-Witten potential $F_0^X(\bm t)$ can be recostructed via a Riemann-Hilbert boundary value problem, see \cite{dubro2,guzzetti2,CDG,CDG1}. Hence, the monodromy data may be interpreted as a \emph{system of coordinates} in the space of solutions of $WDVV$ equations.
 \end{rem}

\subsection{Natural transformations of monodromy data}
The definition of the Stokes and central connection matrices is subordinate to several non-canonical choices:
\begin{enumerate}
\item the choice of an oriented ray $\ell$ in $\widetilde{\C^*}$,
\item the choice of an ordering of canonical coordinates $u_1,\dots, u_n$ on each $\ell$-chamber,
\item the choice of signs in \eqref{fvects}, and hence of the branch of the $\Psi$-matrix on each $\ell$-chamber.
\end{enumerate}
Different choices affect the numerical values of the data $(S,C)$.

For different choices of the oriented ray $\ell$, the transformation of $S$ and $C$ can be described in terms of an action of the braid group $\mc B_n$, described in Section \ref{actbn}.\newline
For different choices of ordering of canonical coordinates, the Stokes and central connection matrices transform as follows:
\[
S\mapsto \Pi S \Pi^{-1},\quad C\mapsto C\Pi^{-1},\quad \Pi\text{ permutation matrix}.
\]
For different choices of the branch of the $\Psi$-matrix, we have a transformation of the following type:
\[
S\mapsto I S I,\quad C\mapsto C I,\quad I={\rm diag}(\pm 1,\dots,\pm 1).
\]
See \cite{CDG,CDG1} for more details.

Moreover,  let us also add that the value of all the monodromy data is affected by different choices of the system of flat coordinates $\bm t$.

\begin{prop}\label{flatcord}
Let $(\tilde t^0,\dots,\tilde t^N)$ be a system of flat coordinates on $\Om$ related to $(t^0,\dots, t^N)$ by the transformations
\[\tilde t^\al= A^\al_\bt t^\bt+c^\al,\quad A^\al_\bt,c^\al\in\C,\quad\al,\bt=0,\dots,N.
\]The monodromy data $(\tilde \mu,\tilde R,\tilde S,\tilde C)$, computed w.r.t. the coordinates $\tilde{\bm t}$, are related to the data $(\mu,R,S,C)$, computed w.r.t. $\bm t$, as follows:
\[\tilde \mu= A\mu A^{-1},\quad \tilde R=A R A^{-1},\quad \tilde S=S,\quad \tilde C=AC.
\]
\end{prop}
\proof
The transformation of $\mu,R$ is due to their tensorial nature: they are (1,1)-tensors on $\Om$. Notice that $\tilde \Psi=\Psi A^{-1}$, $\tilde Z_{R,0}=AZ_{R,0}$ and $\tilde{Z}_{\rm top}=AZ_{\rm top}A^{-1}$ so that
\[
\tilde C=\tilde Z_{\rm top}^{-1}\tilde Z_{R,0}=AZ_{\rm top}^{-1}A^{-1}AZ_{R,0}=AC.
\]Equation \eqref{const2}, together with $\tilde \eta=(A^{-1})^T\eta A^{-1}$, shows that $\tilde S=S.$
\endproof

\begin{rem}\label{defocs}In particular, Proposition \ref{flatcord} applies in the case of deformations of the complex structures of $X$. Consider a smooth proper map $f\colon \mc F\to B$ with a connected base space $B$, and set $X_b:=f^{-1}(b)$ with $b\in B$. Given $b_1,b_2\in B$, there exists a diffeomorphism $\phi\colon X_{b_1}\to X_{b_2}$, which allows to identify (co)homology groups:
\[\phi_*\colon H_\bullet(X_{b_1},\Z)\to H_\bullet(X_{b_2},\Z),\quad \phi^*\colon H^\bullet(X_{b_2},\Z)\to H^\bullet(X_{b_1},\Z).
\]Using the isomorphisms $\phi_*,\phi^*$, and by invoking the Deformation Axiom of Gromov-Witten invariants (see e.g.  \cite[Section 7.3]{cox}), we can identify the quantum cohomologies $QH^\bullet(X_{b_1})$ and $QH^\bullet(X_{b_2})$: the deformation of the complex structure just represents a change of flat coordinates on the same Frobenius manifold.
\end{rem}

\subsection{Action of the braid group $\mathcal B_n$}\label{actbn}
Consider the braid group $\mathcal B_n$ with generators $\beta_1,\dots, \beta_{n-1}$ satisfying the relations 
\beq
\beta_{i}\beta_j=\beta_j\beta_i,\quad |i-j|>1,
\eneq
\beq
\beta_i\beta_{i+1}\beta_i=\beta_{i+1}\beta_i\beta_{i+1}.
\eneq
Let $\mathcal U_n$ be the set of upper triangular $(n\times n)$-matrices with $1$'s along the diagonal.
\begin{defn}Given $U\in\mathcal U_n$ define the matrices $A^{\beta_i}(U)$, with $i=1,\dots,n-1$, as follows
\begin{align}
\left(A^{\beta_i}(U)\right)_{hh}:=1,\quad h&=1,\dots,n,\quad h\neq i,i+1,\\
\left(A^{\beta_i}(U)\right)_{i+1,i+1}&=-U_{i,i+1},\\
\left(A^{\beta_i}(U)\right)_{i,i+1}&=\left(A^{\beta_i}(U)\right)_{i+1,i}=1,
\end{align}
and all other entries of $A^{\beta_i}(U)$ are equal to zero.
\end{defn}

\begin{lem}[\cite{dubro1,dubro2}]\label{actbr}The braid group $\mathcal B_n$ acts on $\mathcal U_n\times GL(n,\mathbb C)$ as follows: 
\begin{align*}
\mathcal B_n\times \mathcal U_n\times GL(n,\mathbb C)&\xrightarrow{\hspace{2cm}}\mathcal U_n\times GL(n,\mathbb C)\\
(\beta_i,U,C)&\xmapsto{\quad\quad} (A^{\beta_i}(U)\cdot U\cdot A^{\beta_i}(U),\ C\cdot A^{\beta_i}(U)^{-1})
\end{align*}
We denote by $(U,C)^{\beta_i}$ the  action of $\beta_i$ on $(U,C)$.\qed
\end{lem}

Fix an oriented ray $\ell_o\equiv\left\{\arg z=\phi_o\right\}$ in $\widetilde{\C^*}$, and denote by $\overline{\ell_o}$ its projection on $\widetilde{\C^*}$. Let $p_o\in O_{\ell_o}$, and let $(S_0,C_0)$ be the Stokes and central connection matrices computed at $p_o$ w.r.t.\,\,$\ell_o$, the $\ell_o$-lexicographical order of canonical coordinates $u_i(p_o)$, and a suitable determination of the $\Psi$-matrix at $p_o$. If we let the oriented ray rotate, so that it crosses some Stokes rays $R_{ij}(p_o)$, the values of $(S_0,C_0)$ will change. We can describe this difference of values in terms of the braid group action of Lemma \ref{actbr}.

\begin{thm}[\cite{CDG1,dubro1,dubro2}]\label{thbr0}
Consider a continuous map $\phi\colon [0,1]\to\R$, with $\phi(0)=\phi_o$, and set $\ell(t):=\{\arg z=\phi(t)\}$ for any $t\in [0,1]$. Assume that
\begin{itemize}
\item the rays $\ell(0)$ and $\ell(1)$ are admissible at $p_o$,
\item there exists a unique $t_o\in]0,1[$ such that $\ell(t_o)$ is not admissible at $p_o$,
\item there exist $i_1,\dots, i_k\in\{1,\dots, n\}$, with $|i_a-i_b|>1$ for $a\neq b$, such that the projected ray $\overline{\ell}(t)\subseteq \C$ crosses the rays $\left(R_{i_j,i_j+1}\right)_{j=1}^k$ in the counterclockwise (resp. clockwise) direction, as $t\to t_o^-$. 
\end{itemize}
Denote by $(S_i,C_i)$, with $i=0,1$, the Stokes and central connection matrices at $p_o$ w.r.t.\,\,the oriented ray $\ell(i)$, with $i=0,1$. We have
\[\puqed
(S_1,C_1)=(S_0,C_0)^\bt,\quad \bt=\prod_{j=1}^k\bt_{i_j}\qquad\left(\text{resp. }\bt=\left(\prod_{j=1}^k\bt_{i_j}\right)^{-1}\right).\qedhere
\poqed
\]
\end{thm}

\begin{rem}
An arbitrary rotation of $\ell$ can be decomposed into the composition of elementary rotations satisfying the assumptions of Theorem \ref{thbr0}.
\end{rem}

Furthermore, the braid group action also describes how the values of Stokes and central connection matrices in different $\ell$-chambers (for a fixed oriented rays $\ell$) are related to each other. 

Fix an oriented ray $\ell\equiv\left\{\arg z=\phi\right\}$ in $\widetilde{\C^*}$, and denote by $\overline{\ell}$ its projection on $\mathbb C^*$. Let $\Om_{\ell,1},\Om_{\ell,2}$ be two $\ell$-chambers and let $p_i\in\Om_{\ell,i}$ for $i=1,2$. The difference of values of the Stokes and central connection matrices $(S_1,C_1)$ and $(S_2,C_2)$, at $p_1$ and $p_2$ respectively, can be described by the action of the braid group $\mathcal B_n$ of Lemma \ref{actbr}. 

\begin{thm}[\cite{CDG1,dubro1,dubro2}]\label{thbr}
Consider a continuous path $\gamma\colon[0,1]\to\Om$ such that
\begin{itemize}
\item $\gamma(0)=p_1$ and $\gamma(1)=p_2$,
\item there exists a unique $t_o\in[0,1]$ such that $\ell$ is not admissible at $\gamma(t_o)$,
\item there exist $i_1,\dots, i_k\in\left\{1,\dots,n\right\}$, with $|i_a-i_b|>1$ for $a\neq b$, such that the rays\footnote{Here the labeling of Stokes rays is the one prolonged from the initial point $t=0$.} $\left(R_{i_j,i_j+1}(t)\right)_{j=1}^r$ (resp. $\left(R_{i_j,i_j+1}(t)\right)_{j=r+1}^k$) cross the ray $\overline{\ell}$ in the clockwise (resp. counterclockwise) direction, as $t\to t_o^-$.
\end{itemize}
Then, we have
\[\puqed
(S_2,C_2)=(S_1,C_1)^\beta,\quad \beta:=\left(\prod_{j=1}^r\beta_{i_j}\right)\cdot\left(\prod_{h=r+1}^{k}\beta_{i_h}\right)^{-1}.\qedhere
\poqed
\]
\end{thm}
\begin{rem}
In the general case, the points $p_1$ and $p_2$ can be connected by concatenations of paths $\gamma$ satisfying the assumptions of Theorem \ref{thbr}. 
\end{rem}

\begin{rem}
The action of $\mathcal B_n$ on $(S,C)$ also describes the analytic continuation of the Frobenius manifold structure on $\Om$, see \cite[Lecture 4]{dubro2}.
\end{rem}

\section{$J$-function and Quantum Lefschetz Theorem}\label{sec5}
\subsection{$J$-function and master functions}

\begin{defn}The $J$-function of $X$ the $H^\bullet(X,\La_X)[\![\hbar^{-1}]\!]$-valued function of ${\bm\tau}\in H^\bullet(X,\mathbb C)$ defined by
\[J_X(\bm\tau):=1+\sum_{\al,\la=0}^N\sum_{n=0}^\infty\hbar^{-(n+1)}\llangle \tau_n T_\alpha,1\rrangle_0\eta^{\alpha\lambda}T_\lambda.
\]
\end{defn}

The following result will be crucial for us. For its proof see Appendix \ref{aproof}.
\begin{thm}\label{topJ}Let $\al=0,\dots, N$ and $\delta\in H^2(X,\C)$. The $(1,\al)$-entry of the matrix $\eta Z_{\rm top}(z,\delta)$ equals
\[z^\frac{\dim X}{2}\int_XT_\alpha\cup J_X(\delta+ \log z\cdot c_1(X))\Big|_{\substack{{\bf Q}=1,\\ \hbar=1}}.
\]
\end{thm}

\begin{cor}\label{masJ}Let $\delta\in H^2(X,\C)$. The components  of the function
\[J(\delta+\log z\cdot c_1(X))\Big|_{\substack{{\bf Q}=1,\\ \hbar=1}},
\]w.r.t. any basis of $H^\bullet(X,\C)$,
span the space of master functions $\mc S_\delta(X)$.
\proof
The functions $z^{-\frac{\dim X}{2}}[\eta Z_{\rm top}(z,\delta)]^1_\al$ define a generating set of the space of master functions $\mc S_\delta(X)$. The claim follows by Theorem \ref{topJ}.
\endproof
\end{cor}

In the notations of Section \ref{notX}, set $\delta=\sum_{i=1}^rt^iT_i$. Any formal differential operator $P\in\C[\![\hbar\frac{\der}{\der t^1},\dots,\hbar\frac{\der}{\der t^r},e^{t^1},\dots, e^{t^r},\hbar]\!]$ such that $P J_X(\delta)=0$ is called a \emph{quantum differential operator}. The equation $PY=0$ is called a \emph{quantum differential equation}, see e.g. \cite[Section 10.3]{cox}. By Corollary \ref{masJ}, the master differential equation, defined as in Section \ref{secmas} at a point $\delta$ of the complement of the $\mc A_\La$-stratum of $QH^\bullet(X)$, is equivalent to a differential equation for master functions
\beq
\label{qdiffeq}
\widetilde P_\delta(\vartheta,z)\Phi=0,\quad \vartheta:=z\frac{d}{dz},
\eneq
for a suitable differentiable operator $\widetilde P_\delta$.

\subsection{Twisted Gromov-Witten invariants} Given a holomorphic vector bundle $E\to X$ and an invertible multiplicative\footnote{A characteristic class $\bm c$ is said to be \emph{multiplicative} if $\bm c(E_1\oplus E_2)=\bm c(E_1)\bm c(E_2)$. It is \emph{invertible} if $\bm c(E)$ is invertible in $H^\bullet(Y,\C)$ for any vector bundle $E$ on a manifold $Y$.} characteristic class $\bm c$, one can introduce a $(E,\bm c)$-\emph{twisted} version of the Gromov-Witten theory of $X$.

Given $E$, there exists a complex $0\to E_{g,n,\bt}^0\to E_{g,n,\bt}^1\to 0$ of locally free orbi-sheaves on $\overline{\mathcal M}_{g,n}(X,\beta)$ whose cohomology sheaves are $R^0{\rm ft}_{n+1,*}({\rm ev}_{n+1}^*E)$ and $R^1{\rm ft}_{n+1,*}({\rm ev}_{n+1}^*E)$ respectively. Here, the forgetful and evaluation morphisms ${\rm ft}_{n+1}, {\rm ev}_{n+1}$ at the last marked point fit in the diagram
\[\xymatrix{
&\overline{\mathcal M}_{g,n+1}(X,\beta)\ar[dl]_{{\rm ft}_{n+1}}\ar[dr]^{{\rm ev}_{n+1}}&\\
\overline{\mathcal M}_{g,n}(X,\beta)&&X
}
\]
Let us introduce an \emph{obstruction $K$-class} $E_{g,n,\bt}\in K^0(\overline{\mathcal M}_{g,n}(X,\beta))$, defined as the $K$-theoretic difference
\[E_{g,n,\bt}:=[E_{g,n,\bt}^0]-[E_{g,n,\bt}^1].
\]It is possible to show that such a difference does not depend on the choice of the complex.

\begin{defn}The $(E,\bm c)$-\emph{twisted Gromov-Witten invariants (with descendants)} of $X$ are the intersection numbers 
\[\langle\tau_1^{d_1}\alpha_1\otimes \dots\otimes \tau_n^{d_n}\alpha_n\rangle_{g,n,\beta}^{X, E,\bm c}:=\int_{[\overline{\mathcal M}_{g,n}(X,\beta)]^{\rm virt}}\bm c(E_{g,n,\bt})\cup \prod_{j=1}^n\psi_j^{d_j}\cup{\rm ev}^*_j(\alpha_j),
\]where $\alpha_1,\dots,\alpha_n\in H^\bullet(X,\mathbb C)$.
\end{defn}
\begin{rem}
If $\bm c$ is the trivial characteristic class, we recover the untwisted Gromov-Witten invariants of $X$. 
\end{rem}

\subsection{Quantum Lefschetz Theorem}Introduce a $\C^*$-action on the total space $E$ defined by fiberwise multiplication. The $\C^*$-equivariant Euler class $\bm e$ is invertible over the field of fractions $\Q(\la)$ of $H^\bullet_{\C^*}({\rm pt})\cong \Q[\la]$. Taking $\bm c=\bm e$ we refer to the twisted Gromov-Witten invariants as \emph{Euler-twisted Gromov-Witten invariants}.

Exactly as in the untwisted case, $(E,\bm c)$-twisted Gromov-Witten invariants can be collected in generating functions. In particular, we can introduce the \emph{Euler-twisted $J$-function} as the $H^\bullet(X,\La_X[\la])[\![\hbar^{-1}]\!]$-valued function on $H^\bullet(X,\C)$ by
\beq
J_{E,\bm e}(\bm\tau)=1+\sum_{\al,k,n,\bt}\hbar^{-n-1}\frac{{\bf Q}^\bt}{k!}\langle\tau_nT_\al,1,\bm\tau,\dots,\bm\tau\rangle^{X,E,\bm e}_{0,k+2,\bt}T^\al.
\eneq

Assume now that $E$ is convex\footnote{Globally generated vector bundles, and direct sums of nef line bundles are automatically convex.}, i.e. $H^1(C,f^*E)=0$ for all stable maps $f\colon C\to X$ with $C$ of genus zero. Let $Y$ be a smooth subvariety of $X$ defined by the zero locus of a regular section of $E$.

\begin{thm}[\cite{coatgiv,qLcoates}]\label{ql1}
The non-equivariant limit $J_{E,\bm e}|_{\la=0}$ exists. Moreover, it is related to the function $J_Y$ by the equation
\beq
\label{J}
\iota^*J_{E,\bm e}|_{\la=0}(\bm \tau)\stackrel{\mathclap{\normalfont\mbox{$\iota_*$}}}{=}J_Y(\iota^*\bm \tau),\quad \tau\in H^\bullet(X,\C),
\eneq
where $\iota\colon Y\hookrightarrow X$ is the inclusion.  
\end{thm}

\begin{rem}
The symbol $\stackrel{\mathclap{\normalfont\mbox{$\iota_*$}}}{=}$ means that identity \eqref{J} holds true after application of the morphism $\iota_*\colon\La_X\to \La_Y$ defined by ${\bf Q}^\bt\mapsto {\bf Q}^{\iota_*\bt}$.
\end{rem}

\begin{rem}
If $\dim_\C X>3$, then $\iota^*$ is an isomorphism, by Hyperplane Lefschetz Theorem.
\end{rem}

Assume that $E=\oplus_{i=1}^sL_i$ where $L_i$ are nef line bundles on $X$ such that $c_1(E)\leq c_1(X)$. In such a case, the Quantum Lefschetz Theorem prescribes how to compute the non-equivariant limit $J_{E,\bm e}(\delta)|_{\la=0}$ at points of the small quantum locus $\delta\in H^2(X,\C)$.

Introduce the \emph{hypergeometric modification} $I_{X,Y}$ of the function $J_X$ as follows: write $J_X=\sum_{\bt} J_\bt {\bf Q}^\bt$, and for $\delta\in H^2(X,\C)$ define
\beq
\label{Ifun}
I_{X,Y}(\delta):=\sum_{\bt}J_{\bt}(\delta){\bf Q}^\bt\prod_{i=1}^s\prod_{m=1}^{\langle c_1(L_i),\bt\rangle}(c_1(L_i)+m\hbar).
\eneq

\begin{thm}[\cite{coatgiv}]\label{ql2}
The function $I_{X,Y}$ admits an expansion of the form
\beq
I_{X,Y}(\delta)=F(\delta)+\frac{1}{\hbar} G(\delta)+O\left(\frac{1}{\hbar^2}\right),\quad \delta\in H^2(X,\C),
\eneq
where $F$ is $H^0(X,\La_X)$-valued and $G$ takes values in $H^0(X,\La_X)\oplus H^2(X,\La_X)$. Moreover, we have
\beq
J_{E,\bm e}(\phi(\delta))|_{\la=0}=\frac{I_{X,Y}(\delta)}{F(\delta)},\quad \phi(\delta):=\frac{G(\delta)}{F(\delta)}.
\eneq
\end{thm}

\begin{prop}[\cite{coatgiv,qperiods}]\label{ql3}
Moreover, if $c_1(X)>c_1(E)$, then we have
\[ F(\delta)\equiv 1,\quad G(\delta)=\delta+H(\delta)\cdot 1,\quad H(\delta)=\sum_{\bt}\left(w_\bt {\bf Q^\bt}e^{\int_\bt\delta}\right)\cdot \delta_{1,\langle\bt,c_1(X)-c_1(E)\rangle},
\]for suitable rational coefficients $w_\bt\in\Q$.
\end{prop}

\proof
The function $I_{X,Y}(\delta)$ is homogeneous of degree 0 w.r.t. the gradings 
\[\deg{\bf Q}^\bt=\int_\bt c_1(X)-\int_\bt c_1(E),\quad \deg\hbar=1,\quad \deg T_\al=k\text{ if }T_\al\in H^{2k}(X,\C).
\]This is easily seen from the expansion of $J_X$ given in Lemma \ref{restJ}. Hence, $F(\delta)$ is given from the only contribution of the term $J_0(\delta)=1+\frac{\delta}{\hbar}+\dots$, and $H(\delta)$ from the terms for which $\deg{\bf Q}^\bt=1$.
\endproof

\subsection{Inequality for dimensions of spaces of master functions} Let $Y\subseteq X$ be the zero locus of a regular section of a vector bundle $E\to X$, sum of nef line bundles, with $c_1(E)< c_1(X)$. Denote by $\iota\colon Y\to X$ the inclusion. We always assume that both $X$ and $Y$ have vanishing odd comology.

For a point $\bm \tau\in QH^\bullet(X)$, denote by $\mc S_{\bm \tau}(X):=\mc S_{\bm\tau}(QH^\bullet(X))$ the space of master functions as $\bm \tau$.

\begin{thm}\label{inmasY}
Let $\delta\in H^2(X,\C)$. We have
\beq\label{ineqS}
\dim_\C \mc S_{\iota^*\delta}(Y)\leq\dim_\C \mc S_{\delta+c}(X),
\eneq
where $c:=c_1(X)-c_1(E)$. 
\end{thm}

\proof
By the adjunction formula, we have $\iota^* c=c_1(Y)$. The components of the function $J_X(\delta+c \log z)\rqh$, w.r.t. any basis of $H^\bullet(X,\C)$, span the space $\mc S_{\delta+c}(X)$. Analogously, the components of the function $J_Y(\iota^*\delta+c_1(Y)\log z)\rqh$, w.r.t. any basis of $H^\bullet(Y,\C)$, span the space $\mc S_{\iota^*\delta}(Y)$.

By Theorems \ref{ql1}, \ref{ql2} and Proposition \ref{ql3}, we have
\[J_Y(\iota^*\delta+c_1(Y)\log z)\rqh=e^{- zH(\delta)}\cdot\iota^*I_{X,Y}(\delta+c \log z)\rqh.
\]
The components of the r.h.s. are obtained by linear combinations and rescaling of the components of $J_X(\delta+c\log z)\rqh$: such a linear combination is due to the hypergeometric modification \eqref{Ifun}, namely the $\cup$-multiplication by an invertible class. 
The claim follows.
\endproof

\begin{thm}Let $Y$ be a hyperplane section of $X$.
Assume that $d:=\dim_\C X$ is odd, and that the following inequalities of Betti numbers hold true:
\beq \label{betneq}
b_{d-1}(X)<\frac{1}{2}b_{d-1}(Y).
\eneq
Then $\iota^*(H^2(X,\C))$ is contained in the $\mc A_\La$-stratum of the Frobenius manifold $QH^\bullet(Y)$. In particular, along $\iota^*(H^2(X,\C))$ the canonical coordinates of $QH^\bullet(Y)$ coalesce.
\end{thm}

\proof
From Hyperplane Lefschetz Theorem we deduce that \eqref{betneq} holds true if and only if $\dim_\C H^\bullet(X,\C)<\dim_\C H^\bullet(Y,\C)$. Then we have $\dim_\C \mc S_{\iota^*\delta}(Y) <\dim_\C H^\bullet(Y,\C)$ for any $\delta\in H^2(X,\C)$, by \eqref{ineqS}. Hence, the master differential equation of $QH^\bullet(Y)$ at $\iota^*\delta$ is not of order $\dim_\C H^\bullet(Y,\C)$. This implies that the denominator of $\det \La$ is identically zero at $\iota^*\delta$. The last statement follows from Lemma \ref{ala} and Theorem \ref{ib}.
\endproof

\section{Borel-Laplace $(\bm \alpha,\bm \beta)$-multitransforms}\label{BL}

\subsection{Algebras of Ribenboim's generalized power series} Let $(M,+,0)$ be a \emph{monoid}, i.e. a commutative semigroup with neutral element. We will say that a partial order relation $\leq$ on $M$ defines a \emph{strictly ordered monoid} $(M,+,0,\leq)$ if the following compatibility condition holds true:
\[\text{if $a<b$, then $a+s<b+s$ for all $s\in M$.}
\]
Let $R$ be a commutative ring with unit. The set $R[\![M]\!]:=R^M$ of all functions $f\colon M\to R$ is equipped with a natural $R$-module structure, w.r.t. pointwise addition and multiplication by scalars. An element $f\in R[\![M]\!]$ will usually be denoted by
\[f=\sum_{a\in M}f(a)Z^a,
\]where $Z$ is an indeterminate. Given two functions $f,g\in R[\![M]\!]$, we could be tempted to define their product as
\begin{equation}\label{prodgenser}f\cdot g:=\sum_{s\in M}\left(\sum_{(p,q)\in X_s(f,g)}f(p)\cdot g(q)\right)Z^{s},
\end{equation}where we set
\[X_s(f,g):=\left\{(p,q)\in M\times M\colon p+q=s,\quad f(p)\neq 0,\quad g(q)\neq 0\right\}.
\]In general the set $X_s(f,g)$ is \emph{not} finite, and consequently the product $f\cdot g$ could be not defined.
\begin{defn}Let $(M,+,0,\leq)$ be a {strictly ordered monoid}.
The $R$-submodule of $R[\![M]\!]$ which consists of all functions $f\colon M\to R$ whose support $${\rm supp}(f):=\left\{s\in M\colon f(s)\neq 0\right\}$$ is 
\begin{enumerate}
\item \emph{Artinian}, i.e. every subset of ${\rm supp}(f)$ admits a minimal element, 
\item \emph{narrow}, i.e. every subset of ${\rm supp}(f)$ of pairwise incomparable elements is finite,
\end{enumerate}
is called the set of \emph{generalized power series} with coefficients in $R$ and exponents in $M$. It is denoted by $R[\![M,\leq]\!]$. 
\end{defn}

\begin{prop}[\cite{riben1,riben2}]
Given $f,g\in R[\![M,\leq]\!]$, the set $X_s(f,g)$ is finite, and the product \eqref{prodgenser} is well-defined. The set $R[\![M,\leq]\!]$ inherits the structure of an associative $R$-algebra.
\end{prop}
\begin{rem}
If $(M,\leq)$ is itself Artinian and narrow, then all its subsets are Artinian and narrow. Consequently $R[\![M,\leq]\!]=R[\![M]\!]$.
\end{rem}

\subsection{The algebra $\mathscr F_{\bm \kappa}(A)$}\label{fka}

Let  ${\bsy\kappa}:=(\kappa_1,\dots,\kappa_h)\in(\C^*)^h$. Consider an associative, commutative, unitary and finite dimensional $\mathbb C$-algebra $(A,+,\cdot,1_A)$. Denote by ${\rm Nil}(A)$ the nilradical of $A$, that is ${\rm Nil}(A):=\{a\in A\colon \exists\,n\in\N\,\,\text{s.t.}\,\,a^n=0\}$.
\vskip3mm
 Set $\N_A:=\{n\cdot 1_A\colon n\in\N\}$. Define the monoid $M_{A,{\bsy\kappa}}$ as the (external) direct sum of monoids
\[M_{A,{\bsy\kappa}}:=\left(\bigoplus_{j=1}^h\kappa_j\mathbb N_A\right)\oplus{\rm Nil}(A).\]
We have two maps $\nu_{\bsy\kappa}\colon M_{A,{\bsy\kappa}}\to \mathbb N^h$ and $\iota_{\bsy\kappa}\colon M_{A,{\bsy\kappa}}\to A$ defined by
\[\nu_{\bsy\kappa}((\kappa_in_i1_A)_{i=1}^h,r):=(n_i)_{i=1}^h,\quad \iota_{\bsy\kappa}((\kappa_in_i1_A)_{i=1}^h,r):=\sum_{i=1}^h\kappa_in_i1_A+r.
\]

On $M_{A,{\bsy\kappa}}$ we can define the partial order 
\[x\leq y \quad\text{iff}\quad \nu_{\bsy\kappa}(x)\leq \nu_{\bsy\kappa}(y),
\]the order on $\mathbb N^h$ being the lexicographical one. This order makes $(M_{A,{\bsy\kappa}},\leq)$ a strictly ordered monoid. 

We denote by $\mathscr F_{\bsy\kappa}(A)$ the ring $A[\![M_{A,{\bsy\kappa}},\leq]\!]$.

By universal property of the direct sums of monoids, the natural inclusions $M_{A,\kappa_i}\to M_{A,{\bsy\kappa}}$ induce a unique morphism
\[\rho_{\bsy\kappa}\colon \bigoplus_{i=1}^hM_{A,\kappa_i}\to M_{A,\bsy\kappa}.
\]

\begin{defn}
Let $r_o\in{\rm Nil}(A)$. We say that an element $f\in \mathscr F_{\bm \kappa}(A)$ is \emph{concentrated at $r_o$} if 
\[{\rm  supp}(f)\subseteq \left(\bigoplus_{i=1}^h\kappa_i\N_A\right)\times\{r_o\}.
\]
\end{defn}

\subsection{Formal Borel-Laplace $(\bm\alpha,\bm\beta)$-multitransforms}\label{formal}

Given two $h$-tuples $\bm \al,\bm\bt\in(\C^*)^h$, we set $\bm\al\cdot\bm\bt:=(a_i\bt_i)_{i=1}^h$, and $\bm\al^{-1}:=\left(\frac{1}{\al_i}\right)_{i=1}^h$.

\begin{defn}
Let $F\in\C[\![x]\!]$ be a formal power series $F(x)=\sum_{k=0}^\infty a_kx^k$. For $\al\in{{\rm Nil}(A)}$ define $F(\al)\in A$ by the finite sum 
\[F(\al)=\sum_{k=0}^\infty a_k\al^k.
\]
If $F$ is invertible, i.e. $a_0\neq 0$, then $F(\al)$ is invertible in $A$.
\end{defn}
In what follows we will usually take $F(x)=\Gamma(\lambda+x)$ with $\lambda\in\C\setminus\mathbb Z_{\leq 0}$, where $\Gamma$ denotes the Euler Gamma function.

\begin{defn}
Let ${\bsy\alpha,\bsy \beta, \bsy \kappa}\in(\mathbb C^*)^h$. We define the \emph{Borel $(\bsy\alpha,\bsy\beta)$-multitransform} as the $A$-linear morphism
\[\mathscr B_{\bsy\alpha,\bsy\beta}\colon \bigotimes_{j=1}^h\mathscr F_{\kappa_j}(A)\to \mathscr F_{\bm\al^{-1}\cdot\bm\bt^{-1}\cdot\bm\kappa}(A),
\]which is defined, on decomposable elements, by
\[\mathscr B_{\bm\alpha,\bm\beta}\left(\bigotimes_{j=1}^h\left(\sum_{s_j\in M_{A,\kappa_j}}f_{s_j}^jZ^{s_j}\right)\right):=\sum_{\substack{s_j\in M_{A,\kappa_j}\\ j=1,\dots, h}}\frac{\prod_{i=1}^hf^i_{s_i}}{\Gamma\left(1+\sum_{\ell=1}^h\iota_{ \kappa_\ell}(s_\ell)\beta_\ell\right)}Z^{\rho_{\bm \kappa}\left(\oplus_{\ell=1}^h \frac{s_\ell}{\al_\ell\bt_\ell }\right)}.
\]
\end{defn}

\begin{defn}
Let ${\bsy\alpha,\bsy \beta, \bsy \kappa}\in(\mathbb C^*)^h$. We define the \emph{Laplace $(\bsy\alpha,\bsy\beta)$-multitransform} as the $A$-linear morphism
\[\mathscr L_{\bsy\alpha,\bsy\beta}\colon \bigotimes_{j=1}^h\mathscr F_{\kappa_j}(A)\to \mathscr F_{\bm\al\cdot\bm\bt\cdot\bm\kappa}(A),
\]which is defined, on decomposable elements, by
\[\mathscr L_{\bm\alpha,\bm\beta}\left(\bigotimes_{j=1}^h\left(\sum_{s_j\in M_{A,\kappa_j}}f_{s_j}^jZ^{s_j}\right)\right):=\sum_{\substack{s_j\in M_{A,\kappa_j}\\ j=1,\dots, h}}\left(\prod_{i=1}^hf^i_{s_i}\right)\Gamma\left(1+\sum_{\ell=1}^h\iota_{ \kappa_\ell}(s_\ell)\beta_\ell\right)Z^{\rho_{\bm \kappa}(\oplus_{\ell=1}^h \al_\ell\bt_\ell s_\ell)}.
\]
\end{defn}

In the case $h=1$, the Borel-Laplace $(\bm \al,\bm\bt)$-multitransform simplify as follows.
\begin{defn}
Given $\alpha,\beta\in\mathbb C^*$, we define two $A$-linear maps 
\[ \mathscr B_{\alpha,\beta}\colon\mathscr F_{\kappa}(A)\to\mathscr F_{\frac{\kappa}{\alpha\beta}}(A),\quad \mathscr L_{\alpha,\beta}\colon\mathscr F_\kappa(A)\to\mathscr F_{\alpha\beta\kappa}(A),\quad \kappa\in\mathbb C^*
\]
called respectively \emph{$(\alpha,\beta)$-Borel and Laplace transforms}, through the formulæ
\begin{align*}
\mathscr B_{\alpha,\beta}\left[\sum_{s\in M_{A,\kappa}}f_sZ^s\right]&:=\sum_{s\in M_{A,\kappa}}\frac{f_s}{\Gamma(1+\beta s)}Z^\frac{s}{\alpha\beta},\\
\mathscr L_{\alpha,\beta}\left[\sum_{s\in M_{A,\kappa}}f_sZ^s\right]&:=\sum_{s\in M_{A,\kappa}}f_s\Gamma(1+\beta s)Z^{\alpha\beta s}.
\end{align*}
\end{defn}

\begin{thm}
The Borel-Laplace $(\al,\bt)$-transform are inverses of each other, i.e. 
\[
\puqed\mathscr B_{\al,\bt}\circ\mathscr L_{\al,\bt}={\rm Id},\quad \mathscr L_{\al,\bt}\circ\mathscr B_{\al,\bt}={\rm Id}.\qedhere
\poqed
\]
\end{thm}

\subsection{Analytic Borel-Laplace $(\bm\al,\bm\bt)$-multitransforms}

\begin{figure}[ht!]
\centering
\def\svgscale{1}
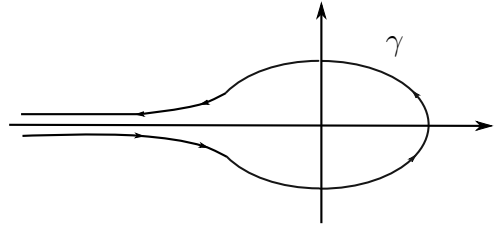
\caption{Hankel-type contour of integration defining Borel $(\bm \al,\bm\bt)$-multitransform.}
\label{gammahankel}
\end{figure}

\begin{defn}
Let $\bm\al,\bm\bt\in(\C^*)^h$. The Borel $(\bm\al,\bm\bt)$-multitransform of an $h$-tuple of $A$-valued functions $(\Phi_1,\dots, \Phi_h)$ is defined, when the integral exists, by
\beq
\mathscr B_{\bm\al,\bm\bt}[\Phi_1,\dots,\Phi_h](z):=\frac{1}{2\pi i}\int_\gamma\prod_{j=1}^h\Phi_j\left(z^\frac{1}{\al_j\bt_j}\la^{-\bt_j}\right)e^\la\frac{d\la}{\la},
\eneq
where $\gamma$ is a Hankel-type contour of integration, see Figure \ref{gammahankel}.
\end{defn}

\begin{defn}
Let $\bsy{\alpha}:=(\alpha_1,\dots,\alpha_h)$, and $\bsy{\beta}:=(\beta_1,\dots,\beta_h)$ be two $h$-tuples in $(\mathbb C^*)^h$. The $(\boldsymbol{\alpha},\boldsymbol{\beta})$-Laplace transform of an $h$-tuple of functions $(\Phi_1,\dots, \Phi_h)$ is defined, when the integral exists, by
\beq
\mathscr L_{\bsy{\alpha},\bsy{\beta}}\left[\Phi_1,\dots,\Phi_h\right](z):=\int_0^\infty\prod_{i=1}^h\Phi_i(z^{\alpha_i\beta_i}\lambda^{\beta_i})\exp(-\lambda)d\lambda.
\eneq
\end{defn}

\begin{prop}
Let $(e_1,\dots, e_n)$ be a basis of $A$ and $\Phi_1,\dots,\Phi_h$ be $A$-valued functions. Write $\Phi_i=\sum_j\Phi_i^je_j$ for $\C$-valued component functions $\Phi_i^j$. The components of $\mathscr B_{\bm \al,\bm\bt}[\Phi_1,\dots,\Phi_h]$ (resp. $\mathscr L_{\bm \al,\bm\bt}[\Phi_1,\dots,\Phi_h]$) are $\C$-linear combinations of the $h\cdot n$ $\C$-valued functions $\mathscr B_{\bm \al,\bm\bt}[\Phi_1^{i_1},\dots,\Phi_h^{i_h}]$ (resp. $\mathscr L_{\bm \al,\bm\bt}[\Phi_1^{i_1},\dots,\Phi_h^{i_h}]$), where $(i_1,\dots, i_h)\in\left\{1,\dots,n\right\}^{\times h}$.
\end{prop}

\proof
Let $c^i_{jk}\in\C$ be the structure constants of the algebra $A$, so that $e_je_k=\sum_ic^i_{jk}e_i$. We have
\[\mathscr B_{\bm \al,\bm\bt}[\Phi_1,\dots,\Phi_h]=\sum_{\bm a,\bm i} c_{i_1i_2}^{a_1}c_{a_1i_3}^{a_2}\dots c_{a_{h-2}i_h}^{a_{h-1}}e_{a_{h-1}}\mathscr B_{\bm\al,\bm\bt}[\Phi_1^{i_1},\dots,\Phi_h^{i_h}].
\]Similarly for the Laplace multitransform.
\endproof

\subsection{Analytification of elements of $\mathscr F_{\bm \kappa}(A)$} Let $s=((\kappa_in_i1_A)_{i=1}^h,r)\in M_{A,\bm\kappa}$. We define the {\it analytification} $\widehat{Z^s}$ of the monomial $Z^s\in\mathscr F_{\bm \kappa}(A)$ to be the $A$-valued holomorphic function 
\[\widehat{Z^s}\colon \widetilde{\C^*}\to A,\quad \widehat{Z^s}(z):=z^{\sum_{i=1}^h\kappa_i n_i}\sum_{j=1}^\infty\frac{r^j}{j!}\log^j z.
\]Notice that the sum is finite, since $r\in{\rm Nil}(A)$.

Let $f\in\mathscr F_{\bm \kappa}(A)$ be a series
\[f(Z)=\sum_{s\in M_{A,\bm\kappa}}f_aZ^s,\quad\text{such that}\quad  {\rm card\ supp}(f)\leq \aleph_0.
\]

The {\it analytification} $\hat f$ of $f$ is the $A$-valued holomorphic function defined, if the series absolutely converges, by
\[\widehat f\colon W\subseteq \widetilde{\C^*}\to A,\quad \widehat f(z):=\sum_{s\in M_{A,\bm\kappa}}f_a\widehat{Z^s}(z).
\]

\begin{thm}\label{BLAF}
Let $f_i\in \mathscr F_{ \kappa_i}(A)$ such that 
\begin{itemize}
\item ${\rm card\ supp}(f_i)\leq\aleph_0$ for $i=1,\dots, h$,
\item the functions $\widehat f_i$ are well defined on $\R_+$.
\end{itemize}
We have
\[\reallywidehat{\mathscr B_{\bm\al,\bm\bt}[\bigotimes_{j=1}^hf_j]}=\mathscr B_{\bm\al,\bm\bt}[\widehat f_1,\dots, \widehat f_h],
\]
\[\reallywidehat{\mathscr L_{\bm\al,\bm\bt}[\bigotimes_{j=1}^hf_j]}=\mathscr L_{\bm\al,\bm\bt}[\widehat f_1,\dots, \widehat f_h],
\]provided that both sides are well-defined.
\end{thm}
\proof
It is sufficient to prove the statement on monomials $Z^{s_1},\dots, Z^{s_h}$. Let $s_j=(\kappa_j n_j 1_A,r_j)$ for $j=1,\dots, h$. We have
\begin{align*}\mathscr B_{\bm \al,\bm \bt}[\otimes_{j=1}^h Z^{s_j}]&=\frac{1}{\Gamma\left(1+\sum_{\ell=1}^h\iota_{ \kappa_\ell}(s_\ell)\beta_\ell\right)}Z^{\rho_{\bm \kappa}\left(\oplus_{\ell=1}^h \frac{s_\ell}{\al_\ell\bt_\ell }\right)}\\
&=\frac{1}{\Gamma\left(1+\sum_{\ell=1}^h(\kappa_\ell n_\ell 1_A+r_\ell)\beta_\ell\right)}Z^{((\frac{\kappa_j}{\al_j\bt_j} n_j 1_A)_{j=1}^h,\ \frac{r_1}{\al_1\bt_1}+\dots+\frac{r_h}{\al_h\bt_h})}.
\end{align*}
Hence, we have
\[
\reallywidehat{\mathscr B_{\bm \al,\bm \bt}[\otimes_{j=1}^h Z^{s_j}]}(z)=\frac{z^{\sum_{i=1}^h\frac{\kappa_i n_i}{\al_i\bt_i}}}{\Gamma\left(1+\sum_{\ell=1}^h(\kappa_\ell n_\ell 1_A+r_\ell)\beta_\ell\right)}\sum_{j=1}^\infty\frac{(\frac{r_1}{\al_1\bt_1}+\dots+\frac{r_h}{\al_h\bt_h})^j}{j!}\log^j z.
\]
On the other hand, we have
\[\widehat{Z^{s_j}}(z)=z^{\kappa_j n_j}\sum_{\ell}\frac{r_j^\ell}{\ell!}\log^\ell z,
\]so that
\begin{align*}
\mathscr B_{\bm\al,\bm\bt}[\widehat{Z^{s_1}},\dots, \widehat{Z^{s_h}}](z)&=\frac{1}{2\pi i}\int_\gamma\prod_{j=1}^h\widehat{Z^{s_j}}\left(z^\frac{1}{\al_j\bt_j}\la^{-\bt_j}\right)e^\la\frac{d\la}{\la}\\
&=\frac{1}{2\pi i}\int_\gamma e^\la\frac{d\la}{\la}\prod_{j=1}^h(z^\frac{1}{\al_j\bt_j}\la^{-\bt_j})^{\kappa_j n_j}\sum_{\ell}\frac{r_j^\ell}{\ell!}\log^\ell (z^\frac{1}{\al_j\bt_j}\la^{-\bt_j})\\
&=\frac{z^{\sum_{i=1}^h\frac{\kappa_i n_i}{\al_i\bt_i}}}{2\pi i}\int_\gamma e^\la\frac{d\la}{\la^{1+\sum_{\ell=1}^h\kappa_\ell n_\ell\bt_\ell}}\prod_{j=1}^h\sum_{\ell}\frac{r_j^\ell}{\ell!}\log^\ell (z^\frac{1}{\al_j\bt_j}\la^{-\bt_j})\\
&=\frac{z^{\sum_{i=1}^h\frac{\kappa_i n_i}{\al_i\bt_i}}}{2\pi i}\int_\gamma e^\la\frac{d\la}{\la^{1+\sum_{\ell=1}^h\kappa_\ell n_\ell\bt_\ell}}\sum_{\ell_1,\dots,\ell_h}\prod_{j=1}^h\frac{r_j^{\ell_j}}{\ell_j!}\log^{\ell_j} (z^\frac{1}{\al_j\bt_j}\la^{-\bt_j}).\\
\end{align*}
We have
\begin{align*}
\prod_{j=1}^h\frac{r_j^{\ell_j}}{\ell_j!}\log^{\ell_j} (z^\frac{1}{\al_j\bt_j}\la^{-\bt_j})&=\prod_{j=1}^h\sum_{w,u=0}^\infty\frac{r_j^{\ell_j}}{w! u!}\left(\frac{\log z}{\al_j\bt_j}\right)^{w}\left(-\bt_j\log\la\right)^{u}\delta_{w+u,\ell_j}\\
&=\sum_{\substack{w_1,\dots,w_h \\ u_1,\dots u_h}}\prod_{j=1}^h\frac{r_j^{\ell_j}}{w_j! u_j!}\left(\frac{\log z}{\al_j\bt_j}\right)^{w_j}\left(-\bt_j\log\la\right)^{u_j}\delta_{w_j+u_j,\ell_j}.
\end{align*}
We have
\[\frac{1}{2\pi i}\int_\gamma e^\la\frac{d\la}{\la^{1+\sum_{\ell=1}^h\kappa_\ell n_\ell\bt_\ell}}(-\log \la)^{u_j}=\left(\frac{1}{\Gamma}\right)^{(u_j)}\left(1+\sum_{\ell=1}^h\kappa_\ell n_\ell\bt_\ell\right),
\]
because of the Hankel formula (see e.g. \cite{NIST})
\[\frac{1}{\Gamma(z)}=\frac{1}{2\pi i}\int_\gamma e^\la \frac{d\la}{\la^z}.
\]Thus, we have
\begin{align*}
&\mathscr B_{\bm\al,\bm\bt}[\widehat{Z^{s_1}},\dots, \widehat{Z^{s_h}}](z)\\ &=z^{\sum_{i=1}^h\frac{\kappa_i n_i}{\al_i\bt_i}}\sum_{\substack{\ell_1,\dots,\ell_h\\w_1,\dots,w_h \\ u_1,\dots u_h}}\prod_{j=1}^h\frac{r_j^{\ell_j}\bt_j^{u_j}}{w_j! u_j!}\left(\frac{\log z}{\al_j\bt_j}\right)^{w_j}\left(\frac{1}{\Gamma}\right)^{(u_j)}\left(1+\sum_{\ell=1}^h\kappa_\ell n_\ell\bt_\ell\right)\delta_{w_j+u_j,\ell_j}.
\end{align*}
This coincides with the formula of $\reallywidehat{\mathscr B_{\bm \al,\bm \bt}[\otimes_{j=1}^h Z^{s_j}]}(z)$. The proof for the Laplace multitransform is similar, based on the identity
\[
\puqed
\Gamma(z)=\int_0^\infty \la^{z-1}e^{-\la}d\la.\qedhere
\poqed
\]

\section{Integral representations of solutions of $qDE$s}\label{sec7}

\subsection{$J_X$-function as element of $\mathscr F_{\bm \kappa}(X)$}
Let $X$ be a variety with nef anticanonical bundle\footnote{We recall that this means $\int_C c_1(X)\geq 0$ for all curves $C$ in $X$. If the strict inequality holds true for any $C$, then $X$ is Fano by Nakai-Moishezon Theorem. Varieties with nef anticanonical bundle can be thought as an interpolation between Fano and Calabi-Yau varieties.}. Introduce the basis $(\bt_1,\dots,\bt_r)$ of $H_2(X,\Z)$ Poincar\'e dual to $(T^1,\dots, T^n)$, so that
\[\int_{\bt_i}T_j=\int_X T^i\cup T_j=\delta_{i,j}.
\]Set 
\[c_1(X)=\sum_{j=1}^{\frak r} c^{\al_{i_j}} T_{\al_{i_j}},\quad  c^{\al_{i_j}}\in\N^*.
\]
\vskip3mm
Consider the $\C$-algebra $H^\bullet(X,\C)$. For brevity, we set $\mathscr F_{\bm \kappa}(X):=\mathscr F_{\bm \kappa}(H^\bullet(X,\C))$ for any $\bm\kappa\in(\C^*)^h$.
\vskip3mm
The $J_X$-function, restricted to the small quantum locus of $QH^\bullet(X)$ admits the following expansion:
\[
\left.J_X(\delta+\log z\cdot c_1(X))\right\rqh=e^\delta z^{c_1(X)}+\sum_\al\sum_{\bt\neq 0}\sum_{k=0}^\infty e^\delta z^{\int_\bt c_1(X)}z^{c_1(X)}\langle\tau_k T_\al,1\rangle_{0,2,\bt}^X T^\al.
\]
Such a series can be seen as an element of $\mathscr F_{\bm \kappa}(X)$ for different choices of $\bm\kappa$. We describe two possible choices. In both cases, we have a series in $\mathscr F_{\bm \kappa}(X)$ concentrated at $c_1(X)$.
\vskip3mm
{\bf Choice 1. } Set $h=1$ and $\kappa=c$, where $c$ is a common divisor of the numbers $c^{\al_{i_1}},\dots, c^{\al_{i_{\frak r}}}$. The series can be rearranged as follows
\[J_X(\delta+\log z\cdot c_1(X))\rqh=\sum_{d\in\N}J_d(\delta)z^{dc+c_1(X)},
\]where
\[J_d(\delta)=e^\delta\sum_{\al,k}\langle\tau_k T_\al,1\rangle_{0,2,d\cdot {\rm PD}(T)}T^\al,\quad d\in\N,\quad T\in H^2(X,\Z),\quad c_1(X)=c T.
\]In particular $J_0(\delta)=e^\delta$.
\vskip3mm
{\bf Choice 2. }Set $h=\frak r$ and $\bm \kappa=(c^{\al_{i_1}},\dots, c^{\al_{i_\frak r}})$. 
By expanding the sum over $\bt$ over the basis $(\bt_1\dots,\bt_r)$, the sum above becomes 
\[
J_X(\delta+\log z\cdot c_1(X))\rqh=\sum_{d\in\N^{\frak r}}J_d(\delta)z^{d_1c^{\al_{i_1}}+\dots+d_{\frak r}c^{\al_{i_{\frak r}}}+c_1(X)},
\]
where
\[J_d(\delta)=e^\delta\sum_{\al,k}\langle\tau_k T_\al,1\rangle_{0,2,d_1\bt_{\al_{i_1}}+\dots+d_{\frak r}\bt_{\al_{i_{\frak r}}}}T^\al,\quad d\in\N^{\frak r}.
\]In particular $J_0(\delta)=e^\delta$.

\subsection{Integral representations of the first kind}
Let $X$ be a Fano smooth projective variety. 
Assume that $\det T_X=L^{\otimes \ell}$ with $L$ ample line bundle. Let $\iota\colon Y\subseteq X$ be a smooth subvariety defined as the zero locus of a regular section of the vector bundle $E=\bigoplus_{j=1}^s L^{\otimes d_j}$, where the numbers $d_j\in\N^*$ are such that $\sum_{j=1}^sd_j<\ell.$

\begin{thm}\label{TH1}
Let $\delta\in H^2(X,\C)$, and $\mc S_\delta(X)$ the corresponding space of master functions of $QH^\bullet(X)$. There exists a complex number $c_\delta\in\C$ such that the space of master functions $\mc S_{\iota^*\delta}(Y)$ is contained in image of the $\C$-linear map $\mathscr S_{(\ell,\bm d)}\colon \mc S_\delta(X)\to \mc O(\widetilde{\C^*})$ defined by 
\[\mathscr S_{(\ell,\bm d)}[\Phi](z):=e^{-c_\delta z}\mathscr L_{\frac{\ell-\sum_{i=1}^s d_i}{d_s},\frac{d_s}{\ell-\sum_{i=1}^{s-1} d_i}}\circ\dots\circ\mathscr L_{\frac{\ell-d_1-d_2}{d_2},\frac{d_2}{\ell-d_1}}\circ\mathscr L_{\frac{\ell-d_1}{d_1},\frac{d_1}{\ell}}[\Phi](z).
\]
In other words, any element of $\mc S_{\iota^*\delta}(Y)$ is of the form
\beq\label{intth1}
e^{-c_\delta z}\int_0^\infty\dots\int_0^\infty\Phi\left(z^{\frac{\ell-\sum_{j=1}^sd_j}{\ell}}\prod_{i=1}^s\zeta_i^\frac{d_i}{\ell}\right)e^{-\sum_{i=1}^s\zeta_i}d\zeta_1\dots d\zeta_s,
\eneq
for some $\Phi\in\mc S_\delta(X)$. Moreover, $c_\delta\neq 0$ only if $\sum_jd_j=\ell-1$.
\end{thm}

\proof
Set $\rho:=c_1(L)$, and $\rho^*\in H_2(X,\Z)$ be its Poincar\'e dual homology class. In particular, we have $c_1(X)=\ell\rho$ and $c_1(E)=(\sum_{i=1}^sd_i)\rho$. By the adjunction formula, we have $c_1(Y)=\iota^*(c_1(X)-c_1(E))$.
From Lemma \ref{restJ}, we have
\beq
\label{JX}
J_X(\delta+\log z\cdot c_1(X))\rqh=\sum_{d\in\N}J_{d\rho^*}(\delta)z^{d\ell+c_1(X)}=\sum_{d\in\N}J_{d\rho^*}(\delta)z^{d\ell+\ell\rho},
\eneq
where $J_{d\rho^*}(\delta)=e^\delta\sum_{\al,k}\langle\tau_k T_\al,1\rangle^X_{0,2,d\rho^{*}}T^\al$.
Analogously, from \eqref{Ifun} we have
\begin{align}
\nonumber
&I_{X,Y}(\delta+(c_1(X)-c_1(E))\log z)\rqh\\
\nonumber
&=\sum_{d\in\N}J_{d\rho^*}(\delta+(c_1(X)-c_1(E))\log z)\prod_{i=1}^s\prod_{m=1}^{\langle d_i\rho,d\rho^*\rangle}(d_i\rho+m)\\
\nonumber
&=\sum_{d\in\N}J_{d\rho^*}(\delta)z^{d(\ell-\sum d_i)+c_1(X)-c_1(E)}\prod_{i=1}^s\prod_{m=1}^{d\cdot d_i}(d_i\rho+m)\\
\label{IXY}
&=\sum_{d\in\N}J_{d\rho^*}(\delta)z^{d(\ell-\sum d_i)+(\ell-\sum d_i)\rho}\prod_{i=1}^s\frac{\Gamma(1+d_i\rho+dd_i)}{\Gamma(1+d_i\rho)}.
\end{align}
On the one hand, from equation \eqref{JX}, one can see that the function $J_X(\delta+\log z\cdot c_1(X)))\rqh$ is the analytification $\widehat{{\rm J}}_X$ of the series ${\rm J}_X\in\mathscr F_\ell(X)$, concentrated at $c_1(X)=\ell\rho$, defined by 
\[{\rm J}_X(Z)=\sum_{d\in\N}J_{d\rho^*}(\delta)Z^{d\ell\oplus c_1(X)}.
\]
On the other hand, one recognises in equation \eqref{IXY} the analytification of the iteration of Laplace transforms
\beq\label{csq}
{\rm I}_{X,Y}:=\prod_{i=1}^s\frac{1}{\Gamma(1+d_i\rho)}\cdot\left(\mathscr L_{\frac{\ell-\sum_{i=1}^s d_i}{d_s},\frac{d_s}{\ell-\sum_{i=1}^{s-1} d_i}}\circ\dots\circ\mathscr L_{\frac{\ell-d_1-d_2}{d_2},\frac{d_2}{\ell-d_1}}\circ\mathscr L_{\frac{\ell-d_1}{d_1},\frac{d_1}{\ell}}[{\rm J}_X]\right),
\eneq
which is an element of $\mathscr F_{\frac{\ell-\sum_{i=1}^s d_i}{\ell}}(X)$. By Theorems \ref{ql1}, \ref{ql2}, \ref{BLAF}, and Proposition \ref{ql3}, we have
\[J_Y(\iota^*\delta+c_1(Y)\log z)\rqh=\iota^*\widehat{\rm I}_{X,Y}(\delta+(c_1(X)-c_1(E)))\exp(-zH(\delta)|_{\bf Q=1}),
\]where $H(\dl)$ is defined in Proposition \ref{ql3}.
Thus, the components of the r.h.s., with respect to any basis of $H^\bullet(Y,\C)$, span the space of master functions $\mc S_{i^*\delta}(Y)$,  by Corollary \ref{masJ}. The factor $\iota^*\prod_{i=1}^s{\Gamma(1+d_i\rho)}^{-1}$ coming from \eqref{csq} can be eliminated by a change of basis of $H^\bullet(Y,\C)$. By $H^\bullet(X,\C)$-linearity of the Laplace $(\al,\bt)$-transforms, the claim follows by setting $c_\delta:=H(\delta)|_{\bf Q=1}$.
\endproof

\begin{rem}
Integral \eqref{intth1} is convergent for any $z\in\widetilde{\C^*}$. This follows from the exponential asymptotics of Theorem \ref{thfs} for $z\to \infty$, the Fano assumption on $Y$ (i.e. $\sum_{j=1}^sd_j<\ell$), and the asymptotics $|\Phi(z)|<C|\log z|^{\dim_\C X}$ for $z\to 0^+$ (see Theorem \ref{topJ} and Corollary \ref{masJ}).
\end{rem}

\begin{rem}
Formula \eqref{csq} generalizes \cite[Lemma 8.1]{galkin2019}.
\end{rem}

\subsection{Integral representations of the second kind} Let $X_1,\dots, X_h$ be Fano smooth projective varieties. 
Assume that $\det T_{X_j}=L_j^{\otimes \ell_j}$ for ample line bundles $L_j$'s. Let $Y$ be a smooth subvariety of $X:=\prod_{j=1}^hX_j$  defined as the zero locus of a regular section of the line bundle $E=\boxtimes_{j=1}^hL_j^{\otimes d_j}$, where the numbers $d_j\in\N^*$ are such that $d_j<\ell_j$ for any $j=1,\dots, h$.

By K\"unneth isomorphism, any element of $H^2(X,\C)$ is of the form $$\bm\delta=\sum_{i=1}^h1\otimes\dots\otimes \delta_i\otimes \dots\otimes 1,\quad\text{with $\delta_i\in H^2(X_i,\C)$.}$$  
Denote by $\iota\colon Y\to X$ the inclusion.
 
 \begin{thm}\label{TH2}
 Let $\bm \delta\in H^2(X,\C)$, $\delta_i\in H^2(X_i,\C)$ as above, and $\mc S_{\delta_i}(X_i)$ the corresponding space of master functions of $QH^\bullet(X_i)$. There exists a rational number $c_{\bm\delta}\in\Q$ such that the space of master functions $\mc S_{\iota^*\delta}(Y)$ is contained in image of the $\C$-linear map $\mathscr P_{(\bm \ell,\bm d)}\colon \otimes_{j=1}^h\mc S_{\delta_j}(X_j)\to \mc O(\widetilde{\C^*})$ defined by 
\[\mathscr P_{(\ell,\bm d)}[\Phi_1,\dots,\Phi_h](z):=e^{-c_\delta z}\mathscr L_{\bm\al,\bm\bt}[\Phi_1,\dots,\Phi_h](z),
\]
where $(\bm\al,\bm\bt)=\left(\frac{\ell_1-d_1}{d_1},\dots,\frac{\ell_h-d_h}{d_h}; \frac{d_1}{\ell_1},\dots,\frac{d_h}{\ell_h}\right)$.
In other words, any element of $\mc S_{\iota^*\delta}(Y)$ is of the form
\beq\label{intth2}
e^{-c_\delta z}\int_0^\infty\prod_{j=1}^h\Phi_j\left(z^{\frac{\ell_j-d_j}{\ell_j}}\la^\frac{d_j}{\ell_j}\right)e^{-\la}d\la,
\eneq
for some $\Phi_j\in\mc S_{\delta_j}(X)$ with $j=1,\dots,h$. Moreover, $c_\delta\neq 0$ only if $d_j=\ell_j-1$ for some $j$.
 \end{thm}
\proof
Set $\rho_i:=c_1(L_i)$ and let $\rho_i^*\in H_2(X_i,\Z)$ be its Poincar\'e dual homology class, for any $i=1,\dots, h$. By K\"unneth isomorphism, and by the universal property of coproduct of algebras (i.e. tensor product), we have injective\footnote{In particular, we have inclusions $\mathscr F_{\bm k}(X_j)\to \mathscr F_{\bm k}(X)$. } maps $H^\bullet(X_i,\C)\to H^\bullet(X,\C)$. In order to ease the computations, in the next formulas we will not distinguish an element of $H^\bullet(X_i,\C)$ with its image in $H^\bullet(X,\C)$. So, for example we will write $c_1(E)=\sum_{p=1}^hd_p\rho_p$. The same will be applied for elements in $H_2(X,\Z)$.

We have
\begin{align}
\nonumber
J_{X}\left(\bm\delta+c_1(X)\log z\right)\rqh&=\bigotimes_{i=1}^h J_{X_i}(\delta_i+c_1(X_i)\log z)\rqh\\
\label{JX2}
&=\bigotimes_{i=1}^h\sum_{k_i\in \N} J_{i,k_i\rho_i^*}(\delta_i)z^{k_i\ell_i+\ell_i\rho_i},
\end{align}
where $J_{i,k_i\rho^*_i}(\delta_i)=e^{\delta_i}\sum_{\al,j}\langle\tau_j T_{\al,i},1\rangle^{X_i}_{0,2,k_i\rho_i^{*}}T^\al_i$.
Analogously, from \eqref{Ifun}, we deduce the formula
\begin{align}
\nonumber
&I_{X,Y}(\bm\delta+(c_1(X)-c_1(E))\log z)\rqh\\
\nonumber
&=\sum_{k_1,\dots,k_h\in\N}\bigotimes_{i=1}^h J_{i,k_i\rho^*_i}(\delta_i) z^{k_i(\ell_i-d_i)+(\ell_i-d_i)\rho_i} \prod_{m=1}^{\langle \sum_{p}d_p\rho_p, \sum_{p}k_p\rho_p^*\rangle}\left(\sum_{p}d_p\rho_p+m\right)\\
\nonumber
&=\sum_{k_1,\dots,k_h\in\N}\bigotimes_{i=1}^h J_{i,k_i\rho^*_i}(\delta_i) z^{k_i(\ell_i-d_i)+(\ell_i-d_i)\rho_i} \prod_{m=1}^{\sum_pd_pk_p}\left(\sum_{p}d_p\rho_p+m\right)\\
\label{IXY2}
&=\sum_{k_1,\dots,k_h\in\N}\bigotimes_{i=1}^h J_{i,k_i\rho^*_i}(\delta_i) z^{k_i(\ell_i-d_i)+(\ell_i-d_i)\rho_i} \frac{\Gamma(1+\sum_{p}d_pk_p+\sum_{p}d_p\rho_p)}{\Gamma(1+\sum_{p}d_p\rho_p)}.
\end{align}
Each element in the tensor product \eqref{JX2} can be easily recognized as the analytification $\widehat{\rm J}_{X_i}$ of a series ${\rm J}_{X_i}\in\mathscr F_{\ell_i}(X)$, for each $i=1,\dots, h$. The function in equation \eqref{IXY2} can be identified with the analytification of the Laplace $(\bm\al,\bm\bt)$-multitransform
\beq
\label{csq2}
{\rm I}_{X,Y}=\left(\bigotimes_{i=1}^h\frac{1}{\Gamma(1+\sum_{p}d_p\rho_p)}\right)\cup_X{\mathscr L_{\al,\bt}[\otimes_{i=1}^h {\rm J}_{X_i}]},\quad 
\eneq
where $(\bm\al,\bm\bt)=\left(\frac{\ell_1-d_1}{d_1},\dots,\frac{\ell_h-d_h}{d_h}; \frac{d_1}{\ell_1},\dots,\frac{d_h}{\ell_h}\right)$. The series ${\rm I}_{X,Y}$ can be seen as an element of $\mathscr F_{\bm \kappa}(X)$, with $\bm \kappa=(\ell_i-d_i)_{i=1}^h$, via the K\"unneth isomorphism. By Theorems \ref{ql1}, \ref{ql2}, \ref{BLAF}, and Proposition \ref{ql3}, we have
\[J_Y(\iota^*\bm\delta+c_1(Y)\log z)\rqh=\iota^*\widehat{\rm I}_{X,Y}(\bm \delta+(c_1(X)-c_1(E)))\exp(-zH(\bm \delta)|_{\bf Q=1}).
\]
Thus, the components of the r.h.s., with respect to any basis of $H^\bullet(Y,\C)$, span the space of master functions $\mc S_{\iota^*\delta}(Y)$,  by Corollary \ref{masJ}. The factor $\iota^*\otimes_{i=1}^s{\Gamma(1+\sum_{p}d_p\rho_p})^{-1}$ coming from \eqref{csq2} can be eliminated by a change of basis of $H^\bullet(Y,\C)$. By $H^\bullet(X,\C)$-linearity of the Laplace $(\bm\al,\bm\bt)$-multitransform, the claim follows by setting $c_{\bm\delta}:=H(\bm \delta)|_{\bf Q=1}$.
\endproof

\begin{rem}
Integral \eqref{intth2} is convergent for any $z\in\widetilde{\C^*}$. This follows from the exponential asymptotics of Theorem \ref{thfs} for $z\to \infty$, the assumption $d_j<\ell_j$ for any $j=1,\dots, h$, and the asymptotics $|\Phi_j(z)|<C|\log z|^{\dim_\C X_j}$ for $z\to 0^+$ (see Theorem \ref{topJ} and Corollary \ref{masJ}).
\end{rem}

\begin{rem}
Formula \eqref{csq2} generalizes \cite[Lemma 8.1]{galkin2019}.
\end{rem}

\subsection{Master functions as Mellin-Barnes integrals} When applied to the case of Fano complete intersections in products of projective spaces, Theorems \ref{TH1} and \ref{TH2} give explicit Mellin-Barnes integral representations of solutions of the $qDE$.

\begin{thm}\label{TH1b}
Let $Y$ be a Fano complete intersection in $\P^{n-1}$ defined by $h$ homogeneous polynomials of degrees $d_1,\dots, d_h$. There exists a unique $c\in\Q$ such that any master functions in $\mc S_0(Y)$ is a linear combination of the Mellin-Barnes integrals
\beq
G_j(z)=\frac{e^{-cz}}{2\pi\sqrt{-1}}\int_\gamma\Gamma(s)^n\prod_{k=1}^h\Gamma\left(1-d_ks\right)z^{-(n-\sum_{k=1}^hd_k)s}\phi_j(s)ds,
\eneq
for $j=0,\dots,n-1$. The path of integration $\gm$ is a parabola of the form ${\rm Re}\,s=-\rho_1({\rm Im}\,s)^2+\rho_2$, for suitable $\rho_1,\rho_2\in\R_+$, such that $\gm$ encircles the poles of $\Gm(s)^n$, and separates them from the poles of the factors $\Gamma\left(1-d_ks\right)$. The functions $\phi_j$ are given by

\noindent{$\bullet$ for $n$ even:}
\beq
\phi_j(s):=\exp\left({2\pi\sqrt{-1}js}\right),\quad j=0,\dots,n-1;
\eneq

\noindent{$\bullet$ for $n$ odd:}
\beq
\phi_j(s):=\exp\left({2\pi\sqrt{-1}js}+\pi\sqrt{-1}s\right),\quad j=0,\dots,n-1.
\eneq
Moreover, $c\neq 0$ only if $\sum_kd_k=n-1$.
\end{thm}
\proof
The functions 
\[g_j(z):=\frac{1}{2\pi\sqrt{-1}}\int_\gamma\Gamma(s)^nz^{-ns}\phi_j(s)ds,\quad j=0,\dots, n-1,
\]
are a basis of the space of master functions $\mc S_0(\P^{n-1})$, see \cite[Lemma 5]{guzzetti1}. The result follows by applying Theorem \ref{TH1} to the case $X=\P^{n-1}$, $\ell=n$.
\endproof

\begin{thm}\label{TH2b}
Let $Y$ be a Fano hypersurface of $\P^{n_1-1}\times\dots\times\P^{n_h-1}$ defined by a homogeneous polynomial of multi-degree $(d_1,\dots,d_h)$. There exists a unique $c\in\Q$ such that any master function in $\mc S_0(Y)$ is a linear combination of the multi-dimensional Mellin-Barnes integrals 
\[
H_{\bm j}(z):=\frac{e^{-cz}}{(2\pi\sqrt{-1})^h}\int_{\bigtimes\gamma_{i}}\left[\prod_{i=1}^h\Gamma(s_i)^{n_i}\phi_{j_i}^{i}(s_i)\right]\Gamma\left(1-\sum_{i=1}^hs_i\right)z^{-\sum_{i=1}^hd_is_i}ds_1\dots ds_h,
\]
for $\bm j=(j_1,\dots, j_h)\in\prod_{i=1}^h\{0,\dots, n_i-1\}$. The paths $\gm_i$ are parabolas of the form ${\rm Re}\,s_i=-\rho_{1,i}({\rm Im}\,s_i)^2+\rho_{2,i}$, for suitable $\rho_{1,i},\rho_{2,i}\in\R_+$, so that they encircle the poles of the factors $\Gamma(s_i)^{n_i}$.
The function $\phi_{j_i}^{i}$ is defined as follows

\noindent$\bullet$ for $n_i$ even:
\[\phi^i_{j_i}(s_i):=\exp\left({2\pi\sqrt{-1}j_is_i}\right),\quad j_i=0,\dots,n_i-1;
\]

\noindent$\bullet$ for $n_i$ odd:
\[\phi^i_{j_i}(s_i):=\exp\left({2\pi\sqrt{-1}j_is_i}+\pi\sqrt{-1}s_i\right),\quad j_i=0,\dots,n_i-1.
\]Moreover, $c\neq 0$ only if $d_i=n_i-1$ for some $i=1,\dots, h$.
\end{thm}

\proof 
The result follows by application of Theorem \ref{TH2} to the case $X_i=\P^{n_i-1}$,  $\ell_i=n_i$. For each factor $\P^{n_i-1}$ a basis of the space $\mc S_0(\P^{n_i-1})$ is given by the integrals
\[\puqed
g_{j_i}^i(z):=\frac{1}{2\pi\sqrt{-1}}\int_{\gamma_i}\Gamma(s)^{n_i}z^{-n_is}\phi_{j_i}^i(s)ds,\quad j_i=0,\dots, n_i-1.\qedhere
\poqed
\]

\begin{example}
Consider the complex Grassmannian $\mathbb G:=\mathbb G(2,4)$: it can be realized as a quadric in $\P^5$, by Pl\"ucker embedding.
It can be shown that the space $\mc S_0(\mathbb G)$ is the space of solutions $\Phi$ of the $qDE$ given by
\beq
\label{qdeg}\vartheta^5\Phi-1024 z^4\vartheta\Phi-2048 z^4\Phi=0,\quad \vartheta:=z\frac{d}{dz}.
\eneq
By Theorem \ref{TH1b}, any solution of \eqref{qdeg} is a linear combination of the functions
\[G_j(z)=\frac{1}{2\pi\sqrt{-1}}\int_\gamma\Gamma(s)^6\Gamma(1-2s)z^{-4s}\exp\left({2\pi\sqrt{-1}js}\right)ds,\quad j=0,\dots,5.
\]Recalling the reflection and duplication formulae for $\Gamma$-function (see e.g. \cite{NIST}),
\[\Gamma(z)\Gamma(1-z)=\frac{\pi}{\sin (\pi z)},\quad \Gamma(2z)=\pi^{-\frac{1}{2}}2^{2z-1}\Gamma(z)\Gamma\left(z+\frac{1}{2}\right),
\]it is easy to see that the function
\[G_0(z)=\frac{2\pi^\frac{3}{2}}{2\pi\sqrt{-1}}\int_\gamma\frac{\Gamma(s)^5}{\Gamma\left(s+\frac{1}{2}\right)}\frac{4^{-s}}{\sin(2\pi s)}z^{-4s}ds
\]
is a solution of \eqref{qdeg}. In \cite[Section 6]{CDG} the solutions 
\[\Phi_1(z)=\frac{1}{2\pi \sqrt{-1}}\int_{\gamma}\frac{\Gamma(s)^5}{\Gamma\left(s+\frac{1}{2}\right)}4^{-s}z^{-4s}ds,
\]
\[\Phi_2(z)=\frac{1}{2\pi \sqrt{-1}}\int_{\gamma}\Gamma(s)^5\Gamma\left(\frac{1}{2}-s\right)e^{i\pi s}4^{-s}z^{-4s}ds,
\]of equation \eqref{qdeg} were found and studied. It is not difficult to see that $\Phi_1$ and $\Phi_2$ are linear combinations of the functions $G_j$'s. 
\end{example}

\begin{rem}
This example can be extended to Grassmannians $\mathbb G(k,n)$ and other families of partial flag varieties. In the case of Grassmannians it gives different integral representations of solutions w.r.t. those obtained from the quantum Satake identification \cite{golymaniv,KS}. More in general, it would be interesting to do a comparison with the integral representations of solutions obtained from the Abelian--Nonabelian correspondence \cite{CFKS}.
\end{rem}

\section{Dubrovin Conjecture}\label{sec8}

\subsection{Exceptional collections and exceptional bases} Let $X$ be a smooth complex projective variety, and denote by $\mc D^b(X)$ the bounded derived category of coherent sheaves on $X$, see \cite{gelman,huy}.  Given $E,F\in{\rm Ob}\left(\mc D^b(X)\right)$, define $\Hom^\bullet(E,F)$ as the $\C$-vector space\footnote{Notice that the category $\mc D^b(X)$ is a $\C$-linear category.}
\[\Hom^\bullet(E,F):=\bigoplus_{k\in\Z}\Hom(E,F[k]).
\]
An object $E\in{\rm Ob}\left(\mc D^b(X)\right)$ is said to be \emph{exceptional} if $\Hom^\bullet(E,E)$ is a one dimensional $\C$-algebra, generated by the identity morphism.

A collection $\frak E=(E_1,\dots, E_n)$ of objects of $\mc D^b(X)$ is said to be an \emph{exceptional collection} if 
\begin{enumerate}
\item each object $E_i$ is exceptional,
\item we have $\Hom^\bullet(E_j,E_i)=0$ for $j>i$.
\end{enumerate}
Moreover, an exceptional collection $\frak E$ is \emph{full} if it generates $\mc D^b(X)$, i.e. any triangular subcategory containing all objects of $\frak E$ is equivalent to $\mc D^b(X)$ via the inclusion functor.

Consider the Grothendieck group $K_0(X)\equiv K_0(\mathcal D^b(X))$, and let $\chi$ to be the Grothen\-dieck-Euler-Poincar\'e bilinear form
\beq
\chi([V],[F]):=\sum_{k}(-1)^k\dim_{\mathbb C}\Hom(V,F[k]),\quad V,F\in\mathcal D^b(X).
\eneq
\begin{defn}
A basis $(e_i)_{i=1}^n$ of $K_0(X)_{\mathbb C}$ is called \emph{exceptional} if $\chi(e_i,e_i)=1$ for $i=1,\dots, n$, and $\chi(e_j,e_i)=0$ for $1\leq i<j\leq n$.
\end{defn}
\begin{lem}
Let $(E_i)_{i=1}^n$ be a full exceptional collection in $\mathcal D^b(X)$. The $K$-classes $([E_i])_{i=1}^n$ form an exceptional basis of $K_0(X)_{\mathbb C}$.\qed
\end{lem}

\subsection{Mutations and helices} Let $\frak E=(E_1,\dots, E_n)$ be an exceptional collection in $\mc D^b(X)$. For any $i=1,\dots, n-1$ define the collections
\begin{align*}
\L_i\frak E:=(E_1,\dots, E_{i-1},E_{i+1}',E_{i},E_{i+2},\dots, E_n),\\
\R_i\frak E=(E_1,\dots, E_{i-1},E_{i+1},E''_{i},E_{i+2},\dots, E_n),
\end{align*}
where the objects $E_{i+1}', E_{i}''$ sit in the distinguished triangles
\[\xymatrix{
E_{i+1}'[-1]\ar[r]&\Hom^\bullet(E_i,E_{i+1})\otimes E_i\ar[r]&E_{i+1}\ar[r]&E_{i+1}'
}\]
\[\xymatrix{
E_{i}''\ar[r]&E_{i}\ar[r]&\Hom^\bullet(E_i,E_{i+1})^*\otimes E_{i+1}\ar[r]&E_{i}''[1].
}
\]
\begin{rem}The object $E_{i+1}'$ (resp. $E_{i}''$) is uniquely defined \emph{up to unique isomorphism}, because of the exceptionality of $E_i$ (resp. $E_{i+1}$), see \cite[Section 3.3]{CDG1}.
\end{rem}
\begin{prop}\cite{BK,helix}
For any $i$, with $0<i<n$, the
 collections $\mathbb L_i\frak E, \mathbb R_i\frak E$ are exceptional.
  The mutation operators $\mathbb L_i,\mathbb R_i$ satisfy the following identities:
\[
\mathbb L_i\mathbb R_i=\mathbb R_i\mathbb L_i={\rm Id},
\]
\[
\mathbb R_i\mathbb R_j=\mathbb R_j\mathbb R_i,\quad \text{if}\ \
|i-j|>1,\quad \mathbb R_{i+1}\mathbb R_i\mathbb R_{i+1}=\mathbb R_i\mathbb R_{i+1}\mathbb R_i.
\]Moreover, if $\frak E$ is full, then also $\mathbb L_i\frak E$ and $\mathbb R_i\frak E$ are full. \qed
\end{prop}

Denote by $\bt_1,\dots, \bt_{n-1}$ the generators of the braid group $\mathcal B_n$, satisfying the relations 
\[\bt_i\bt_{i+1}\bt_i=\bt_{i+1}\bt_i\bt_{i+1},\quad \bt_i\bt_j=\bt_j\bt_i,\quad
\text{if}\ \ |i-j|>1.
\]
 We define the left action of $\mathcal B_n$ on the set of exceptional collections of length $n$ by identifying the action of  $\bt _i$ with $\mathbb L_i$.
 
 \begin{defn}
 Let $\frak E=(E_1,\dots,E_n)$ be a full exceptional collection. We define the \emph{helix} generated by $\frak E$ to be the infinite family $(E_i)_{i\in\Z}$ of exceptional objects such that
 \[(E_{1-kn},E_{2-kn},\dots, E_{n-kn})=\frak E^\bt,\quad \bt=(\bt_{n-1}\dots\bt_1)^{kn},\quad k\in\Z.
 \]Any family of $n$ consecutive exceptional objects $(E_{i+k})_{k=1}^n$ is called a \emph{foundation} of the helix.
 \end{defn}

 \begin{lem}[\cite{helix}]\label{sper}$\quad$\begin{enumerate}
\item Any foundation is a full exceptional collection.
\item For $i,j\in\mathbb Z$, we have $\Hom^\bullet(E_i,E_j)\cong \Hom^\bullet(E_{i-n},E_{j-n})$.\qed
\end{enumerate}
\end{lem}
 
 The action of the braid group on the set of exceptional collections in $\mathcal D^b(X)$ admits a $K$-theoretical analogue on the set of exceptional bases of $K_0(X)_{\mathbb C}$, see \cite{helix,CDG1}.

\subsection{$\Gamma$-classes and graded Chern character}Let $V$ be a complex vector bundle on $X$ of rank $r$, and let $\delta_1,\dots,\delta_r$ be its Chern roots, so that $c_j(V)=s_j(\delta_1,\dots, \delta_r)$, where $s_j$ is the $j$-th elementary symmetric polynomial.
\begin{defn}Let $Q$ be an indeterminate, and $F\in\mathbb C[\![Q]\!]$ be of the form $F(Q)=1+\sum_{n\geq 1}\alpha_nQ^n$. The $F$-\emph{class} of $V$ is the charcateristic class $\widehat F_V\in H^\bullet(X)$ defined by  $\widehat F_V:=\prod_{j=1}^rF(\delta_j).$
\end{defn}
\begin{defn}
The $\Gamma^{\pm}$-\emph{classes} of $V$ are the characteristic classes associated with the Taylor expansions
\beq
\Gamma(1\pm Q)=\exp\left(\mp\gamma Q+\sum_{m=2}^\infty(\mp1)^m\frac{\zeta(m)}{m}Q^m\right)\in\mathbb C[\![Q]\!],
\eneq
where $\gamma$ is the Euler-Mascheroni constant and $\zeta$ is the Riemann zeta function.
\end{defn}
If $V=TX$, then we denote $\widehat{\Gamma}_X^\pm$ its $\Gamma$-classes.
\begin{defn}
The \emph{graded Chern character} of $V$ is the characteristic class ${\rm Ch}(V)\in H^\bullet(X)$ defined by ${\rm Ch}(V):=\sum_{j=1}^r\exp(2\pi\sqrt{-1}\delta_j)$.
\end{defn}
\subsection{Statement of the conjecture}Let $X$ be a Fano variety. In \cite{dubro0} Dubrovin conjectured that many properties of the $qDE$ of $X$, in particular its monodromy, Stokes and central connection matrices, are encoded in the geometry of exceptional collections in $\mathcal D^b(X)$. The following conjecture is a refinement of the original version in \cite{dubro0}. 
\begin{conj}[\cite{CDG1}]\label{conj}
Let $X$ be a smooth Fano variety of Hodge-Tate type.
\begin{enumerate}
\item The quantum cohomology $QH^\bullet(X)$ has semisimple points if and only if there exists a full exceptional collection in $\mathcal D^b(X)$.
\item If $QH^\bullet(X)$ is generically semisimple, for any oriented ray $\ell$ of slope $\phi\in[0,2\pi[$ 
there is a map from the set of $\ell$-chambers to the set of helices with a marked foundation.
\item Let $\Omega_\ell$ be an $\ell$-chamber and $\frak E_\ell=(E_1,\dots, E_n)$ the corresponding exceptional collection (the marked foundation). Denote by $S$ and $C$ Stokes and central connection matrices computed in $\Om_\ell$ w.r.t. a basis $(T_\alpha)_{\alpha=1}^n$ of $H^\bullet(X,\C)$.
\begin{enumerate}
\item The matrix $S$ is the inverse of the Gram matrix of the $\chi$-pairing in $K_0(X)_{\mathbb C}$ wrt the exceptional basis $[\frak E_\ell]$,
\beq
(S^{-1})_{ij}=\chi(E_i,E_j);
\eneq
\item The matrix $C$ coincides with the matrix associated with the $\mathbb C$-linear morphism
\begin{align}
\textnormal{\textcyr{D}}_X^-\colon K_0(X)_\mathbb C\longrightarrow& H^\bullet(X)\\ 
F\xmapsto{\quad\quad}& \frac{(\sqrt{-1})^{\overline{d}}}{(2\pi)^{\frac{d}{2}}}\widehat\Gamma^-_X\exp(-\pi\sqrt{-1}c_1(X)){\rm Ch}(F),
\end{align}
where $d:=\dim_{\mathbb C}X$, and $\overline d$ is the residue class $d\,({\rm mod\,}2)$. The matrix is computed wrt the exceptional basis $[\frak E_\ell]$ and the pre-fixed basis $(T_\alpha)_{\alpha=1}^n$.
\end{enumerate}
\end{enumerate}
\end{conj}

\begin{rem}
If point (3.b) holds true, then automatically also point (3.a) holds true. This follows from the identity \eqref{const2} and Hirzebruch-Riemann-Roch Theorem, see \cite[Corollary 5.8]{CDG1}.
\end{rem}

\begin{rem}
In \cite{bay4}, A.\,Bayer suggested to drop any reference to $X$ being Fano in the formulation of Dubrovin Conjecture. He proved indeed that the semisimplicity of the quantum cohomology preserves under blow-ups at any number of points. It follows that point (1) of Conjecture \ref{conj} (the \emph{qualitative} part) still holds true after blowing up $X$ at an arbitrary number of points, which may yield  a non-Fano variety. To the best of our knowledge, however, there is no non-Fano example for which both the Stokes and central connection matrices have been explicitly computed. In Sections \ref{Df2k} and \ref{Df2k1} we will provide the first example, in the case of Hirzebruch surfaces.
\end{rem}

\begin{rem}
Assume the validity of points (3.a) and (3.b) of Conjecture \ref{conj}. The action of the braid group $\mathcal B_n$ on the Stokes and central connection matrices (Lemma \ref{actbr}) is compatible with the action of $\mathcal B_n$ on the marked foundations attached at each $\ell$-chambers. Different choices of the branch of the $\Psi$-matrix correspond to shifts of objects of the marked foundation. The matrix $M_0^{-1}$ is identified with the canonical operator $\kappa\colon K_0(X)_{\mathbb C}\to K_0(X)_{\mathbb C},\ [F]\mapsto (-1)^{d}[F\otimes \omega_X]$. Equations \eqref{allsc} imply that the connection matrices $C^{(m)}$, with $m\in \mathbb Z$, correspond to the matrices of the morphism \textcyr{D}$^-_X$ wrt the foundations $(\frak E_\ell\otimes \omega_X^{\otimes m})[md]$. The statement  $S^{(m)}=S$ coincides with the $\Hom$-periodicity described in point (2) of Lemma \ref{sper}, see \cite[Theorem 5.9]{CDG1}.
\end{rem}

\begin{rem}
Conjecture \ref{conj} relates two different aspects of the geometry of $X$, namely its \emph{symptectic structure} (Gromov-Witten theory) and its \emph{complex structure} (the derived category $\mathcal D^b(X)$). Heuristically, Conjecture \ref{conj} follows from Homological Mirror Symmetry Conjecture of M.\,Kontsevich, see \cite[Section 5.5]{CDG1}. 
\end{rem}
\begin{rem}
In the paper \cite{KKP} it was underlined the role of $\Gamma$-classes for refining the original version of Dubrovin's conjecture \cite{dubro0}. Subsequently, in \cite{dubro4} and \cite[$\Gamma$-conjecture II]{gamma1} two equivalent versions of point (3.b) above were given. However, in both these versions, different choices of solutions in Levelt form of the $qDE$ at $z=0$ are chosen wrt the natural ones in the theory of Frobenius manifolds, see  \cite[Section 5.6]{CDG1}.
\end{rem}

\begin{rem}
Point (3.b) of Conjecture \ref{conj} allows to identify $K$-classes with solutions of the joint system of equations \eqref{eq1}, \eqref{qde}. Under this identification, Stokes fundamental solutions correspond to exceptional bases of $K$-theory. In the approach of \cite{TV,CV}, where the equivariant case is addressed, such an identification is more fundamental and \emph{a priori}: it is defined via explicit integral representations of solutions of the joint system of $qDE$ and $qKZ$ equations.
\end{rem}

\begin{rem}
The existence of a map between $\ell$-chambers and helices with a marked foundation, discussed in point (2) of Conjecture \ref{conj}, is an important aspect of Dubrovin conjecture. A careful study of such correspondence may hide several delicate open problems. Consider, for instance, the study of injectivity and surjectivity of such a map. This study is closely related (possibly equivalent) to the study of the freeness and transitivity of the braid group action on the set of exceptional collections. These are well-known open problems, whose answer is known in a few special cases only, see \cite{helix}. In the remaining sections of this paper, we will address the study of point (3) of the Conjecture \ref{conj}, but not of point (2).
\end{rem}

\section{Quantum cohomology of Hirzebruch surfaces}\label{sec9}

\subsection{Preliminaries on Hirzebruch surfaces}
Hirzebruch surfaces $\mathbb F_k$, with $k\in\mathbb Z$, are defined as the total space of $\mathbb P^1$-projective bundles on $\mathbb P^1$, namely
\[\mathbb F_k:=\mathbb P\left(\mathcal O\oplus\mathcal O(-k)\right),\quad k\in\mathbb Z,
\]where $\mathcal O(n)$ are line bundles on $\mathbb P^1$.
More explicitly, they can be defined as hypersurfaces in $\mathbb P^2\times \mathbb P^1$ by 
\begin{equation}\label{hirze1}\mathbb F_k:=\left\{([a_0:a_1:a_2],[b_1:b_2])\in \mathbb P^2\times \mathbb P^1\colon a_1b_1^k=a_2b_2^k\right\},\quad k\in\mathbb N.
\end{equation}
Hirzebruch surfaces have the following properties:
\begin{itemize}
\item the surfaces $(\Fb_{2k})_{k\in\N}$ are all diffeomorphic;
\item the surfaces $(\Fb_{2k+1})_{k\in\N}$ are all diffeomorphic;
\item the surfaces $\Fb_n$ and $\Fb_m$ with $n\neq m$ are not biholomorphic;
\item the only Fano Hirzebruch surfaces are $\Fb_0\cong \P^1\times\P^1$ and $\Fb_1\cong {\rm Bl}_{pt}\P^2$;
\item the surfaces $\Fb_n$ and $\Fb_m$ are deformation equivalent if and only if $n$ and $m$ have the same parity.
\end{itemize}
See \cite{hirzsup,beauville}.
\begin{rem}
Let $0\leq m\leq \frac{1}{2}n$. Consider the family $\mc F$ defined by the equation
\[\mc F:=\left\{([a_0:a_1:a_2],[b_1:b_2],t)\in \mathbb P^2\times \mathbb P^1\times\C\colon a_1b_1^n-a_2b_2^n+ta_0b_1^{n-m}b_2^m=0\right\}.
\]The central fiber over $t=0$ is $\Fb_n$. Any non-central fiber over $t\neq 0$ is isomorphic to $\Fb_{n-2m}$. See \cite[Example 2.16]{kodaira}. See also \cite{suwa} and \cite[Example 0.1.10]{namba}.
\end{rem}

\begin{rem}
The only possible complex structures on $\mathbb S^2\times\mathbb S^2$ are the even Hirzebruch surfaces $\Fb_{2k}$, with $k\in\N$, and  the only possible complex structures on the connected sum $\P^2\#\overline{\P^2}$ are the odd Hirzebruch surfaces $\Fb_{2k+1}$, with $k\in\N$, see \cite{qin}.
\end{rem}

\subsection{Classical cohomology of Hirzebruch surfaces}
Using the explicit polynomial description \eqref{hirze1} of the Hirzebruch surfaces, let us define the following subvarieties of $\mathbb F_k$:
\begin{align*}
\Sigma_1^k&:=\left\{a_1=a_2=0\right\},\\
\Sigma_2^k&:=\left\{a_2=b_1=0\right\},\\
\Sigma_3^k&:=\left\{a_1=b_2=0\right\},\\
\Sigma_4^k&:=\left\{a_0=0\right\}.
\end{align*} 
Each of these subvarieties naturally define a cycle in $H_2(\mathbb F_k,\mathbb Z)$. Notice that, under the identification
\[\mathbb F_k\equiv\mathcal O(-k)\cup\infty \text{ section},
\]we can 
\begin{enumerate}
\item identify $\Sigma_1^k$ with the $0$-section of $\mathcal O(-k)$,
\item identify $\Sigma_4^k$ with the $\infty$-section,
\item identify both $\Sigma_2^k$ and $\Sigma_3^k$ with (the compactification of) two fibers of $\mathcal O(-k)$.
\end{enumerate}
Using the original notations of Hirzebruch, we denote by
\begin{itemize}
\item $\tau_k\in H_2(\mathbb F_k,\mathbb C)$ the homology class defined by $\Sigma_1^k$,
\item $\varepsilon_k\in H_2(\mathbb F_k,\mathbb C)$ the homology class defined by $\Sigma_4^k$,
\item $\nu_k\in H_2(\mathbb F_k,\mathbb C)$ the homology class defined by both $\Sigma_2^k$ and $\Sigma_3^k$.
\end{itemize}
As it is easily seen, the three classes $\tau_k,\varepsilon_k,\nu_k$ are not $\mathbb Z$-linearly independent. They are indeed related by the equation
\beq\label{etn}
\varepsilon_k=\tau_k+k\nu_k.
\eneq
Finally, let us also introduce a homogeneous basis $(T_{0,k},T_{1,k},T_{2,k},T_{3,k})$ of the classical cohomology $H^\bullet(\mathbb F_k,\mathbb Z)$, where
$$
T_{0,k}:= 1,\quad T_{1,k}:={\rm PD}(\varepsilon_k),\quad T_{2,k}:={\rm PD}(\nu_k),\quad T_{3,k}:={\rm PD}(\rm pt),
$$
where PD$(\alpha)$ denotes the Poincaré dual class of $\alpha\in H_\bullet(\mathbb F_k,\mathbb Z)$. We denote the corresponding dual coordinates by $(t^{0,k},t^{1,k},t^{2,k},t^{3,k})$.

By Leray-Hirsch Theorem, the classical cohomology algebra is generated by the classes $(T_{1,k},T_{2,k})$. More precisely we have the following result.
\begin{thm}
In the classical cohomology ring $H^\bullet(\mathbb F_k,\mathbb Z)$, the following identities hold true:
\begin{enumerate}
\item $T_{1,k}^2=k\cdot T_{3,k}$,
\item $T_{2,k}^2=0$,
\item $T_{1,k}T_{2,k}=T_{3,k}$.
\end{enumerate}
Hence, the following presentation of algebras holds:
\[H^\bullet(\mathbb F_k,\mathbb C)\cong\frac{\mathbb C[T_{1,k},T_{2,k}]}{\langle T_{2,k}^2,T_{1,k}^2-k\cdot T_{1,k}T_{2,k}\rangle}.
\]
\end{thm}
The Poincaré metric in the basis $(T_{i,k})_{i=0}^3$ is given by
\beq
\label{pmetr}\eta_k=\begin{pmatrix}
0&0&0&1\\
0&k&1&0\\
0&1&0&0\\
1&0&0&0
\end{pmatrix}.
\eneq

\begin{prop}[\cite{KN90}]
Let $k\in\N$. The collection $(\mc O,\mc O(\Si_2^k),\mc O(\Si_4^k),\mc O(\Si_2^k+\Si_4^k))$ is a full exceptional collection in $\mc D^b(\Fb_k)$. The corresponding Gram matrix of the $\chi$-pairing is
\[\begin{pmatrix}
1&2&2+k&4+k\\
0&1&k&2+k\\
0&0&1&2\\
0&0&0&1
\end{pmatrix}.
\]
\end{prop}
\proof
The Gram matrix can be easily computed by Hirzebruch-Riemann-Roch Theorem.
\endproof

\subsection{Quantum cohomology of Hirzebruch surfaces}

There exist only two classes of deformation equivalence of Hirzebruch surfaces, namely $(\Fb_{2k})_{k\in\N}$ and $(\Fb_{2k+1})_{k\in\N}$. Hence, by the Deformation Axiom of Gromov-Witten invariants \cite{cox}, the quantum cohomology algebra of $\Fb_{2k}$ (resp. $\Fb_{2k+1}$) can be identified with the one of $\Fb_0$ (resp. $\Fb_{1}$), as explained in Remark \ref{defocs}. Notice that, the quantum cohomology algebras of $\Fb_0$ and $\Fb_1$ coincide with the corresponding Batyrev rings \cite{batyrev}. This does not hold true for other Hirzebruch surfaces $\Fb_{k}$ with $k\neq 0,1$, being not Fano \cite{spiehirz}. See also \cite{audin} for a presentation of the quantum cohomology algebra of $\Fb_1$.

\subsubsection{Case of $\Fb_{2k}$} The diffeomorphism $\phi_{2k}\colon \Fb_{2k}\to \Fb_{0}$ induces isomorphisms in homology and cohomology. We have $(\phi_{2k})_*(\tau_{2k})=\tau_{0}$ and $(\phi_{2k})_*(\nu_{2k})=\nu_{0}$, so that from equations \eqref{etn} and \eqref{pmetr} we deduce 
\begin{align}
\label{phi2k1}
\phi_{2k}^*(T_{0,0})&=T_{0,2k},\\
\phi_{2k}^*(T_{1,0})&=T_{1,2k}-kT_{2,2k},\\
\phi_{2k}^*(T_{2,0})&=T_{2,2k},\\
\label{phi2k2}
\phi_{2k}^*(T_{3,0})&=T_{3,2k}.
\end{align}
Thus, we can identify the quantum cohomologies $QH^\bullet(\Fb_0)$ and $QH^\bullet(\Fb_{2k})$ via the change of coordinates 
\beq\label{chcord}
t^{0,2k}=t^{0,0},\quad\quad
t^{1,2k}=t^{1,0},\quad\quad
t^{2,2k}=t^{2,0}-k t^{1,0},\quad\quad
t^{3,2k}=t^{3,0}.
\eneq

\begin{thm}
For any $k\geq 0$, the following isomorphism of algebras holds true:
\[QH^\bullet(\mathbb F_{2k})\cong \frac{\mathbb C[T_{1,2k},T_{2,2k},q_1,q_2]}{\langle T_{2,2k}^{\circ 2}-q_1^kq_2, (T_{1,2k}-k\cdot T_{2,2k})^{\circ 2}-q_1\rangle},
\]where $q_1=\exp(t^{1,2k})$ and $q_2=\exp(t^{2,2k})$.
\end{thm}
\proof It follows from the presentation of the quantum cohomology algebra of $QH^\bullet(\Fb_0)$ $\cong QH^\bullet(\P^1)\otimes QH^\bullet(\P^1)$, and formulae \eqref{phi2k1}-\eqref{phi2k2}, \eqref{chcord}. 
\endproof

\begin{lem}
For all $k\geq 0$ we have that
\beq
\label{t1t2q}
T_{1,2k}\circ T_{2,2k}=T_{3,2k}+kq_1^kq_2.
\eneq
\end{lem}
\proof
By homogeneity
Let $\la_{0,2k},\la_{1,2k},\la_{2,2k},\la_{3,2k}$ be the dual basis of $H_\bullet(\Fb_{2k},\C)$ of the basis $(T_{i,2k})_{i=0}^3$. By the Deformation Axiom of Gromov-Witten invariants, for any $r,s\in\N$, we have
\begin{align*}
&\langle T_{1,2k}, T_{2,2k}, T_{3,2k}\rangle_{0,3,r\la_{1,2k}+s\la_{1,2k}}^{\Fb_{2k}}=\langle T_{1,0}+k T_{2,0}, T_{2,0}, T_{3,0}\rangle^{\Fb_0}_{0,3,r\la_{0,0}+s(\la_{1,0}-k\la_{0,0})}\\
&=\langle T_{1,0}, T_{2,0}, T_{3,0}\rangle^{\Fb_0}_{0,3,(r-sk)\la_{0,0}+s\la_{1,0}} +k\langle T_{2,0}, T_{2,0}, T_{3,0}\rangle^{\Fb_0}_{0,3,(r-sk)\la_{0,0}+s\la_{1,0}}\\
&={\langle\si,1,\si\rangle^{\P^1}_{0,3, (r-sk)H}\langle 1,\si,\si\rangle^{\P^1}_{0,3, (r-sk)H}}+ k\langle 1,1,\si\rangle^{\P^1}_{0,3, (r-sk)H}\langle \si,\si,\si\rangle^{\P^1}_{0,3, (r-sk)H}\\
&=k\cdot \delta_{1,2(r-ks)+1}\delta_{3,2s+1}.
\end{align*}
Here we used the class $H\in H_2(\P^1,\Z)$ to be the hyperplane class, and $\si\in H^2(\P^2,\Z)$ to be its dual. This gives the quantum correction in \eqref{t1t2q}. 
\endproof

\subsubsection{Case of $\Fb_{2k+1}$} The diffeomorphism $\phi_{2k+1}\colon \Fb_{2k+1}\to \Fb_1$ induces an isomorphism $\phi_{2k+1}^*$ in cohomology given by
\begin{align}
\label{f11}
\phi_{2k+1}^*(T_{0,1})&=T_{0,2k+1},\\
\phi_{2k+1}^*(T_{1,1})&=T_{1,2k+1}-kT_{2,2k+1},\\
\phi_{2k+1}^*(T_{2,1})&=T_{2,2k+1},\\
\label{f12}
\phi_{2k+1}^*(T_{3,1})&=T_{3,2k+1}.
\end{align}

We can identify the quantum cohomologies $QH^\bullet(\Fb_1)$ and $QH^\bullet(\Fb_{2k+1})$ via the change of coordinates 
\beq\label{chcordubro1}
t^{0,2k+1}=t^{0,1},\quad\quad
t^{1,2k+1}=t^{1,1},\quad\quad
t^{2,2k+1}=t^{2,1}-k t^{1,1},\quad\quad
t^{3,2k+1}=t^{3,1}.
\eneq

\begin{thm}
For any $k\geq 0$, the following isomorphism of algebras holds true:
\beq
QH^\bullet(\Fb_{2k+1})\cong \frac{\C[T_{1,2k+1}, T_{2,2k+1}, q_1, q_2]}{\Big\langle\begin{aligned}
T_{2,2k+1}^{\circ 2}&-(T_{1,2k+1}-(k+1)T_{2,2k+1})q_1^kq_2,\\  (T_{1,2k+1}-k&T_{2,2k+1})\circ (T_{1,2k+1}-(k+1)T_{2,2k+1})-q_1
\end{aligned}\Big\rangle},
\eneq
where $q_1:=\exp(t^{1,2k+1})$ and $q_2:=\exp(t^{2,2k+1})$.
\end{thm}
\proof
The following presentation for $QH^\bullet(\Fb_1)$ holds true:
\beq
QH^\bullet(\Fb_1)\cong \frac{\C[T_{1,1}, T_{2,1}, q_1, q_2]}{\langle T_{2,1}^{\circ 2}-(T_{1,1}-T_{2,1})q_2,\  T_{1,1}^{\circ 2}-T_{1,1}\circ T_{2,1}-q_1\rangle}.
\eneq
The result follows by formulae \eqref{f11}-\eqref{f12} and \eqref{chcordubro1}.
\endproof

\section{Dubrovin Conjecture for Hirzebruch Surfaces $\Fb_{2k}$}\label{Df2k}

\subsection{$\mc A_\La$-stratum and Maxwell stratum of $QH^\bullet(\Fb_{2k})$}
Fix a point $p=t^{1,2k}T_{1,2k}+t^{2,2k}T_{2,2k}$ of the small quantum cohomology of $\Fb_{2k}$. 
The matrix form of the tensor $\mathcal U$ 
is given by
\[\mathcal U(p)=\begin{pmatrix}
0&2q_1+2kq_1^kq_2&2q_1^kq_2&0\\
2&0&0&2q_1^kq_2\\
2-2k&0&0&2q_1-2kq_1^kq_2\\
0&2+2k&2&0
\end{pmatrix}.
\]The canonical coordinates are given by
\begin{align*}
u_1(p)=-2 \left(q_1^{\frac{1}{2}}-q_1^\frac{k}{2}q_2^\frac{1}{2}\right),&\quad u_2(p)=2 \left(q_1^{\frac{1}{2}}-q_1^\frac{k}{2}q_2^\frac{1}{2}\right),\\
u_3(p)=-2 \left(q_1^{\frac{1}{2}}+q_1^\frac{k}{2}q_2^\frac{1}{2}\right),&\quad u_4(p)=2 \left(q_1^{\frac{1}{2}}+q_1^\frac{k}{2}q_2^\frac{1}{2}\right).
\end{align*}

The $\Psi$-matrix at the point $p$ is given by
\begin{align*}&\Psi(p)=\left(
\begin{array}{cccc}
 -\frac{i q_1^{\frac{1}{2} \left(-\frac{k}{2}-\frac{1}{2}\right)}}{2 \sqrt[4]{q_2}} & \frac{i q_1^{\frac{1}{2} \left(-\frac{k}{2}-\frac{1}{2}\right)} \left(\sqrt{q_1}-k q_1^{k/2} \sqrt{q_2}\right)}{2 \sqrt[4]{q_2}} & -\frac{1}{2} i q_1^{\frac{k-1}{4}} \sqrt[4]{q_2} & \frac{1}{2} i q_1^{\frac{k+1}{4}} \sqrt[4]{q_2} \\
 -\frac{i q_1^{\frac{1}{2} \left(-\frac{k}{2}-\frac{1}{2}\right)}}{2 \sqrt[4]{q_2}} & -\frac{i q_1^{\frac{1}{2} \left(-\frac{k}{2}-\frac{1}{2}\right)} \left(\sqrt{q_1}-k q_1^{k/2} \sqrt{q_2}\right)}{2 \sqrt[4]{q_2}} & \frac{1}{2} i q_1^{\frac{k-1}{4}} \sqrt[4]{q_2} & \frac{1}{2} i q_1^{\frac{k+1}{4}} \sqrt[4]{q_2} \\
 \frac{q_1^{\frac{1}{2} \left(-\frac{k}{2}-\frac{1}{2}\right)}}{2 \sqrt[4]{q_2}} & -\frac{q_1^{\frac{1}{2} \left(-\frac{k}{2}-\frac{1}{2}\right)} \left(k \sqrt{q_2} q_1^{k/2}+\sqrt{q_1}\right)}{2 \sqrt[4]{q_2}} & -\frac{1}{2} q_1^{\frac{k-1}{4}} \sqrt[4]{q_2} & \frac{1}{2} q_1^{\frac{k+1}{4}} \sqrt[4]{q_2} \\
 \frac{q_1^{\frac{1}{2} \left(-\frac{k}{2}-\frac{1}{2}\right)}}{2 \sqrt[4]{q_2}} & \frac{q_1^{\frac{1}{2} \left(-\frac{k}{2}-\frac{1}{2}\right)} \left(k \sqrt{q_2} q_1^{k/2}+\sqrt{q_1}\right)}{2 \sqrt[4]{q_2}} & \frac{1}{2} q_1^{\frac{k-1}{4}} \sqrt[4]{q_2} & \frac{1}{2} q_1^{\frac{k+1}{4}} \sqrt[4]{q_2} \\
\end{array}
\right).
\end{align*}

\begin{prop}
The small quantum cohomology of $\Fb_{2k}$ is contained in the $\mc I^0_\La$-stratum of $QH^\bullet(\Fb_{2k})$. Moreover, the point $p$ is in the $\mc A_\La$-stratum of $QH^\bullet(\Fb_{2k})$ if and only if $q_1=q_1^kq_2$.
\end{prop}
\proof
By Theorem \ref{strLa}, the function $\det\La$ takes the form
\[\det\La(z,p)=\frac{z^2}{z^2 A_0(p)+zA_1(p)+A_2(p)},
\]where $A_0,A_1,A_2$ are holomorphic functions on $QH^\bullet(\Fb_{2k})$. If $p$ is a point of the small quantum locus, an explicit computation shows that 
$\det\La(z,p)=-\frac{1}{256} \left(q_1-q_2 q_1^k\right)^{-1}
$, so that $A_1(p)=A_2(p)=0$. The claim follows.
\endproof

\begin{cor}
Along the small quantum locus of $QH^\bullet(\Fb_{2k})$ the $\mc A_\La$-stratum coincides with the Maxwell stratum $\mc M_{\Fb_{2k}}$.
\end{cor}
\proof
If $q_1=q_1^kq_2$, then we have coalescences of canonical coordinates $u_1,u_2,u_3,u_4$. Any point of the small quantum locus, however, is semisimple.
\endproof

\subsection{Small $qDE$ of $\Fb_{2k}$} In the coordinates $(t^{\al,2k})_{\al=0}^3$, the grading tensor $\mu$ has matrix $\mu={\rm diag}(-1,0,0,1)$.
The isomonodromic system \eqref{qdedu} is
 \begin{empheq}[left=\mathcal H_k^{\rm ev}\colon \empheqlbrace]{align*} 
        \frac{\partial \xi_1}{\partial z}&=(2-2 k) \xi _3+2 \xi _2+\frac{1}{z}\xi _1,\\
        \frac{\partial \xi_2}{\partial z}&=(2 k+2) \xi _4+\xi _1 \left(2 k q_1^kq_2 +2 q_1\right),\\
        \frac{\partial \xi_3}{\partial z}&=2 \xi _1 q_1^kq_2 +2 \xi _4,\\
        \frac{\partial \xi_4}{\partial z}&=2 \xi _2 q_1^kq_2 +\xi _3 \left(2 q_1-2 k  q_1^kq_2\right)-\frac{1}{z}\xi _4.
    \end{empheq}
In the complement of the $\mc A_\La$-stratum, it can be reduced to the single equation in $\xi_1$, the \emph{master differential equation} 
\begin{equation}\label{qdiffpari1}z^4\frac{\partial^4\xi_1}{\partial z^4}-z^2 \left[z^2(8 q_1^kq_2 +8 q_1)-1\right]\frac{\partial^2\xi_1}{\partial z^2}-3z\frac{\partial\xi_1}{\partial z}- \left(-16 z^4 \left(q_1-q_1^kq_2 \right){}^2-3\right)\xi_1=0.
\end{equation}
Given a solution $\xi_1(z,t)$ of equation \eqref{qdiffpari1}, we can reconstruct a solution of the system $\mathcal H_k^{\rm ev}$ through the fomulae
\begin{align*}
\xi_2=&-\frac{\left(-4 (k+1) q_2 z^2 q_1^k+4 (k+1) q_1 z^2+k-1\right)}{16 z^3 \left(q_1-q_2 q_1^k\right)}\xi_1\\ &-\frac{\left(4 (3 k-1) q_2 z^2 q_1^k+4 (k-3) q_1 z^2-k+1\right)}{16 z^2 \left(q_1-q_2 q_1^k\right)}\frac{\partial\xi_1}{\partial z}\\&+\frac{(k-1)}{16 \left(q_1-q_2 q_1^k\right)}\frac{\partial^3\xi_1}{\partial z^3},\\
\xi_3=&-\frac{\left(-4 q_2 z^2 q_1^k+4 q_1 z^2+1\right)}{16 z^3 \left(q_1-q_2 q_1^k\right)}\xi_1\\&-\frac{\left(12 q_2 z^2 q_1^k+4 q_1 z^2-1\right) }{16 z^2 \left(q_1-q_2 q_1^k\right)}\frac{\partial\xi_1}{\partial z}\\&+\frac{1}{16 \left(q_1-q_2 q_1^k\right)}\frac{\partial^3\xi_1}{\partial z^3},\\
\xi_4=&-\frac{\left(4 q_2 z^2 q_1^k+4 q_1 z^2-1\right)}{8 z^2}\xi_1-\frac{1}{8 z}\frac{\partial\xi_1}{\partial z}+\frac{1}{8} \frac{\partial^2\xi_1}{\partial z^2}.
\end{align*}
Looking for solution of the form
\[\xi_1(z,t)=z\cdot\Phi(z,t),
\]the equation \eqref{qdiffpari1} can be rewritten as the \emph{(small) quantum differential equation}
\[
z\left(\vartheta^4\Phi-2\vartheta^3\Phi\right)-8 z^3 \left(q_1+q_1^kq_2\right)\left[\vartheta^2\Phi+\vartheta\Phi\right]+16z^5\left(q_1-q_1^kq_2\right)^2\Phi=0,\quad\vartheta:=z\frac{\partial}{\partial z}.
\]

\subsection{Proof for $QH^\bullet(\Fb_{2k})$.} Let us specialize the system $\mathcal H_k^{\rm ev}$ at the point $0\in QH^\bullet(\mathbb F_{2k})$, for which $q_1=q_2=1$:
 \begin{empheq}[left=\mathcal H_k'\colon \empheqlbrace]{align*} 
        \frac{\partial \xi_1}{\partial z}&=(2-2 k) \xi _3+2 \xi _2+\frac{1}{z}\xi _1,\\
        \frac{\partial \xi_2}{\partial z}&=(2 k+2) \xi _4+\xi _1 \left(2 k +2 \right),\\
        \frac{\partial \xi_3}{\partial z}&=2 \xi _1  +2 \xi _4,\\
        \frac{\partial \xi_4}{\partial z}&=2 \xi _2  +\xi _3 \left(2 -2 k  \right)-\frac{1}{z}\xi _4.
    \end{empheq}
The point $p=0$ is in the $\mc A_\La$-stratum of $QH^\bullet(\Fb_{2k})$, and so in  the Maxwell stratum. Hence, the study of monodromy data of the system of differential equations $\mc H'_k$ fits in the analysis developed in \cite{CDG0,CDG}. In particular, the isomonodromy property is justified by \cite[Theorem 4.5]{CDG}. As explained in Remark \ref{defocs}, we can reduce the computation of the monodromy data of the system $\mc H'_k$ to the single case of $\mc H'_0$. The system $\mc H'_0$ can in turn be integrated using solutions of the isomonodromic system of $QH^\bullet(\P^1)$ \cite[Lemma 4.10]{dubro2}.  

\begin{prop}
Let $(\phi_1^{(i)},\phi_2^{(i)})$ with $i=1,2$ be two solutions of the system \eqref{qdedu} for the quantum cohomology of $\P^1$, specialized at $0\in H^2(\P^1,\C)$, i.e.
 \begin{empheq}[left=\empheqlbrace]{align*} 
        \frac{\partial \phi_1}{\partial z}&=2\phi_2+\frac{1}{2z}\phi _1,\\
        \frac{\partial \phi_2}{\partial z}&=2\phi _1-\frac{1}{2z}\phi_2.
          \end{empheq}
Then the tensor product
\[\binom{\phi_1^{(1)}}{\phi_2^{(1)}}\otimes \binom{\phi_1^{(2)}}{\phi_2^{(2)}}=\begin{pmatrix}
\phi_1^{(1)}\cdot \phi_1^{(2)}\\
\phi_1^{(1)}\cdot \phi_2^{(2)}\\
\phi_2^{(1)}\cdot \phi_1^{(2)}\\
\phi_2^{(1)}\cdot \phi_2^{(2)}
\end{pmatrix}
\]is a solution of the system $\mathcal H'_0$.\qed
\end{prop}

\begin{rem}In order to explicitly compute the monodromy data of $\mc H'_{\rm ev}$ one could still develop the study of solutions of the small quantum differential equation, and then reconstruct the Stokes solutions of $\mc H'_k$ doing a similar argument to the one developed in \cite[Section 6]{CDG} for the quantum cohomology of $\mathbb G(2,4)$.
\end{rem}

\begin{thm}\label{DubconjF2k}
The central connection matrix of $QH^\bullet(\mathbb F_{2k})$, computed at the point $0\in QH^\bullet(\mathbb F_{2k})$, w.r.t. an oriented admissible line $\ell$ of slope $\phi\in]\frac{\pi}{2},\frac{3\pi}{2}[$ and for a suitable choice of the determination of the $\Psi$-matrix, is equal to
\[C_k=\left(
\begin{array}{cccc}
 \frac{1}{2 \pi } & \frac{1}{2 \pi } & \frac{1}{2 \pi } & \frac{1}{2 \pi } \\
 -i+\frac{\gamma }{\pi } & -i+\frac{\gamma }{\pi } & \frac{\gamma }{\pi } & \frac{\gamma }{\pi } \\
 -\frac{(k-1) (\gamma -i \pi )}{\pi } & \frac{i \pi  k-\gamma  k+\gamma }{\pi } & -i+\frac{\gamma -\gamma  k}{\pi } & \frac{\gamma -\gamma  k}{\pi } \\
 \frac{2 (\gamma -i \pi )^2}{\pi } & \frac{2 \gamma  (\gamma -i \pi )}{\pi } & \frac{2 \gamma  (\gamma -i \pi )}{\pi } & \frac{2 \gamma ^2}{\pi } \\
\end{array}
\right),
\]and the corresponding Stokes matrix is equal to
\[S=\left(
\begin{array}{cccc}
 1 & -2 & -2 & 4 \\
 0 & 1 & 0 & -2 \\
 0 & 0 & 1 & -2 \\
 0 & 0 & 0 & 1 \\
\end{array}
\right).
\]The matrix $C_k$ is the matrix associated with the morphism
\[\textnormal{\textcyr{D}}^-_{\mathbb F_{2k}}\colon K_0(\mathbb F_{2k})_\mathbb C\to H^\bullet(\mathbb F_{2k},\mathbb C)\colon [\mathscr F]\mapsto \frac{1}{2\pi}\widehat{\Gamma}^-_{\mathbb F_{2k}}\cup e^{-\pi i c_1(\mathbb F_{2k})}\cup {\rm Ch}(\mathscr F),
\]w.r.t. an exceptional basis $\frak E:=(E_i)_{i=1}^4$ of $K_0(\mathbb F_{2k})_{\mathbb C}$ and the basis $(T_{i,2k})_{i=0}^3$ of $H^\bullet(\mathbb F_{2k},\mathbb C)$. The exceptional basis $\frak E$  is the one obtained by acting on the exceptional basis
\[\left([\mathcal O],[\mathcal O(\Sigma_2^{2k})],[\mathcal O(\Sigma_4^{2k})],[\mathcal O(\Sigma_2^{2k}+\Sigma_4^{2k})]\right),
\]with the element $(J_k^{-1},b_k)\in(\mathbb Z/2\mathbb Z)^4\rtimes \mathcal B_4$, where
\begin{multicols}{2}
 \begin{empheq}[left={J_k:=} \empheqlbrace]{align*} 
 &(1,1,(-1)^{p+1},(-1)^{p}),\quad\text{if }k=2p+1,\\
 \\
&(1,1,(-1)^p,(-1)^p),\quad\text{if }k=2p,
\end{empheq}\break
\begin{align*}
b_k:&=\bt_3^k.\\
& 
\end{align*}
\end{multicols}

\end{thm}

\proof
{\bf Step 1:} Let us show that for suitable choices of oriented line $\ell$ and $\Psi$-matrix, the central connection matrix computed at $0\in QH^\bullet(\mathbb F_0)$ is
\beq\label{c0F0}C_0:=\left(
\begin{array}{cccc}
 \frac{1}{2 \pi } & \frac{1}{2 \pi } & \frac{1}{2 \pi } & \frac{1}{2 \pi } \\
 -i+\frac{\gamma }{\pi } & -i+\frac{\gamma }{\pi } & \frac{\gamma }{\pi } & \frac{\gamma }{\pi } \\
 -i+\frac{\gamma }{\pi } & \frac{\gamma }{\pi } & -i+\frac{\gamma }{\pi } & \frac{\gamma }{\pi } \\
 \frac{2 (\gamma -i \pi )^2}{\pi } & \frac{2 \gamma  (\gamma -i \pi )}{\pi } & \frac{2 \gamma  (\gamma -i \pi )}{\pi } & \frac{2 \gamma ^2}{\pi } \\
\end{array}
\right).
\eneq According to \cite[Corollary 6.11]{CDG1}, the central connection matrix $C$ of $QH^\bullet(\mathbb P^1)$ computed at the point $0$, w.r.t. an oriented line $\ell$ of slope $\phi\in]\frac{\pi}{2},\frac{3\pi}{2}[$ and w.r.t. the following choice of $\Psi$-matrix
\[\Psi_0=\left(
\begin{array}{cc}
 \frac{1}{\sqrt{2}} & \frac{1}{\sqrt{2}} \\
 \frac{i}{\sqrt{2}} & -\frac{i}{\sqrt{2}} \\
\end{array}
\right),
\]equals 
\[C:=\frac{i }{\sqrt{2 \pi }}\left(
\begin{array}{cc}
 1 & 1 \\
 2 (\gamma -\pi  i) & 2 \gamma  \\
\end{array}
\right).
\]
This is the matrix associated with the morphism
\[\textcyr{D}^-_{\mathbb P^1}\colon K_0(\mathbb P^1)_{\mathbb C}\to H^\bullet(\mathbb P^1,\mathbb C)\colon [\mathscr F]\mapsto \frac{i}{(2\pi)^\frac{1}{2}}\widehat{\Gamma}^-_{\mathbb P^1}\cup e^{-\pi i c_1(\mathbb P^1)}\cup {\rm Ch}(\mathscr F),
\]w.r.t. the bases
\begin{itemize}
\item $([\mathcal O],[\mathcal O(1)])$ of $K_0(\mathbb P^1)_{\mathbb C}$ (the Beilinson basis),
\item $(1,\sigma)$ of $H^\bullet(\mathbb P^1,\mathbb C)$.
\end{itemize}
By taking the Kronecker tensor square $C^{\otimes 2}$, we obtain the central connection matrix of $QH^\bullet(\mathbb P^1\times\mathbb P^1)$ computed at the point $0$, w.r.t. the same line $\ell$ (which is still admissible) and w.r.t. the choice of the $\Psi$-matrix given by the Kronecker  tensor square $\Psi_0^{\otimes 2}$:
\[C^{\otimes 2}=\left(
\begin{array}{cccc}
 -\frac{1}{2 \pi } & -\frac{1}{2 \pi } & -\frac{1}{2 \pi } & -\frac{1}{2 \pi } \\
 -\frac{\gamma -i \pi }{\pi } & -\frac{\gamma }{\pi } & -\frac{\gamma -i \pi }{\pi } & -\frac{\gamma }{\pi } \\
 -\frac{\gamma -i \pi }{\pi } & -\frac{\gamma -i \pi }{\pi } & -\frac{\gamma }{\pi } & -\frac{\gamma }{\pi } \\
 -\frac{2 (\gamma -i \pi )^2}{\pi } & -\frac{2 \gamma  (\gamma -i \pi )}{\pi } & -\frac{2 \gamma  (\gamma -i \pi )}{\pi } & -\frac{2 \gamma ^2}{\pi } \\
\end{array}
\right).
\]By changing all the signs of the rows of $\Psi_0^{\otimes 2}$, i.e. acting with $(-1,-1,-1,-1)\in(\mathbb Z/2\mathbb Z)^4$ on $C^{\otimes 2}$, we obtain the matrix $-C^{\otimes 2}$ associated with the morphism
\[\textcyr{D}^-_{\mathbb P^1\times\mathbb P^1}\colon K_0(\mathbb P^1\times \mathbb P^1)_{\mathbb C}\to H^\bullet(\mathbb P^1\times \mathbb P^1,\mathbb C)\colon [\mathscr F]\mapsto \frac{1}{2\pi}\widehat{\Gamma}^-_{\mathbb P^1\times \mathbb P^1}\cup e^{-\pi i c_1(\mathbb P^1\times\mathbb P^1)}\cup {\rm Ch}(\mathscr F),
\]written w.r.t. the bases
\begin{itemize}
\item $([\mathcal O],[\mathcal O(1,0)],[\mathcal O(0,1)],[\mathcal O(1,1)])$ of $K_0(\mathbb P^1\times\mathbb P^1)_{\mathbb C}$,
\item $(1,\sigma\otimes 1,1\otimes\sigma,\sigma\otimes\sigma)$ of $H^\bullet(\mathbb P^1\times \mathbb P^1,\mathbb C)\cong H^\bullet(\mathbb P^1,\mathbb C)^{\otimes 2}$.
\end{itemize}
See \cite[Proposition 5.11]{CDG1}.
In the notations introduced before for Hirzebruch surfaces, this exceptional collection is 
$$\left(\mathcal O,\mathcal O(\Sigma_4^{0}),\mathcal O(\Sigma_2^{0}),\mathcal O(\Sigma_2^{0}+\Sigma_4^{0})\right).$$ 
It is a 3-block exceptional collection\footnote{An exceptional collection $(E_1,\dots, E_n)$ is a $k$-\emph{block exceptional collection} if it is possible to decompose it into $k$ exceptional sub-collections $\frak B_1,\dots, \frak B_k$, called {\it blocks}, such that
\begin{itemize}
\item they are consecutive, i.e. of the form $\frak B_1=(E_1,\dots, E_{j_1}), \frak B_2=(E_{j_1+1},\dots, E_{j_2}),\dots, \frak B_k=(E_{j_{k-1}+1},\dots, E_{j_k})$, with $1\leq j_1<j_2<\dots<j_k\leq n$, 
\item we have $\Hom^\bullet(E_j,E_i)=0$ if $E_i$ and $E_j$ belong to a same block $\frak B_h$.
\end{itemize}
In particular, inside each block $\frak B_h$, mutations act as permutations of exceptional objects. See \cite[Section 3.6.4]{CDG1}, and references therein.}, coherently with the fact that $0\in QH^\bullet(\mathbb F_0)$ is a semisimple coalescing point, see \cite[Section 6]{CDG} and \cite[Remark 5.4]{CDG1}. In particular, the braids $\beta_{2,3}$ and $\beta_{2,3}^{-1}$ act as a mere permutation of the central objects, and of the two central columns of the matrix $-C^{\otimes 2}$. Such a permuted matrix is exactly the matrix $C_0$ in \eqref{c0F0}, and it corresponds to the matrix associated with the morphism $\textcyr{D}^-_{\mathbb F_0}$ w.r.t. the collection
\[\left(\mathcal O,\mathcal O(\Sigma_2^{0}),\mathcal O(\Sigma_4^{0}),\mathcal O(\Sigma_2^{0}+\Sigma_4^{0})\right).
\]
{In conclusion, we have proved that, for suitable choices of $\ell$ and $\Psi$, the central connection matrix computed at $0\in QH^\bullet(\mathbb F_0)$ is
\[C_0=\left(
\begin{array}{cccc}
 \frac{1}{2 \pi } & \frac{1}{2 \pi } & \frac{1}{2 \pi } & \frac{1}{2 \pi } \\
 -i+\frac{\gamma }{\pi } & -i+\frac{\gamma }{\pi } & \frac{\gamma }{\pi } & \frac{\gamma }{\pi } \\
 -i+\frac{\gamma }{\pi } & \frac{\gamma }{\pi } & -i+\frac{\gamma }{\pi } & \frac{\gamma }{\pi } \\
 \frac{2 (\gamma -i \pi )^2}{\pi } & \frac{2 \gamma  (\gamma -i \pi )}{\pi } & \frac{2 \gamma  (\gamma -i \pi )}{\pi } & \frac{2 \gamma ^2}{\pi } \\
\end{array}
\right),
\]which coincides with the matrix associated with the collection
\[\left(\mathcal O,\mathcal O(\Sigma_2^{0}),\mathcal O(\Sigma_4^{0}),\mathcal O(\Sigma_2^{0}+\Sigma_4^{0})\right).
\]
{\bf Step 2:} Equations \eqref{chcord} and Proposition \ref{flatcord} imply that the central connection matrix computed at $0\in QH^\bullet(\mathbb F_{2k})$, w.r.t. the same choices of $\ell$ and $\Psi$, equals
\[C_k=\left(
\begin{array}{cccc}
 \frac{1}{2 \pi } & \frac{1}{2 \pi } & \frac{1}{2 \pi } & \frac{1}{2 \pi } \\
 -i+\frac{\gamma }{\pi } & -i+\frac{\gamma }{\pi } & \frac{\gamma }{\pi } & \frac{\gamma }{\pi } \\
 -\frac{(k-1) (\gamma -i \pi )}{\pi } & \frac{i \pi  k-\gamma  k+\gamma }{\pi } & -i+\frac{\gamma -\gamma  k}{\pi } & \frac{\gamma -\gamma  k}{\pi } \\
 \frac{2 (\gamma -i \pi )^2}{\pi } & \frac{2 \gamma  (\gamma -i \pi )}{\pi } & \frac{2 \gamma  (\gamma -i \pi )}{\pi } & \frac{2 \gamma ^2}{\pi } \\
\end{array}
\right).
\]The corresponding Stokes matrix is independent of $k$, and it is equal to
\begin{equation}\label{stokesf2k}S=\left(
\begin{array}{cccc}
 1 & -2 & -2 & 4 \\
 0 & 1 & 0 & -2 \\
 0 & 0 & 1 & -2 \\
 0 & 0 & 0 & 1 \\
\end{array}
\right).
\end{equation}
{\bf Step 3:} Let us define the matrix $J_k\in (\mathbb Z/2\mathbb Z)^4$ as follows:
 \begin{empheq}[left={J_k:=} \empheqlbrace]{align*} 
 &(1,1,(-1)^{p+1},(-1)^{p}),\quad\text{if }k=2p+1,\\
 \\
&(1,1,(-1)^p,(-1)^p),\quad\text{if }k=2p.
\end{empheq}
We claim that by acting on $C_kJ_k$ with the braid $\bt_3^{-k}$ we obtain the matrix associated with $\textcyr{D}^-_{\mathbb F_{2k}}$ and w.r.t. the exceptional collection
\[\left(\mathcal O,\mathcal O(\Sigma_2^{2k}),\mathcal O(\Sigma_4^{2k}),\mathcal O(\Sigma_2^{2k}+\Sigma_4^{2k})\right),
\]namely the matrix
\[E_k:=\left(
\begin{array}{cccc}
 \frac{1}{2 \pi } & \frac{1}{2 \pi } & \frac{1}{2 \pi } & \frac{1}{2 \pi } \\
 -i+\frac{\gamma }{\pi } & -i+\frac{\gamma }{\pi } & \frac{\gamma }{\pi } & \frac{\gamma }{\pi } \\
 -\frac{(k-1) (\gamma -i \pi )}{\pi } & \frac{i \pi  k-\gamma  k+\gamma }{\pi } & -\frac{(k-1) (\gamma -i \pi )}{\pi } & \frac{i \pi  k-\gamma  k+\gamma }{\pi } \\
 \frac{2 (\gamma -i \pi )^2}{\pi } & \frac{2 \gamma  (\gamma -i \pi )}{\pi } & \frac{2 \gamma  (i \pi  (k-1)+\gamma )}{\pi } & \frac{2 \gamma  (i \pi  k+\gamma )}{\pi } \\
\end{array}
\right).
\]Notice that the claim is equivalent to the following statement: the matrix $A^\beta(J_k\cdot S\cdot J_k)$, with $\beta=\bt_3^{-k}$ and $S$ as in \eqref{stokesf2k}, is equal to
\begin{equation}\label{treccef2k}E_k^{-1}C_kJ_k=\left(
\begin{array}{cccc}
 1 & 0 & 0 & 0 \\
 0 & 1 & 0 & 0 \\
 0 & 0 & k+1 & k \\
 0 & 0 & -k & 1-k \\
\end{array}
\right)\cdot J_k.
\end{equation}

Given a generic $4\times 4$ unipotent upper triangular matrix $X$, the action of subsequent powers of the braid $\bt_3$, or of its inverse $\bt_3^{-1}$, simply changes the sign of the entry  in position $(3,4)$: more precisely,we have that
\[[X^{\beta}]_{3,4}=(-1)^n[X]_{3,4},\quad \text{if }\beta=\bt_3^{\pm n}.
\]
For example, by acting twice with the braid $\bt_3$ we have
\[\left(
\begin{array}{cccc}
 1 & a & b & c \\
 0 & 1 & d & e \\
 0 & 0 & 1 & f \\
 0 & 0 & 0 & 1 \\
\end{array}
\right)\mapsto \left(
\begin{array}{cccc}
 1 & a & c & b-c f \\
 0 & 1 & e & d-e f \\
 0 & 0 & 1 & -f \\
 0 & 0 & 0 & 1 \\
\end{array}
\right)\mapsto \left(
\begin{array}{cccc}
 1 & a &b-c f & c+f (b-c f) \\
 0 & 1 & d-e f & e+f (d-e f) \\
 0 & 0 & 1 & f \\
 0 & 0 & 0 & 1 \\
\end{array}
\right).
\]In particular, the matrix $A^\beta(X)$, with $\beta=\bt_3^{-k}$, is equal to
\[\prod_{j=1}^k\left(
\begin{array}{cccc}
 1 & 0 & 0 & 0 \\
 0 & 1 & 0 & 0 \\
 0 & 0 & (-1)^{j}x & 1 \\
 0 & 0 & 1 & 0 \\
\end{array}
\right),\quad x=X_{3,4}.
\]
In the case $X=J_k\cdot S\cdot J_k$, we have 
\[x=(-1)^{k+1} 2.
\]
So, in conclusion, we have to prove that the following identity holds for all $k\geq 0$:
\[\prod_{j=1}^k\left(
\begin{array}{cccc}
 1 & 0 & 0 & 0 \\
 0 & 1 & 0 & 0 \\
 0 & 0 & (-1)^{j+k+1}2 & 1 \\
 0 & 0 & 1 & 0 \\
\end{array}
\right)=\left(
\begin{array}{cccc}
 1 & 0 & 0 & 0 \\
 0 & 1 & 0 & 0 \\
 0 & 0 & k+1 & k \\
 0 & 0 & -k & 1-k \\
\end{array}
\right)\cdot J_k.
\]
We prove the claim by induction on $k$. The base case $k=0$ is evidently true. Let us assume that the statement holds true for $k-1$, and let us prove it for $k$. We have that
\begin{align*}
\prod_{j=1}^k\left(
\begin{array}{cccc}
 1 & 0 & 0 & 0 \\
 0 & 1 & 0 & 0 \\
 0 & 0 & (-1)^{j+k+1}2 & 1 \\
 0 & 0 & 1 & 0 \\
\end{array}
\right)=&\left[\prod_{j=1}^{k-1}\left(
\begin{array}{cccc}
 1 & 0 & 0 & 0 \\
 0 & 1 & 0 & 0 \\
 0 & 0 & (-1)^{j+k+1}2 & 1 \\
 0 & 0 & 1 & 0 \\
\end{array}
\right)\right]\cdot\left(
\begin{array}{cccc}
 1 & 0 & 0 & 0 \\
 0 & 1 & 0 & 0 \\
 0 & 0 & -2 & 1 \\
 0 & 0 & 1 & 0 \\
\end{array}
\right)\\
=&\left(
\begin{array}{cccc}
 1 & 0 & 0 & 0 \\
 0 & 1 & 0 & 0 \\
 0 & 0 & k & k-1 \\
 0 & 0 & 1-k & 2-k \\
\end{array}
\right)\cdot J_{k-1}\cdot \left(
\begin{array}{cccc}
 1 & 0 & 0 & 0 \\
 0 & 1 & 0 & 0 \\
 0 & 0 & -2 & 1 \\
 0 & 0 & 1 & 0 \\
\end{array}
\right),
\end{align*}
and in both cases $k$ even/odd, the last term is easily seen to be equal to \eqref{treccef2k}.
\endproof


\section{Dubrovin conjecture for Hirzebruch surfaces $\Fb_{2k+1}$}\label{Df2k1}

\subsection{$\mc A_\La$-stratum and Maxwell stratum of $QH^\bullet(\Fb_{2k+1})$}Fix a point $$p=t^{1,2k+1}T_{1,2k+1}+t^{2,2k+1}T_{2,2k+1}$$ of the small quantum cohomology of $\Fb_{2k+1}$. The matrix associated to the $\bm{\mc U}$-tensor at $p$ is 
\[
\mc U(p)=\left(
\begin{array}{cccc}
 0 & 2 q_1 & 0 & 3 q_1^{k+1} q_2 \\
 2 & k q_1^k q_2 & q_1^k q_2 & 0 \\
 1-2 k & k \left(-k q_2 q_1^k-q_1^k q_2\right) & -k q_2 q_1^k-q_1^k q_2 & 2 q_1 \\
 0 & 2 k+3 & 2 & 0 \\
\end{array}
\right).
\]
The canonical coordinates are the roots $u_1(p),u_2(p),u_3(p),u_4(p)$ of the polynomial
\beq
j(u):=u^4+ u^3 q_1^kq_2-8 q_1 u^2-36  u q_1^{k+1}q_2-27 q_2^2 q_1^{2 k+1}+16 q_1^2.
\eneq
Hence the bifurcation set $\mc B_{\Fb_{2k+1}}$, along the small quantum cohomology, is defined by the zero locus of the discriminant of $j(u)$, i.e.
\beq
\mc B_{\Fb_{2k+1}}=\left\{ p\colon  q_1^{2 k+2}q_2^2 \left(27 q_2^2 q_1^{2 k}+256 q_1\right)^3=0\right\}.
\eneq
Since any point of the small quantum cohomology of $\Fb_{2k+1}$ is semisimple, the set above actually coincides with the Maxwell stratum $\mc M_{\Fb_{2k+1}}$.
The determinant of the $\La$-matrix is given by
\beq
\det \La(z,p)=-\frac{z}{(27   q_1^{2 k}q_2^2+256 q_1) z-24 q_2 q_1^k}.
\eneq
Hence, the $\mc A_\La$-stratum is given by
\beq
\label{Alf1}
\mc A_\La:=\left\{p\colon 27   q_1^{2 k}q_2^2+256 q_1=0\right\}.
\eneq
Also in this case, the Maxwell stratum and the $\mc A_\La$-stratum coincide along the small quantum cohomology of $\Fb_{2k+1}$.

\subsection{Small $qDE$ of $\Fb_{1}$} At the point $p$, the grading operator $\mu$ has matrix $\mu={\rm diag}(-1,0,0,1)$. Hence the isomonodromic system of differential equations \eqref{qdedu} for $QH^\bullet(\Fb_{2k+1})$ is given by
\begin{empheq}[left=\mathcal H^{\rm od}_k\colon \empheqlbrace]{align*}
	\frac{\partial \xi_1}{\partial z}&=(1-2 k) \xi _3+2 \xi _2+\frac{\xi _1}{z},\\
        \frac{\partial \xi_2}{\partial z}&=(2 k+3) \xi _4+k \xi _2 q_2 q_1^k+k \xi _3 \left(-k q_2 q_1^k-q_2 q_1^k\right)+2 \xi _1 q_1,\\
        \frac{\partial \xi_3}{\partial z}&=\xi _2 q_2 q_1^k+\xi _3 \left(-k q_2 q_1^k-q_2 q_1^k\right)+2 \xi _4,\\
        \frac{\partial \xi_4}{\partial z}&=3 \xi _1 q_2 q_1^{k+1}+2 \xi _3 q_1-\frac{\xi _4}{z}.
\end{empheq}
As explained in Remark \ref{defocs}, the computation of the monodromy data of $\mc H^{\rm od}_k$ can be reduced to the single case $\mc H^{\rm od}_{0}$. 

The point $0\in QH^\bullet(\Fb_1)$ is not in the $\mc A_\La$-stratum, as it follows from equation \eqref{Alf1}. At the point $0\in QH^\bullet(\Fb_1)$, indeed, the system $\mc H^{\rm od}_0$ can be reduced to the \emph{small quantum differential equation}
\begin{align}\label{qdiffdisp2bis}
&(283 z-24)\vartheta^4\Phi+\left(283 z^2-590 z+24\right)\vartheta^3\Phi+ \left(-2264 z^2+192 z+3\right)\vartheta^2\Phi\\
\nonumber&-4 z^2 \left(2547 z^2+350 z-104\right)\vartheta\Phi+z^2 \left(-3113 z^3-9924 z^2+1476 z+192\right)\Phi=0.
\end{align}
Given a solution $\Phi(z)$ of \eqref{qdiffdisp2bis}, the corresponding solution of the system $\mc H^{\rm od}_0$ can be reconstructed by the formulae
\begin{align}
\label{rec1}
\xi_1(z)&=z\cdot \Phi(z),\\
\nonumber
\xi_2(z)&=\frac{1}{z^2 (283 z-24)}\left(169 z^3 \xi _1'(z)+z^3 \xi _1''(z)+204 z^3 \xi _1(z)-8 z^3 \xi _1{}^{(3)}(z)-9 z^2 \xi _1'(z)\right.\\
\label{rec2}
&\left.-105 z^2 \xi _1(z)-8 z \xi _1'(z)+9 z \xi _1(z)+8 \xi _1(z)\right),\\
\nonumber
\xi_3(z)&=\frac{1}{z^2 (283 z-24)}\left(-55 z^3 \xi _1'(z)-2 z^3 \xi _1''(z)-408 z^3 \xi _1(z)+16 z^3 \xi _1{}^{(3)}(z)-6 z^2 \xi _1'(z)\right.\\
\label{rec3}
&\left.-73 z^2 \xi _1(z)+16 z \xi _1'(z)+6 z \xi _1(z)-16 \xi _1(z)\right),\\
\nonumber
\xi_4(z)&=\frac{1}{z^2 (283 z-24)}\left(-28 z^3 \xi _1'(z)+35 z^3 \xi _1''(z)-218 z^3 \xi _1(z)+3 z^3 \xi _1{}^{(3)}(z)-35 z^2 \xi _1'(z)\right.\\
\label{rec4}
&\left.-3 z^2 \xi _1''(z)+16 z^2 \xi _1(z)+6 z \xi _1'(z)+35 z \xi _1(z)-6 \xi _1(z)\right).
\end{align}

These formulae are obtained by the identity
\[\xi=\Lambda^T\begin{pmatrix}
\xi_1\\
\xi_1'\\
\xi_1''\\
\xi_1^{(3)}
\end{pmatrix},
\]where the $\La$-matrix at $0\in QH^\bullet(\Fb_1)$ is 
\[\La(z,0)=\left(
\begin{array}{cccc}
 1 & \frac{204 z^3-105 z^2+9 z+8}{z^2 (283 z-24)} & \frac{-408 z^3-73 z^2+6 z-16}{z^2 (283 z-24)} & \frac{-218 z^3+16 z^2+35 z-6}{z^2 (283 z-24)} \\
 0 & \frac{169 z^2-9 z-8}{z (283 z-24)} & \frac{-55 z^2-6 z+16}{z (283 z-24)} & \frac{-28 z^2-35 z+6}{z (283 z-24)} \\
 0 & \frac{z}{283 z-24} & -\frac{2 z}{283 z-24} & \frac{35 z-3}{283 z-24} \\
 0 & -\frac{8 z}{283 z-24} & \frac{16 z}{283 z-24} & \frac{3 z}{283 z-24} \\
\end{array}
\right).
\]
\begin{rem}The quantum differential equation \eqref{qdiffdisp2bis} has one apparent singularity at $z=\frac{24}{283}$. This coincides with the zero of the denominator of the determinant of the $\La$-matrix:
\[\det \La(z,0)=\frac{z}{24-283 z}.
\]
\end{rem}
The $\Psi$-matrix at the point $0\in QH^\bullet(\Fb_1)$ is given by
\beq
\Psi=\begin{pmatrix}
\al_1^\frac{1}{2}\eps_1&\al_1^\frac{1}{2}\dl_1&\al_1^\frac{1}{2}\si_1&\al_1^\frac{1}{2}\ups_1\\
\al_2^\frac{1}{2}\eps_2&\al_2^\frac{1}{2}\dl_2&\al_2^\frac{1}{2}\si_2&\al_2^\frac{1}{2}\ups_2\\
\al_3^\frac{1}{2}\eps_3&\al_3^\frac{1}{2}\dl_3&\al_3^\frac{1}{2}\si_3&\al_3^\frac{1}{2}\ups_3\\
\al_4^\frac{1}{2}\eps_4&\al_4^\frac{1}{2}\dl_4&\al_4^\frac{1}{2}\si_4&\al_4^\frac{1}{2}\ups_4
\end{pmatrix},
\eneq
where the numbers $\al_i,\eps_i,\dl_i,\si_i,\ups_i$ satisfy the algebraic equations
\begin{align*}
\al_i^4+\al_i^3-6\al_i^2-283&=0,\\
283 \eps_i^4+6\eps_i^2-\eps_i-1&=0,\\
283\dl_i^4-2\dl_i^2-9\dl_i-1&=0,\\
283\si_i^4-32\si_i^2-\si_i+1&=0,\\
283\ups_i^4-283\ups_i^3+105\ups_i^2-17\ups_i+1&=0.
\end{align*}
Their numerical approximations are
\begin{align*}
\al_1&\approx 4.21193,& \eps_1&\approx 0.237421,\\
\al_2&\approx -0.204399-3.73457 i,& \eps_2&\approx -0.0146116 + 0.266969 i,\\
\al_3&\approx -0.204399+3.73457 i,& \eps_3&\approx -0.0146116 - 0.266969 i,\\
\al_4&\approx -4.80313,& \eps_4&\approx -0.208197,\\
\dl_1&\approx 0.353808,&\si_1&\approx 0.194489,\\
\dl_2&\approx -0.122264 - 0.276482 i,&\si_2&\approx -0.240929 - 0.0719476 i,\\
\dl_3&\approx -0.122264 + 0.276482 i,&\si_3&\approx -0.240929 + 0.0719476 i,\\
\dl_4&\approx -0.10928,&\si_4&\approx 0.28737,\\
\ups_1&\approx 0.28983,\\
\ups_2&\approx 0.279666 - 0.0511337 i,\\
\ups_3&\approx 0.279666 + 0.0511337 i,\\
\ups_4&\approx 0.150837.\\
\end{align*}
The reader can check that $\Psi^T\Psi=\eta$, and that
\beq
\Psi\mc U\Psi^{-1}={\rm diag}(x_1,x_2,x_3,x_4),
\eneq
where the canonical coordinates $x_i$'s are the roots of the polynomial 
\beq
x^4+x^3-8x^2-36x-11=0.
\eneq
Their numerical approximations are
\begin{align}
x_1&\approx 3.7996,\\
x_2&\approx -2.23455 + 1.94071 i,\\
x_3&\approx -2.23455 - 1.94071 i,\\
x_4&\approx -0.3305.
\end{align}

\subsection{Coordinates on $\mc S(\P^1)\otimes \mc S(\P^2)$}
Consider the spaces $\mathcal S(\mathbb P^1),\mathcal S(\mathbb P^2)$ of solutions of the $qDE$'s of $\mathbb P^1$ and $\mathbb P^2$ specialized at the origins of $H^2(\P^1,\C)$ and $H^2(\P^2,\C)$, respectively: these equations are
\begin{align}
\label{qdep1}
\vartheta^2\Phi_1=&4z^2\Phi_1,\\
\label{qdep2}
\vartheta^3\Phi_2=&27z^3\Phi_2.
\end{align}
Solutions $\Phi_1(z)$ of equation \eqref{qdep1} have the following expansion at $z=0$:
\beq\label{phi1}
\Phi_1(z)=\sum_{m=0}^\infty(A_{m,1}+A_{m,0}\log z)\frac{z^{2m}}{(m!)^2},
\eneq
where $A_{0,0},A_{0,1}$ are arbitrary complex numbers, and the other coefficients are uniquely determined by the difference equations
\begin{align}
\label{02.12.18-1}A_{m-1,0}&=A_{m,0},\\
\label{02.12.18-2}A_{m-1,1}&=\frac{A_{m,0}}{m}+A_{m,1}.
\end{align}
in particular, notice that from the equation \eqref{02.12.18-2} we deduce that
\begin{equation}\label{02.12.18-6}A_{m,1}=A_{0,1}-A_{0,0}H_m,\quad m\geq 0,
\end{equation}where $H_m:=\sum_{i=1}^m\frac{1}{i}$ denotes the $m$-th harmonic number.

Analogously, solutions $\Phi_2(z)$ of equation \eqref{qdep2} have the following expansion at $z=0$:
\beq\label{phi2}
\Phi_2(z)=\sum_{n=0}^\infty(B_{n,2}+B_{n,1}\log z+B_{n,0}\log^2z)\frac{z^{3n}}{(n!)^3},
\eneq
where $B_{0,0},B_{0,1},B_{0,2}$ are arbitrary complex numbers, and the other coefficients are uniquely determined by the difference equations
\begin{align}
\label{02.12.18-3} B_{n-1,0}&=B_{n,0},\\
\label{02.12.18-4} B_{n-1,1}&=\frac{2}{n}B_{n,0}+B_{n,1},\\
\label{02.12.18-5} B_{n-1,2}&=\frac{2}{3n^2}B_{n,0}+\frac{1}{n}B_{n,1}+B_{n,2}.
\end{align}
From the difference equation \eqref{02.12.18-4} we deduce that
\begin{align}
\label{02.12.18-7} B_{n,1}&=B_{0,1}-2B_{0,0}H_n.
\end{align}
The products $A_{0,i}B_{0,j}$, with $i=0,1$ and $j=0,1,2$, define a natural system of coordinates on the tensor product $\mathcal S(\mathbb P^1)\otimes_{\mathbb C}\mathcal S(\mathbb P^2)$.

\subsection{Solutions of the $qDE$ of $\Fb_1$ as Laplace $(1,2;\frac{1}{2},\frac{1}{3})$-multitransforms}\label{LapF1} According to Theorem \ref{TH2}, the space of solutions of the quantum differential equation \eqref{qdiffdisp2bis} can be reconstructed from the spaces of solutions of the $qDE$'s \eqref{qdep1} and \eqref{qdep2}. From the polynomial equation \eqref{hirze1}, indeed, it follows that Theorem \ref{TH2} applies with the specialization of the parameters $h=2$, $\bm \ell=(2,3)$, $\bm d=(1,1)$.

Hence, we expect to reconstruct the solutions of the differential equation \eqref{qdiffdisp2bis} via a $\C$-bilinear operator $\mathscr P\colon \mc S(\P^1)\otimes \mc S(\P^2)\to \mc O(\widetilde{\C^*})$ involving the Laplace $(1,2;\frac{1}{2},\frac{1}{3})$-multitransform:
\[
\mathscr P[\Phi_1,\Phi_2](z):=e^{-cz}\mathscr L_{(1,2;\frac{1}{2},\frac{1}{3})}[\Phi_1,\Phi_2]=e^{-cz}\int_0^\infty\Phi_1\left(z^\frac{1}{2}\la^\frac{1}{2}\right)\Phi_2\left(z^\frac{2}{3}\la^\frac{1}{3}\right)e^{-\la}d\la,
\]
for a suitable number $c\in\Q$ to be determined. 

\begin{lem}
We have $c=1$.
\end{lem}
\proof
Along the locus of small quantum cohomology, the $J$-function of $\P^{n-1}$ is $$J_{\P^{n-1}}(\delta)=e^{\frac{\delta}{\hbar}}\sum_{d=0}^\infty{\bf Q}^de^{dt}\frac{1}{\left(\prod_{k=1}^d(H+k\hbar)\right)^{n}},\quad \delta=tH,$$ where $H\in H^2(\P^{n-1},\C)$ denotes the hyperplane class. Hence, the $I$-function $I_{\P^1\times\P^2,\Fb_1}$ equals
\begin{align*}
&I_{\P^1\times\P^2,\Fb_1}(\delta_1\otimes 1+1\otimes\delta_2)=e^\frac{\delta_1}{\hbar}\otimes e^{\frac{\delta_2}{\hbar}}\cdot\\
&\cdot\sum_{d_1,d_2\geq 0}{\bf Q}_1^{d_1}{\bf Q}_2^{d_2}\frac{e^{t^1d_1}}{\left(\prod_{k=1}^{d_1}(H_1+k\hbar)\right)^{2}}\otimes \frac{e^{t^2d_2}}{\left(\prod_{k=1}^{d_2}(H_2+k\hbar)\right)^{3}}\prod_{j=1}^{d_1+d_2}(H_1\otimes 1+1\otimes H_2+j\hbar)\\
&=1+\frac{1}{\hbar}\left({\bf Q}_1^{d_1}e^{t^1}+\delta_1\otimes 1+1\otimes \delta_2\right)+O\left(\frac{1}{\hbar^2}\right),
\end{align*}
where we set
\begin{itemize}
\item $H_1\in H^2(\P^1,\C)$ and $H_2\in H^2(\P^2,\C)$ are the hyperplane classes, 
\item  $\delta_1=t^1H_1$ and $\delta_2=t^2 H_2$ with $t^1,t^2\in\C$,
\item  ${\bf Q}_i={\bf Q}^{\bt _i}$, $\bt_i$ being the dual homology class of $H_i$, for $i=1,2$.
\end{itemize}
In the notations of Proposition \ref{ql3}, we have $H(\delta_1\otimes 1+1\otimes \delta_2)={\bf Q}_1^{d_1}e^{t^1}$. The number $c$ equals
\[\puqed
c=H(0)|_{\bf Q=1}=1.\qedhere
\poqed
\]

For brevity, in all the remaining part of this section, we will simply write $\mathscr L$ to denote the Laplace $(1,2;\frac{1}{2},\frac{1}{3})$-multitransform.

\subsubsection{The subspace $\mc H$} The space $\mc S(\P^1)\otimes \mc S(\P^2)$ has dimension 6. We are going to identify a subspace $\mc H$ of dimension 4 which is isomorphically mapped to the space $\mc S(\Fb_1)$ via the operator $\mathscr P$.

\begin{thm}
Let $\Phi_1(z)$ and $\Phi_2(z)$ be two solutions of the quantum differential equations of $\mathbb P^1$ and $\mathbb P^2$ respectively, namely
\[\vartheta^2\Phi_1(z)=4z^2\Phi_1(z),\quad \vartheta^3\Phi_2(z)=27z^3\Phi_2(z). 
\]The function
\[\Phi(z):=e^{-z}\mathscr L[\Phi_1,\Phi_2;z]
\]is a solution of the quantum differential equation of $\mathbb F_1$ if the following vanishing conditions are satisfied:
\[\mathscr D_1[\Phi_1,\Phi_2;z]=0,\quad \mathscr D_2[\Phi_1,\Phi_2;z]=0,
\]where
\begin{align*}\mathscr D_1[\Phi_1,\Phi_2;z]:=&2z^2\mathscr L[\vartheta\Phi_1,\Phi_2;z]-\frac{2}{9}\mathscr L[\vartheta\Phi_1,\vartheta^2\Phi_2;z]+\frac{4}{9}z\mathscr L[\Phi_1,\vartheta^2\Phi_2;z],\\
\mathscr D_2[\Phi_1,\Phi_2;z]:=&z^3\mathscr L[\Phi_1,\Phi_2;z]-\frac{z^2}{3}\mathscr L[\Phi_1,\vartheta\Phi_2;z]\\
&-\frac{z}{9}\mathscr L[\Phi_1,\vartheta^2\Phi_2;z]+\frac{z}{6}\mathscr L[\vartheta\Phi_1,\vartheta\Phi_2;z].
\end{align*}
\end{thm}

\proof Let us look for solutions of the equation \eqref{qdiffdisp2bis} in the form
\[\Phi(z)=e^{-z}\mathscr L_{(1,2;\frac{1}{2},\frac{1}{3})}[\Phi_1,\Phi_2;z],
\]where $\Phi_1$ and $\Phi_2$ are solutions of the quantum differential equation for $\mathbb P^1$ and $\mathbb P^2$ respectively, that is
\begin{align}\label{26.11.18-1}\vartheta^2\Phi_1&=4z^2\Phi_1,\\
\label{26.11.18-2} \vartheta^3\Phi_2&=27z^3\Phi_2.
\end{align} Given arbitrary functions  $f$ and $g$, we have
\begin{align*}\mathscr L[s^2f(s),g(s);z]=&z\left\{\mathscr L[f(s),g(s);z]+\frac{1}{2}\mathscr L[\vartheta_sf(s),g(s);z]\right.\\
&\left.+\frac{1}{3}\mathscr L[f(s),\vartheta_sg(s);z]-\mathcal I(f,g)\right\},
\end{align*}with
\begin{equation}\label{28.11.18-1}\mathcal I(f,g):=\left.\lambda\cdot f(z^{\frac{1}{2}}\lambda^{\frac{1}{2}})g(z^{\frac{2}{3}}\lambda^{\frac{1}{3}})e^{-\lambda}\right|_{\lambda=0}^{\lambda=\infty}.
\end{equation}
Applying the previous identity to $\Phi_1$ and $\Phi_2$, and using equations \eqref{26.11.18-1},\eqref{26.11.18-2}, we deduce the following identities:
\begin{align*}
\mathscr L[\vartheta^2\Phi_1,\Phi_2; z]=&\ 4z\left\{\mathscr L[\Phi_1,\Phi_2;z]+\frac{1}{2}\mathscr L[\vartheta\Phi_1,\Phi_2;z]+\frac{1}{3}\mathscr L[\Phi_1,\vartheta \Phi_2;z]\right\}+\mathcal R_1,\\
\mathscr L[\vartheta^3\Phi_1,\Phi_2;z]=&\ 8(z+z^2)\mathscr L[\Phi_1,\Phi_2;z]+(8z+4z^2)\mathscr L[\vartheta\Phi_1,\Phi_2;z]\\&+\frac{8}{3}\left(z+z^2\right)\mathscr L[\Phi_1,\vartheta\Phi_2;z]+\frac{4}{3}z\mathscr L[\vartheta\Phi_1,\vartheta\Phi_2;z]+\mathcal R_2,\\
\mathscr L[\vartheta^4\Phi_1,\Phi_2;z]= &\ 16(z+4z^2+z^3)\mathscr L[\Phi_1,\Phi_2;z]+8(3z+5z^2+z^3)\mathscr L[\vartheta\Phi_1,\Phi_2;z]\\
&+\frac{16}{3}(z+5z^2+z^3)\mathscr L[\Phi_1,\vartheta\Phi_2;z]+\frac{16}{3}(z+z^2)\mathscr L[\vartheta\Phi_1,\vartheta\Phi_2;z]\\
&+\frac{16}{9}z^2\mathscr L[\Phi_1,\vartheta^2\Phi_2;z]+\mathcal R_3,\\
\mathscr L[\Phi_1,\vartheta^3\Phi_2;z]=&\ 27z^2\left\{\mathscr L[\Phi_1,\Phi_2;z]+\frac{1}{2}\mathscr L[\vartheta\Phi_1,\Phi_2;z]+\frac{1}{3}\mathscr L[\Phi_1,\vartheta \Phi_2;z]\right\}+\mathcal R_4,\\
\mathscr L[\Phi_1,\vartheta^4\Phi_2;z]=&\ \frac{9}{2}z^2\Bigl\{ 18\mathscr L[\Phi_1,\Phi_2;z]+12\mathscr L[\Phi_1,\vartheta \Phi_2;z]+2\mathscr L[\Phi_1,\vartheta^2\Phi_2;z]\\
&+9\mathscr L[\vartheta\Phi_1,\Phi_2;z]+3\mathscr L[\vartheta\Phi_1,\vartheta \Phi_2;z]\Bigr\}+\mathcal R_5,\\
\mathscr L[\vartheta\Phi_1,\vartheta^3\Phi_2;z]=&\ 54z^3\mathscr L[\Phi_1,\Phi_2;z]+27(z^2+z^3)\mathscr L[\vartheta\Phi_1,\Phi_2;z]+18z^3\mathscr L[\Phi_1,\vartheta\Phi_2;z]\\ &+9z^2\mathscr L[\vartheta\Phi_1,\vartheta\Phi_2;z]+\mathcal R_6,\\
\mathscr L[\vartheta^2\Phi_1,\vartheta^2\Phi_2;z]=&\ 36z^3\mathscr L[\Phi_1,\Phi_2;z]+18z^3\mathscr L[\vartheta\Phi_1,\Phi_2;z]+12 z^3 \mathscr L[\Phi_1,\vartheta\Phi_2;z]\\
&+4 z \mathscr L[\Phi_1,\vartheta^2\Phi_2;z]+2 z \mathscr L[\vartheta\Phi_1,\vartheta^2\Phi_2;z]+\mathcal R_7,\\
\mathscr L[\vartheta^3\Phi_1,\vartheta\Phi_2;z]=&\ \ 8(z+z^2)\mathscr L[\Phi_1,\vartheta\Phi_2;z]+(8z+4z^2)\mathscr L[\vartheta\Phi_1,\vartheta\Phi_2;z]\\&+\frac{8}{3}\left(z+z^2\right)\mathscr L[\Phi_1,\vartheta^2\Phi_2;z]+\frac{4}{3}z\mathscr L[\vartheta\Phi_1,\vartheta^2\Phi_2;z]+\mathcal R_8,\\
\mathscr L[\vartheta^2\Phi_1,\vartheta\Phi_2; z]=&\ 4z\left\{\mathscr L[\Phi_1,\vartheta\Phi_2;z]+\frac{1}{2}\mathscr L[\vartheta\Phi_1,\vartheta\Phi_2;z]+\frac{1}{3}\mathscr L[\Phi_1,\vartheta^2 \Phi_2;z]\right\}+\mathcal R_9,
\end{align*}
where $\mathcal R_j$ with $j=1,\dots, 9$ denote some negligible boundary terms due to the cumulations of terms like \eqref{28.11.18-1}. Using these identities, after some computations, we can rewrite the quantum differential equation \eqref{qdiffdisp2bis} as follows
\[
\puqed(-72+1674z+283z^2)\mathscr D_1[\Phi_1,\Phi_2]+(36+724z+4811 z^2)\mathscr D_2[\Phi_1,\Phi_2]=0.\qedhere
\poqed
\]

An explicit computation shows that $\mathscr D_1[\Phi_1,\Phi_2;z]$ and $\mathscr D_2[\Phi_1,\Phi_2;z]$ have the following expansions
\begin{align}
\mathscr D_1[\Phi_1,\Phi_2;z]&=\Theta_1(z)\log^3z+\Theta_2(z)\log^2z+\Theta_3(z)\log z+\Theta_4(z),\\
\mathscr D_2[\Phi_1,\Phi_2;z]&=\Lambda_1(z)\log^3z+\Lambda_2(z)\log^2z+\Lambda_3(z)\log z+\Lambda_4(z),
\end{align}
where the functions $\Theta_i(z)$ and $\Lambda_i(z)$ are of the form
\begin{align}
\label{theta}
\Theta_i(z)=&\sum_{m=0}^\infty\sum_{n=0}^\infty\frac{(m+n)!}{(m!)^2(n!)^3}\left(\mathcal A_1^{(i)}(m,n)+\mathcal A_2^{(i)}(m,n)z+\mathcal A_3^{(i)}(m,n)z^2\right)z^{m+2n},\\
\label{lambda}
\Lambda_i(z)=&\sum_{m=0}^\infty\sum_{n=0}^\infty\frac{(m+n)!}{(m!)^2(n!)^3}\left(\mathcal B_1^{(i)}(m,n)+\mathcal B_2^{(i)}(m,n)z+\mathcal B_3^{(i)}(m,n)z^2\right)z^{m+2n+1},
\end{align}
for $i=1,2,3,4$. See Appendix \ref{appa} for the explicit expressions of the coefficients $\mathcal A_j^{(i)},\mathcal B_j^{(i)}$. 

\begin{lem}
For all $m,n\geq 1$ and $i=1,2,3,4$, the following identities hold true
\begin{align}
\label{idubro1}
(m+n)\mathcal A_1^{(i)}(m,n)+m^2\mathcal A_1^{(i)}(m-1,n)+n^3\mathcal A_1^{(i)}(m,n-1)&=0,\\
\label{idubro2}
(m+n)\mathcal B_1^{(i)}(m,n)+m^2\mathcal B_1^{(i)}(m-1,n)+n^3\mathcal B_1^{(i)}(m,n-1)&=0,\\
\label{id5}
\mathcal A_1^{(i)}(m,0)+m\mathcal A_2^{(i)}(m-1,0)&=0,\\
\label{id6}
\mathcal B_1^{(i)}(m,0)+m\mathcal B_2^{(i)}(m-1,0)&=0,\\
\label{idubro4}
\mathcal A_1^{(i)}(0,n)+n^2\mathcal A_3^{(i)}(0,n-1)&=0,\\
\label{id4}
\mathcal B_1^{(i)}(0,n)+n^2\mathcal B_3^{(i)}(0,n-1)&=0.
\end{align}
\end{lem}
\proof
The reader can check the validity of these identities using the explicit expressions in Appendix \ref{appa}, equations \eqref{02.12.18-1}, \eqref{02.12.18-2}, \eqref{02.12.18-3}, \eqref{02.12.18-4}, \eqref{02.12.18-5}, and the following identities (see e.g. \cite{NIST}):
\begin{align*}
\psi^{(k)}(z+1)&=\psi^{(k)}(z)+\frac{(-1)^kk!}{z^{k+1}},\quad k\geq 0,\\
\psi^{(0)}(n)&=H_{n-1}-\gamma,\quad n\geq 1,\quad \psi(z):=\frac{\Gamma'(z)}{\Gamma(z)}.\qedhere
\end{align*}

\begin{thm}\label{teosolf1}
Let $\Phi_1(z)\in\mathcal S(\mathbb P^1), \Phi_2(z)\in\mathcal S(\mathbb P^2)$ be as in equations \eqref{phi1} and \eqref{phi2}, respectively. Then the function $\Phi(z):=e^{-z}\mathscr L[\Phi_1,\Phi_2;z]$ is a solution of the $qDE$ of $\mathbb F_1$ if
\beq
\label{constr}
A_{0,0}B_{0,0}=0,\quad 4A_{0,1}B_{0,0}=3A_{0,0}B_{0,1}.
\eneq
\end{thm}
\proof
Let us rearrange the double series \eqref{theta} as follows:
\begin{align*}
\Theta_i(z)=&\Biggl\{\mathcal A_1^{(i)}(0,0)+\cancel{\sum_{m=1}^\infty\sum_{n=1}^\infty\frac{(m+n)!}{(m!)^2(n!)^3}\mathcal A_1^{(i)}(m,n)z^{m+2n}}\\
&+\cancel{{\color{red}\sum_{m=1}^\infty\frac{1}{m!}\mathcal A_1^{(i)}(m,0)z^{m}}}+\cancel{{\color{blue}\sum_{n=1}^\infty\frac{1}{(n!)^2}\mathcal A_1^{(i)}(0,n)z^{2n}}}\\
&+\cancel{\sum_{m=0}^\infty\sum_{n=1}^\infty\frac{(m+n)!}{(m!)^2(n!)^3}\mathcal A_2^{(i)}(m,n)z^{1+m+2n}}+\cancel{{\color{red}\sum_{m=0}^\infty\frac{1}{m!}\mathcal A_2^{(i)}(m,0)z^{1+m}}}\\
&+\cancel{\sum_{m=1}^\infty\sum_{n=0}^\infty\frac{(m+n)!}{(m!)^2(n!)^3}\mathcal A_3^{(i)}(m,n)z^{2+m+2n}}+\cancel{{\color{blue}\sum_{n=0}^\infty\frac{1}{(n!)^2}\mathcal A_3^{(i)}(0,n)z^{2+2n}}}\Biggr\},
\end{align*}
where 
\begin{enumerate}
\item the black summands cancel by equation \eqref{idubro1},
\item the red summands cancel by equation \eqref{id5},
\item the blue summands cancels by equation \eqref{idubro4}.
\end{enumerate}
The proof for $\Lambda_i(z)$ is identical.
\endproof

\begin{defn}
Let $\mc H$ be the 4-dimensional subspace of $\subseteq \mc S(\P^1)\otimes \mc S(\P^2)$ defined by the linear equations \eqref{constr}.
\end{defn}
\begin{cor}
The space $\mc H$ is isomorphic to the space of solutions $\mc S(\Fb_1)$ via the operator $\mathscr P$.\qed
\end{cor}
\subsubsection{Bases of $\mc S(\P^1)$} Define
\beq\label{g}
g(z):=\frac{1}{2\pi i}\int_{\mc L_1}\Gm\left(\frac{s}{2}\right)^2z^{-s}ds,
\eneq
where $\mc L_1$ is a (positively oriented) parabola ${\rm Re\ }s=-c\cdot ({\rm Im\ }s)^2+c'$, for suitable $c,c'\in\R_+$ so that it encircles all the poles of the integrand at $s\in2\Z_{\leq 0}$. It is easy to see that the integral \eqref{g} converges for all $z\in\widetilde{\C^*}$ and that its value does not depend on the particular choice of $c,c'$.

\begin{prop}
The functions $g(e^{-i\pi}z),g(z)$
define a basis of solutions of the $qDE$ of $\P^1$.
\qed
\end{prop}

Define the bases $(g_1(z),g_2(z))$ and $(s_1(z), s_2(z))$ of $\mc S(\P^1)$ by 
\beq
\left(
\begin{array}{c}
g_1(z)\\
g_2(z)
\end{array}\right)=M_1\left(\begin{array}{c}
g(e^{-\pi i}z)\\
g(z)
\end{array}\right),\quad
\left(
\begin{array}{c}
s_1(z)\\
s_2(z)
\end{array}\right)=M_2\left(\begin{array}{c}
g(e^{-\pi i}z)\\
g(z)
\end{array}\right),
\eneq
where 
\beq
M_1:=\left(
\begin{array}{cc}
 -\frac{i \gamma }{4 \pi } & \frac{i (\gamma +i \pi )}{4 \pi } \\
 \frac{i}{4 \pi } & -\frac{i}{4 \pi } \\
\end{array}
\right),\quad 
M_2:=\begin{pmatrix}
-1&2\\
0&1
\end{pmatrix}.
\eneq
\begin{lem}\label{gi}
For $z\to 0$, the following asymptotic expansions hold true:
\begin{align}
g_1(z)&=\log z +O(z^2\log z),\\
g_2(z)&=1+O(z^2\log z).
\end{align}
\end{lem}
\proof
The proof is a simple computation of residues: by modifying the paths of integration $\mc L_1$, one obtains the asymptotic expansions of $g$ as a sum of residues of the integrand.
\endproof

\begin{lem}\label{ai1}
We have
\[g(z)\sim\frac{2\pi^\frac{1}{2}}{z^\frac{1}{2}}e^{-2z},\quad z\to\infty
\]in the sector $|\arg z|<\frac{3}{2}\pi$.
\end{lem}
\proof
The estimate follows from application of steepest descent method.
\endproof

\subsubsection{Bases of $\mc S(\P^2)$} Define
\beq
\label{h}
h(z):=\frac{1}{2\pi i}\int_{\mc L_2}\Gm\left(\frac{s}{3}\right)^3e^{\frac{\pi i s}{3}}z^{-s}ds,
\eneq
where $\mc L_2$ is a (positively oriented) parabola ${\rm Re\ }s=-c\cdot ({\rm Im\ }s)^2+c'$, for suitable $c,c'\in\R_+$ so that it encircles all the poles of the integrand at $s\in3\Z_{\leq 0}$. It is easy to see that the integral \eqref{h} converges for all $z\in\widetilde{\C^*}$ and that its value does not depend on the particular choice of $c,c'$.

\begin{prop}
The functions $h(e^{-\frac{2i\pi}{3}}z),h(z),h(e^\frac{2i\pi}{3}z)$
define a basis of solutions of the $qDE$ of $\P^2$.
\qed
\end{prop}

Define the bases $(h_1(z), h_2(z), h_3(z))$ and $(p_1(z),p_2(z),p_3(z))$ of $\mc S(\P^2)$ by
\beq\label{baseshp}
\left(
\begin{array}{c}
h_1(z)\\
h_2(z)\\
h_3(z)
\end{array}\right)=N_1\left(\begin{array}{c}
h(e^{-\frac{2i\pi}{3}}z)\\
h(z)\\
h(e^\frac{2i\pi}{3}z)
\end{array}\right),
\quad
\left(
\begin{array}{c}
p_1(z)\\
p_2(z)\\
p_3(z)
\end{array}\right)=N_2\left(\begin{array}{c}
h(e^{-\frac{2i\pi}{3}}z)\\
h(z)\\
h(e^\frac{2i\pi}{3}z)
\end{array}\right),
\eneq
where
\beq
N_1:=
\left(
\begin{array}{ccc}
 \frac{-18 \gamma ^2-\pi ^2}{216 \pi ^2} & \frac{-18 \gamma ^2-24 i \gamma  \pi +7 \pi ^2}{216 \pi ^2} & \frac{18 \gamma ^2+12 i \gamma  \pi +5 \pi ^2}{108 \pi ^2} \\
 \\
 \frac{\gamma }{12 \pi ^2} & \frac{3 \gamma +2 i \pi }{36 \pi ^2} & \frac{-3 \gamma -i \pi }{18 \pi ^2} \\
 \\
 -\frac{1}{12 \pi ^2} & -\frac{1}{12 \pi ^2} & \frac{1}{6 \pi ^2} \\
\end{array}
\right),\quad N_2:=\begin{pmatrix}
-1&3&-3\\
0&1&0\\
0&0&-1
\end{pmatrix}.
\eneq
The basis $(p_1,p_2,p_3)$ will be studied later, in Section \ref{stbqdef1}, where it will be used to construct Stokes bases of solutions. We now focus on the properties of the basis $(h_1,h_2,h_3)$.
\begin{lem}\label{hi}
For $z\to 0$, the following asymptotic expansions hold true:
\begin{align}
h_1(z)&=\log^2 z+O(z^3\log^2 z),\\
h_2(z)&=\log z+O(z^3\log^2 z),\\
h_3(z)&=1+O(z^3\log^2 z).
\end{align}
\end{lem}
\proof
The proof is a simple computations of residues: by modifying the paths of integration $\mc L_2$, one obtains the asymptotic expansions of $h$ as a sum of residues of the integrand.
\endproof

\begin{lem}\label{ai2}
We have
\[h(z)\sim e^{-\frac{5}{3}\pi i}\frac{\sqrt{3}}{z}\exp\left(3e^{\frac{2\pi i}{3}}z\right),\quad z\to\infty,
\]in the sector $-\pi<\arg z<\frac{5}{3}\pi$.
\end{lem}
\proof
The estimate follows from the steepest descent method.
\endproof

\subsection{Basis of solutions $\Ups$ of $\mc S(\Fb_1)$}

\begin{thm}
The tensors 
\beq
\label{basis}
\frac{1}{3}g_1\otimes h_2+\frac{1}{4}g_2\otimes h_1,\quad g_1\otimes h_3,\quad g_2\otimes h_2,\quad g_2\otimes h_3.
\eneq
define a basis of the subspace $\mc H$.
\end{thm}

\proof
Each of the vectors \eqref{basis} satisfy the constraints \eqref{constr}, by Lemmata \ref{gi} and \ref{hi}. 
\endproof

\begin{cor}\label{corUps}
The functions
\begin{align}
\Ups_1&:=\mathscr P\left(\frac{1}{3}g_1\otimes h_2+\frac{1}{4}g_2\otimes h_1\right),\\
\Ups_2&:=\mathscr P\left( g_1\otimes h_3\right),\\
\Ups_3&:=\mathscr P\left( g_2\otimes h_2\right),\\
\Ups_4&:=\mathscr P\left( g_2\otimes h_3\right)
\end{align}
define a basis of solutions of the $qDE$ of $\Fb_1$.\qed
\end{cor}

\begin{rem}
Explicit double Mellin-Barnes integral representations of solutions $\Ups_1,\dots,\Ups_4$ can be obtained:  for any $j,k$ we have
\begin{align*}&\mathscr P\left(g(e^{\pi k i}z)\otimes h(e^{\frac{2\pi j i}{3}}z)\right)\\
&=\frac{e^{-z}}{(2\pi i)^2}\int_{\mc L_1\times \mc L_2}\Gamma\left(\frac{s}{2}\right)^2\Gamma\left(\frac{t}{3}\right)^3\Gamma\left(1-\frac{s}{2}-\frac{t}{3}\right)e^{-\pi i k s+\frac{\pi i }{3}t(1-2j)}z^{-\frac{s}{2}-\frac{2t}{3}}dtds.
\end{align*}
The functions $\Ups_i$'s are linear combinations of the integrals above, in accordance with Theorem \ref{TH2b}.
\end{rem}

\subsection{Asymptotics of Laplace $(1,2;\frac{1}{2},\frac{1}{3})$-multitransforms}\label{aLap}Consider the integral
\beq
\mc I(z):=\int_0^\infty\Phi_1(z^\frac{1}{2}\la^\frac{1}{2})\Phi_2(z^\frac{2}{3}\la^\frac{1}{3})e^{-\la}d\la,
\eneq
where
\beq
\Phi_1(z)=z^{D_1}\exp(zu_1),\quad \Phi_2(z)=z^{D_2}\exp(zu_2),
\eneq
with $D_1,D_2,u_1,u_2\in\C$. The integral $\mc I(z)$ is convergent for all $z\in\widetilde{\C^*}$.

Set $z=re^{i\ka}$ with $r>0$, and change variable of integration $\la=\al z$:
\beq
\mc I(z)=z^{1+D_1+D_2}\int_0^{e^{-i\ka}\infty}\al^{\frac{D_1}{2}+\frac{D_2}{3}}\exp\left\{z(-\al+u_1\al^\frac{1}{2}+u_2\al^\frac{1}{3})\right\}d\al.
\eneq
Change variable $\al=\bt^6$, by taking the principal determination of the sixth root:
\beq\label{intI}
\mc I(z)=6z^{1+D_1+D_2}\int_0^{e^{-\frac{i\ka}{6}}\infty}\bt^{5+3D_1+2D_2}\exp\left\{z(-\bt^6+u_1\bt^3+u_2\bt^2)\right\}d\bt.
\eneq
Define 
\beq
f(\bt;u_1,u_2):=-\bt^6+u_1\bt^3+u_2\bt^2,\quad\text{ for }\bt\in\C,
\eneq and consider the $z$-dependent downward flow in the $\bt$-plane defined by
\beq
\label{flow}
\frac{d\bt}{dt}=-\overline{z}\frac{\partial \overline f}{\partial \overline{\bt}},\quad \frac{d\overline\bt}{dt}=-{z}\frac{\partial f}{\partial {\bt}}.
\eneq
The equilibria points $\bt_c$ are the critical points of $f$, that is$$\left.\frac{\partial f}{\partial \bt}\right|_{\bt=\bt_c}=0.$$
For a fixed $z$, we associate to each critical point $\bt_c$ a curve $\mc L_c$, a \emph{Lefschetz thimble}, defined as the set-theoretic union of the trajectories of the flow \eqref{flow} starting at $\bt_c$ for $t\to-\infty$. Morse and Picard-Lefschetz Theory guarantees that the cycles $\mc L_c$ are smooth one-dimensional submanifolds of $\C$, piecewise smoothly dependent on the parameter $z$, and they represent a basis for the inverse limit of relative homology groups 
\[\varprojlim_T H_1(\C,\C_{T,z}),\quad \C_{T,z}:=\left\{\bt\in\C\colon {\rm Re}(zf(\bt;u_1,u_2))<-T\right\},\quad T\in\R_+.
\]
\begin{lem}
The Lefschetz thimble $\mc L_c$ is the steepest descent path at $\bt_c$: the function $t\mapsto {\rm Im}(zf(\bt;u_1,u_2))$ is constant on $\mc L_c$, and the function $t\mapsto {\rm Re}(zf(\bt;u_1,u_2))$ is strictly decreasing along the flow.
\end{lem}
\proof We have
\[
\frac{d}{dt}[{\rm Im}(zf)]=\left(\frac{d\bt}{d t}\frac{\der }{\der \bt}+\frac{d\overline\bt}{d t}\frac{\der }{\der \overline\bt}\right)\left[\frac{zf-\overline{zf}}{2i}\right]=0,\]
\[
\puqed
\frac{d}{dt}[{\rm Re}(zf)]=\left(\frac{d\bt}{d t}\frac{\der }{\der \bt}+\frac{d\overline\bt}{d t}\frac{\der }{\der \overline\bt}\right)\left[\frac{zf+\overline{zf}}{2}\right]=-\left|{z}\frac{\partial  f}{\partial {\bt}}\right|^2.\qedhere
\poqed
\]
We are interested in the following cases, by Lemmata \ref{ai1} and \ref{ai2}:
\beq
\label{u1u2}
u_1=\pm 2,\quad u_2=3\zeta_3^k,\quad \zeta_3:=\exp\frac{2\pi i}{3},\quad k=0,1,2.
\eneq

\begin{table}
\[
\begin{array}{|c|c|c|c|c|}
\hline
u_1&u_2&\beta_c&f(\bt_c)&f(\bt_c)-1\\
\hline
\hline
 2 & 3 & -0.724492 & 0.6695 & -0.3305 \\
 2 & 3 & 0. & 0. & -1. \\
 2 & 3 & \clr{1.22074} & 4.7996 & 3.7996 \\
 2 & 3 & -0.248126-1.03398 i & -1.23455+1.94071 i & -2.23455+1.94071 i \\
 2 & 3 & -0.248126+1.03398 i & -1.23455-1.94071 i & -2.23455-1.94071 i \\
 \hline
 2 & 3 e^{\frac{2 i \pi }{3}} & 0. & 0. & -1. \\
 2 & 3 e^{\frac{2 i \pi }{3}} & -0.771392-0.731875 i & -1.23455-1.94071 i & -2.23455-1.94071 i \\
 2 & 3 e^{\frac{2 i \pi }{3}} & -0.610372+1.0572 i & 4.7996\,  & 3.7996\,  \\
 2 & 3 e^{\frac{2 i \pi }{3}} & 0.362246\, -0.627428 i & 0.6695\,  & -0.3305 \\
 2 & 3 e^{\frac{2 i \pi }{3}} & \clr{1.01952\, +0.302108 i }& -1.23455+1.94071 i & -2.23455+1.94071 i \\
 \hline
 2 & 3 e^{-\frac{1}{3} (2 i \pi )} & 0. & 0. & -1. \\
 2 & 3 e^{-\frac{1}{3} (2 i \pi )} & -0.771392+0.731875 i & -1.23455+1.94071 i & -2.23455+1.94071 i \\
 2 & 3 e^{-\frac{1}{3} (2 i \pi )} & -0.610372-1.0572 i & 4.7996\, & 3.7996\, \\
 2 & 3 e^{-\frac{1}{3} (2 i \pi )} & 0.362246\, +0.627428 i & 0.6695\, & -0.3305\\
 2 & 3 e^{-\frac{1}{3} (2 i \pi )} & \clr{1.01952\, -0.302108 i }& -1.23455-1.94071 i & -2.23455-1.94071 i \\
 \hline
 -2 & 3 & -1.22074 & 4.7996 & 3.7996 \\
 -2 & 3 & 0. & 0. & -1. \\
 -2 & 3 & \clr{0.724492} & 0.6695 & -0.3305 \\
 -2 & 3 & 0.248126\, -1.03398 i & -1.23455-1.94071 i & -2.23455-1.94071 i \\
 -2 & 3 & 0.248126\, +1.03398 i & -1.23455+1.94071 i & -2.23455+1.94071 i \\
 \hline
 -2 & 3 e^{\frac{2 i \pi }{3}} & 0. & 0. & -1. \\
 -2 & 3 e^{\frac{2 i \pi }{3}} & -1.01952-0.302108 i & -1.23455+1.94071 i & -2.23455+1.94071 i \\
 -2 & 3 e^{\frac{2 i \pi }{3}} & -0.362246+0.627428 i & 0.6695\,  & -0.3305 \\
 -2 & 3 e^{\frac{2 i \pi }{3}} & 0.610372\, -1.0572 i & 4.7996\,  & 3.7996\,  \\
 -2 & 3 e^{\frac{2 i \pi }{3}} & \clr{0.771392\, +0.731875 i} & -1.23455-1.94071 i & -2.23455-1.94071 i \\
 \hline
 -2 & 3 e^{-\frac{1}{3} (2 i \pi )} & 0. & 0. & -1. \\
 -2 & 3 e^{-\frac{1}{3} (2 i \pi )} & -1.01952+0.302108 i & -1.23455-1.94071 i & -2.23455-1.94071 i \\
 -2 & 3 e^{-\frac{1}{3} (2 i \pi )} & -0.362246-0.627428 i & 0.6695\, & -0.3305\\
 -2 & 3 e^{-\frac{1}{3} (2 i \pi )} & 0.610372\, +1.0572 i & 4.7996\, & 3.7996\, \\
 -2 & 3 e^{-\frac{1}{3} (2 i \pi )} & \clr{0.771392\, -0.731875 i} & -1.23455+1.94071 i & -2.23455+1.94071 i \\
 \hline
\end{array}
\]
\caption{For any possible value of the pair $(u_1,u_2)$, we list the corresponding critical points $\bt_c$ of the function $f(\bt;u_1,u_2)$, and the corresponding critical values $f(\bt_c)$. Notice that the numbers $f(\bt_c)-1$, with $\bt_c\neq 0$, equal all possible values of the canonical coordinates $x_1,x_2,x_3,x_4$ at the origin of $QH^\bullet(\Fb_1)$. In red, we represent the critical point $\bt_+$ with maximal real part.}\label{table1}
\end{table}

For any possible pair $(u_1,u_2)$, define $\bt_+$ as the critical point of $f(\bt;u_1,u_2)$ with maximal real part (the red one in Table \ref{table1}).

\begin{lem}\label{lemasI}
We have
\begin{align*}
\mc I(z)
&\sim 6z^{\frac{1}{2}+D_1+D_2}\bt_+^{5+2D_1+3D_2}\left(\frac{2\pi}{9u_1\bt_++8u_2}\right)^\frac{1}{2}\exp{z(-\bt_+^6+u_1\bt_+^3+u_2\bt_+^2)},
\end{align*}
for $|z|\to\infty$ in the sector $|\arg z-\arg \overline{f(\bt_+)}|<\pi$.
\end{lem}

\proof After choosing an orientation for each Lefschetz thimble, the path of integration $\gamma_z\equiv e^{-i\frac{\kappa}{6}}\cdot\R_+$, defining the function $\mc I$ in equation \eqref{intI}, can be expressed as integer combination, $\gamma_z=\sum_{j=1}^5 n_j(z)\mc L_j$ with $n_j\in\Z$, of the thimbles $\mc L_c$ for any value of $z$ not on a Stokes ray $\mc R_{ij}$, defined by
\[ \mc R_{ij}:=\left\{z\in\widetilde{\C^*}\colon z=r(\overline{f(\bt_{c,i})}-\overline{f(\bt_{c,j})}),\quad r\in\R_+\right\},\quad i,j=1,\dots,5,\]where $\bt_{c,i}$ are the critical points of \eqref{flow}. If we let $z$ vary, the Lefschetz thimbles change. When $z$ crosses a Stokes ray $\mc R_{ij}$, Lefschetz thimbles jump discontinuously: in particular, for $z$ on a Stokes ray there exists a flow line of \eqref{flow} connecting two critical points $\bt_c$'s. A detailed analysis of the phase portrait of the flow \eqref{flow}, for each pair $(u_1,u_2)$ as in \eqref{u1u2}, shows that in the sector $|\arg z-\arg \overline{f(\bt_+)}|<\pi$ we have $\mc I=\pm \mc L_{\bt_+}\pm \mc L_{0}^1\pm\mc L'$, where $\mc L_0^1$ is only one half of the Lefschetz thimble $\mc L_0$, and $\mc L'$ denote the sum of Lefschetz thimbles attached to other critical points $\bt_c$. Hence, we have three contributions in the asymptotics of $\mc I(z)$: one from the integration along $\mc L_{\bt_+}$, one from other critical points, the last one from the integration along $\mc L_0^1$. The last two contributions are easily seen to be negligible w.r.t. the first one. So, by the steepest descent method, we obtain the estimate
\[
\puqed
\mc I(z)\sim \pm 6z^{\frac{1}{2}+D_1+D_2}\bt_+^{5+2D_1+3D_2}\left(-\frac{2\pi}{f''(\bt_+)}\right)^\frac{1}{2}\exp{zf(\bt_+)}.\qedhere
\poqed
\]

\begin{rem}
Note that the arbitrariness of the orientations of the Lefschetz thimbles can be incorporated in the choice of the entries of the  $\Psi$-matrix. Consequently, it will affect the monodromy data by the action of the group $(\Z/2\Z)^4$.
\end{rem}

\begin{figure}
\includegraphics[width=0.45\textwidth, angle=0]{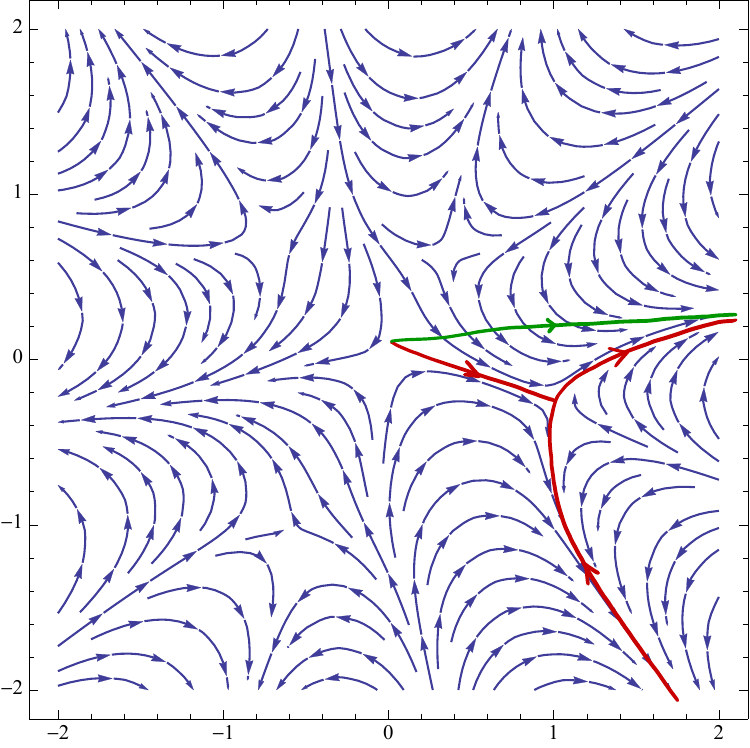}
\includegraphics[width=0.45\textwidth, angle=0]{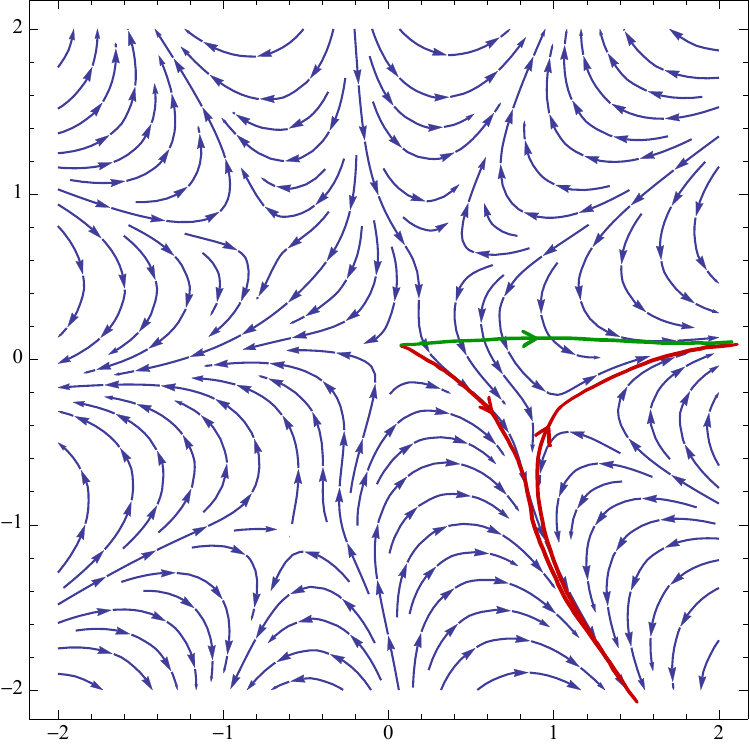}
\includegraphics[width=0.45\textwidth, angle=0]{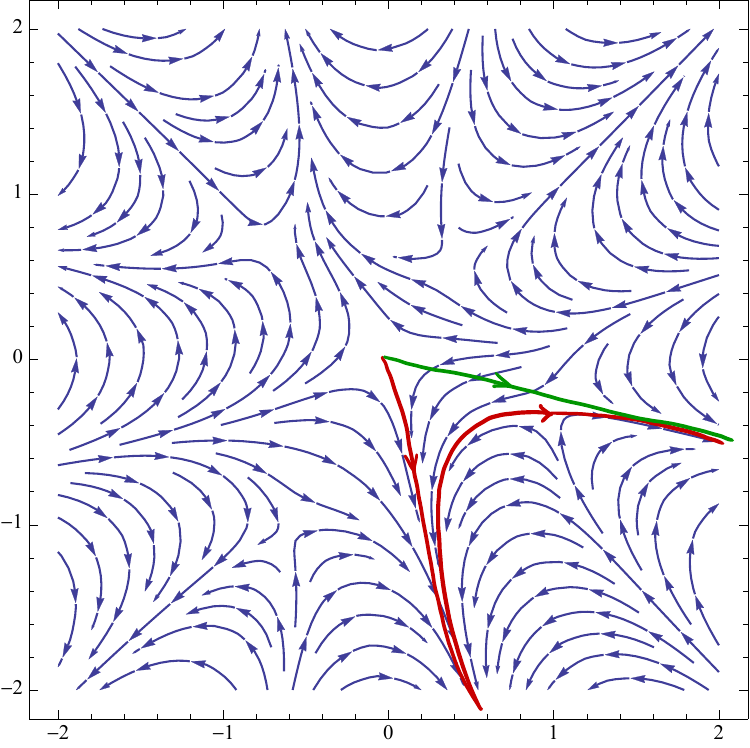}
\includegraphics[width=0.45\textwidth, angle=0]{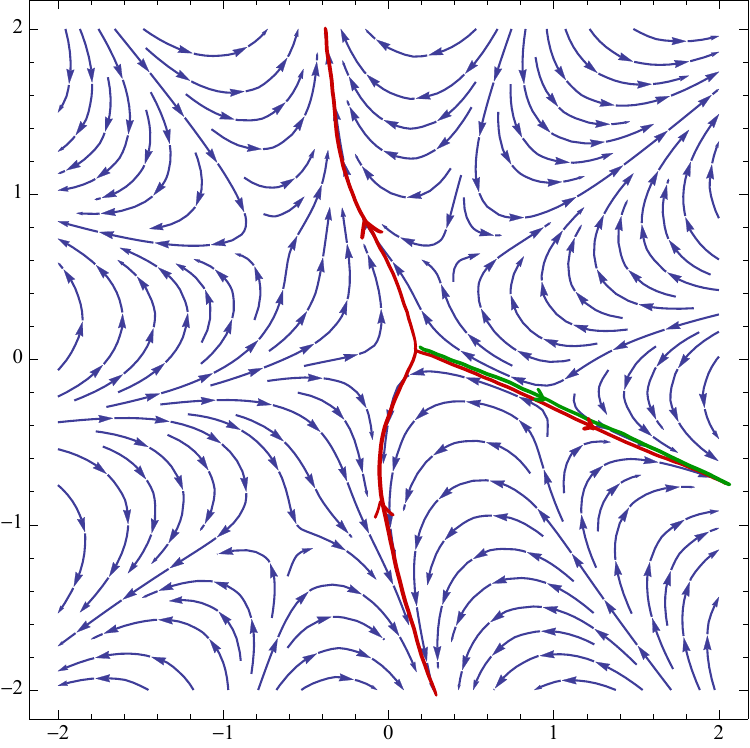}
\includegraphics[width=0.45\textwidth, angle=0]{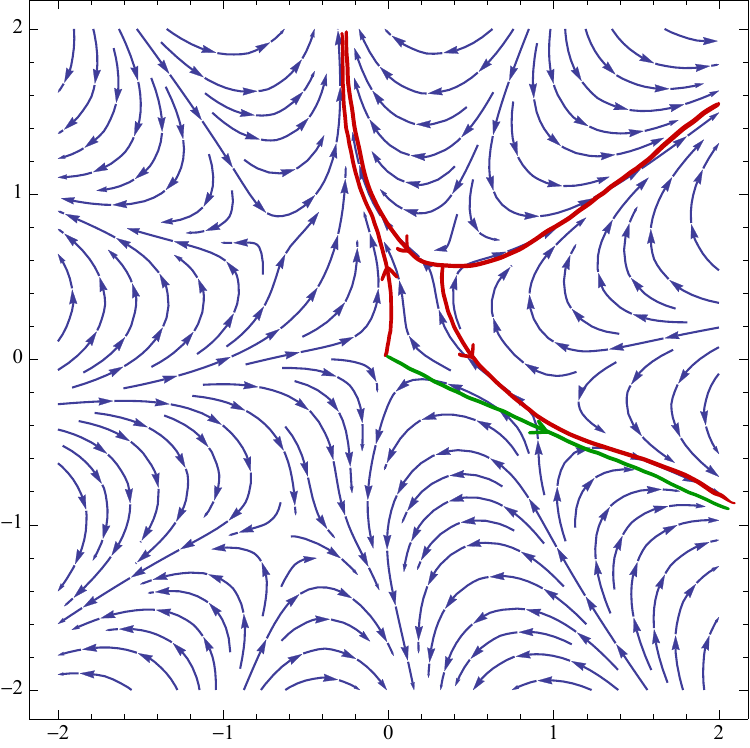}
\includegraphics[width=0.45\textwidth, angle=0]{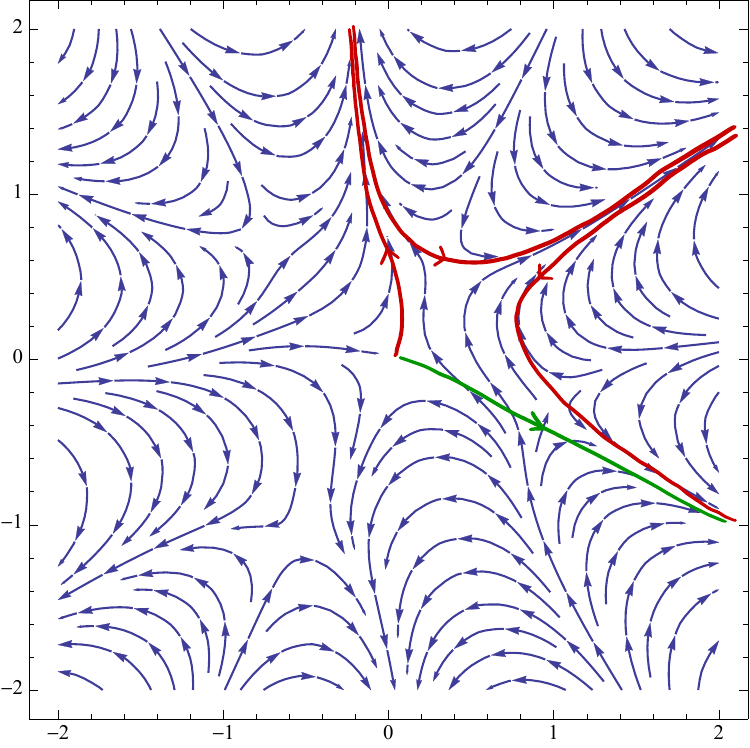}
\caption{In this figure we represent the downward flow \eqref{flow} and the mutations of Lefschetz thimbles for $|z|= 10^5$, and $|\arg z-\arg \overline{f(\bt_+)}|<\pi$ for the pair $(u_1,u_2)=(2,3e^{\frac{4\pi i}{3}})$. Lefschetz thimbles are in red. The path of integration in equation \eqref{intI} is drawn in green. {\it Continues in the next page.}}
\end{figure}

\begin{figure}
\includegraphics[width=0.45\textwidth, angle=0]{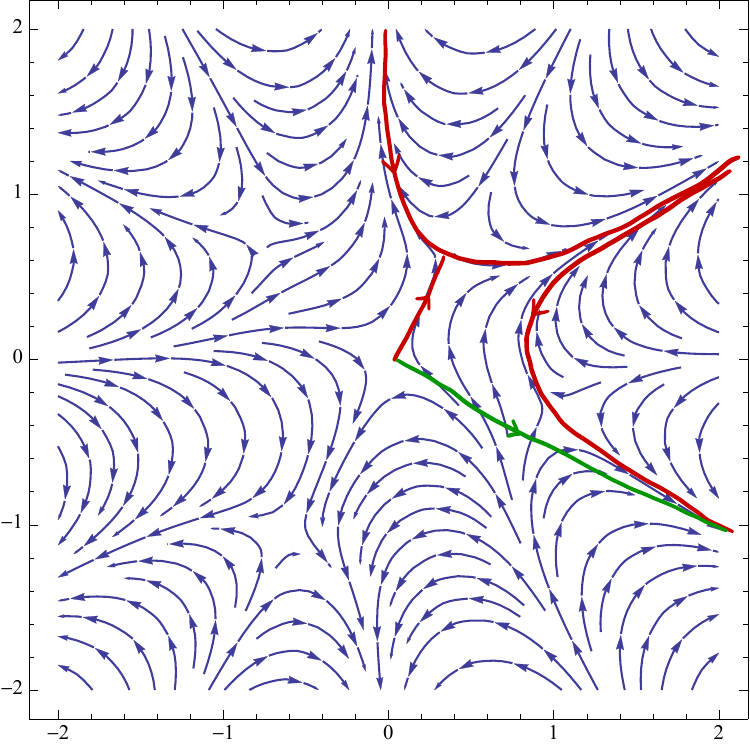}
\includegraphics[width=0.45\textwidth, angle=0]{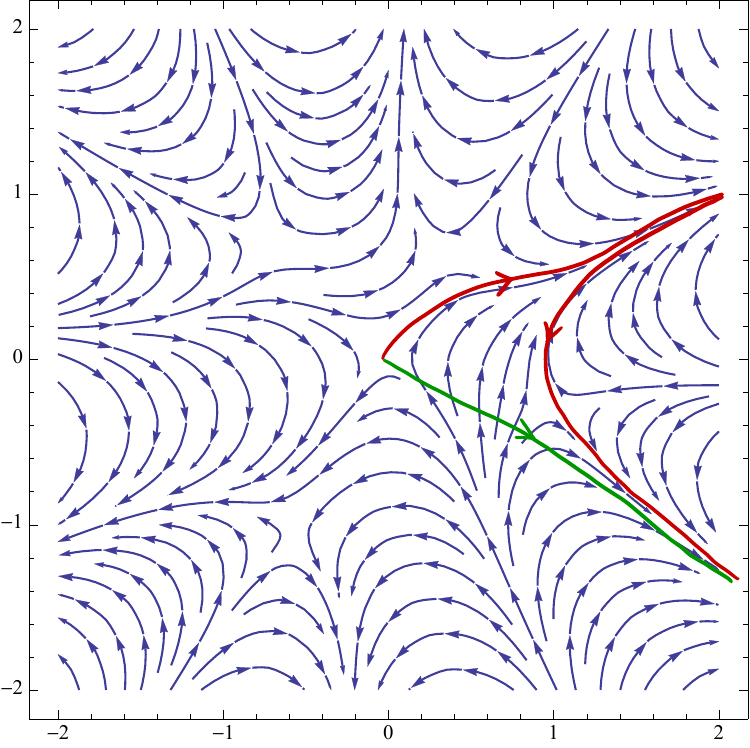}
\includegraphics[width=0.45\textwidth, angle=0]{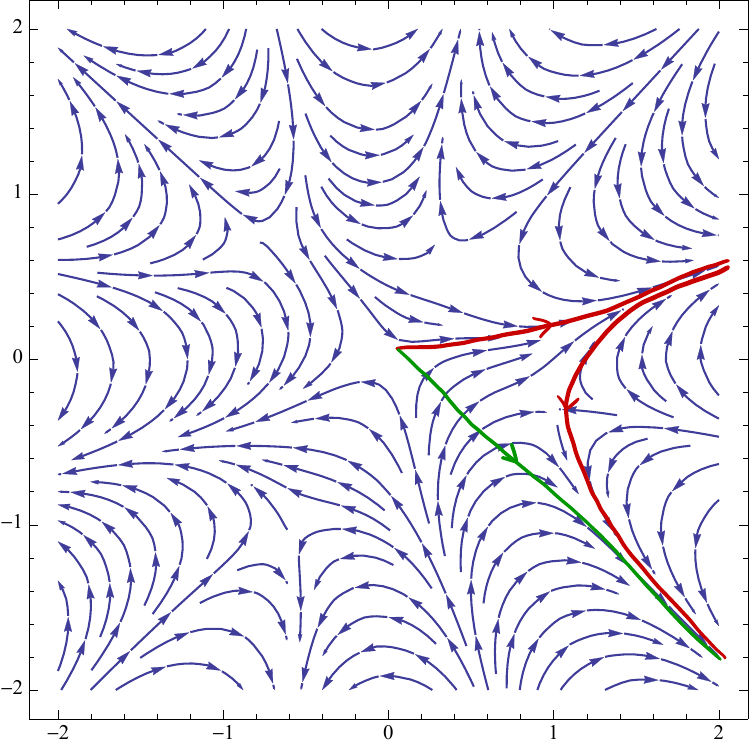}
\includegraphics[width=0.45\textwidth, angle=0]{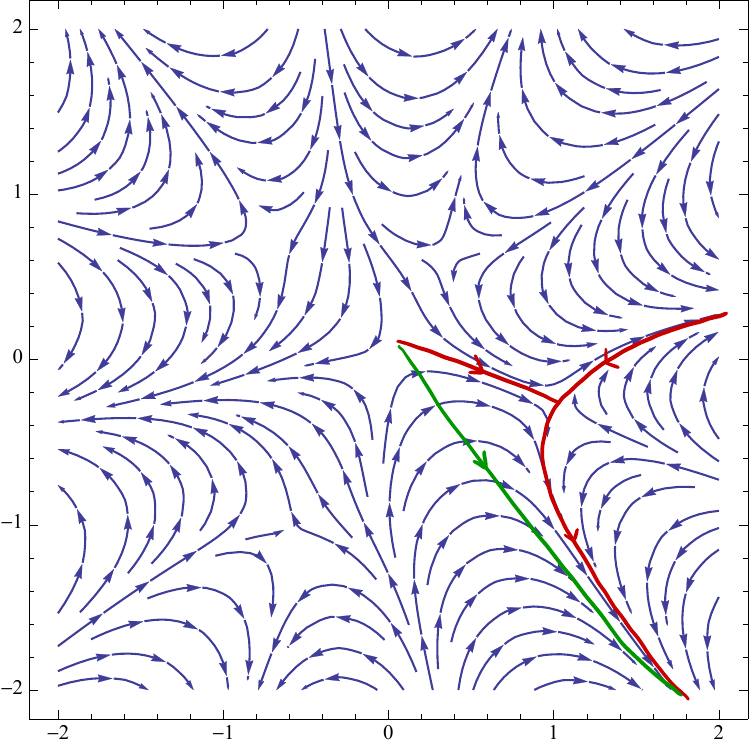}
\caption{Notice that, for a certain range of values of $\arg z$, there is also a contribution in the asymptotic expansion coming from a third critical point. Such a term is negligible, since it is dominated by the exponential term from the critical point $\bt_+$.}
\end{figure}

\begin{prop}\label{A'B'asy}
Let now $\Phi_1,\Phi_2$ be two functions with asymptotic expansions
\beq\label{asP1P2}
\Phi_1(z)\sim z^{D_1}\exp(zu_1),\quad \Phi_2(z)\sim z^{D_2}\exp(zu_2),
\eneq
for $|z|\to\infty$ in the sectors 
\beq\label{zone1}
A_1<\arg z<B_1,\quad A_2<\arg z< B_2,
\eneq
respectively.  We have that
\[\mathscr L_{(1,2;\frac{1}{2},\frac{1}{3})}[\Phi_1,\Phi_2;z]\sim 
Cz^{\frac{1}{2}+D_1+D_2}\exp{z(-\bt_+^6+u_1\bt_+^3+u_2\bt_+^2)},
\]
where
\[C:= 6\bt_+^{5+2D_1+3D_2}\left(\frac{2\pi}{9u_1\bt_++8u_2}\right)^\frac{1}{2},
\]
for $|z|\to\infty$ in the sector $A'<\arg z< B'$, where
\begin{align}
A'&:=\max\left\{A_1-3\arg \bt_+,\  A_2-2\arg \bt_+,\ \arg \overline{f(\bt_+)}-\pi\right\},\\
B'&:=\min\left\{B_1-3\arg \bt_+,\  B_2-2\arg \bt_+,\ \arg \overline{f(\bt_+)}+\pi\right\}.
\end{align}
\end{prop}
\proof
The statement follows by application of the steepest descent path method and Lemma \ref{lemasI}. Notice that the sector $A'<\arg z<B'$ is chosen so that the critical point of the logarithm of the integrand lies in the region \eqref{zone1} of validity of the asymptotic expansions \eqref{asP1P2}.
\endproof

\begin{table}
\[
\begin{array}{|c|c|c|c|}
\hline
u_1&u_2&A'&B'\\
\hline
\hline
2&3&-\pi&\frac{\pi}{3}\\
\hline
2&3e^{\frac{2\pi i}{3}}&-3.71775&0.471036\\
\hline
2&3e^{\frac{4\pi i}{3}}&-1.00423&1.62336\\
\hline
-2&3&-\pi &\frac{\pi }{3}\\
\hline
-2&3e^{\frac{2\pi i}{3}}&-1.00423& -0.706554\\
\hline
-2&3e^{\frac{4\pi i}{3}}&-1.62336& 1.00423\\
\hline
\end{array}
\]
\caption{In this table we represent the values $A'$ and $B'$ predicted in Proposition \ref{A'B'asy} for all possible values of $u_1$ and $u_2$.}\label{table2}
\end{table}

\subsection{Stokes basis of the $qDE$ of $\Fb_1$}\label{stbqdef1}Set 
\beq
s_{ij}:=s_i\otimes p_j\in\mc S(\P^1)\otimes \mc S(\P^2),
\eneq
for $i=1,2$ and $j=1,2,3$. See equation \eqref{baseshp}.
\begin{thm}
 The following linear combinations of the tensors $s_{ij}$'s define a basis of $\mc H$:
\beq
s_{11}-5s_{22}-6s_{23},\quad s_{12}+s_{23},\quad s_{13}-s_{22}-2s_{23},\quad s_{21}-4s_{22}-5s_{23}.
\eneq
\end{thm}
\proof
Define the column vectors 
\begin{itemize}
\item $\bm g=(g_1,g_2)^T$ and $\bm s=(s_1,s_2)^T$, bases of $\mc S(\P^1)$,
\item $\bm h=(h_1,h_2,h_3)^T$ and $\bm p=(p_1,p_2,p_3)^T$, bases of $\mc S(\P^2)$, respectively. 
\end{itemize}
In what follow we denote by $A\otimes B$ the Kronecker tensor product of two matrices $A$ and $B$. Hence we denote
\begin{itemize}
\item by 
$\bm g\otimes \bm h$ the basis $(g_i\otimes h_j)_{i,j}$ of $\mc S(\P^1)\otimes \mc S(\P^2)$.
\item by
$\bm s\otimes \bm p$ the basis $(s_i\otimes p_j)_{i,j}$ of $\mc S(\P^1)\otimes \mc S(\P^2)$.
\end{itemize}
We have
\beq
\label{gh0}
\bm g\otimes \bm h=[(M_1M_2^{-1})\otimes (N_1N_2^{-1})] \bm s\otimes \bm p,
\eneq
where we represent the basis $\bm g\otimes \bm h$ and $\bm s\otimes \bm p$ as column vectors. Multiply on the left both sides of \eqref{gh0} by the matrix 
\[
E_1:=\left(
\begin{array}{cccccc}
 1 & 0 & 0 & 0 & 0 & 0 \\
 0 & 1 & 0 & 0 & 0 & 0 \\
 0 & 0 & 1 & 0 & 0 & 0 \\
 0 & \frac{1}{3} & 0 & \frac{1}{4} & 0 & 0 \\
 0 & 0 & 0 & 0 & 1 & 0 \\
 0 & 0 & 0 & 0 & 0 & 1 \\
\end{array}
\right).
\]
We thus obtain the relation
\beq
\label{gh1}
\bm s\otimes \bm p= X\begin{pmatrix}
g_1\otimes h_1\\
g_1\otimes h_2\\
g_1\otimes h_3\\
\frac{1}{3}g_1\otimes h_2+\frac{1}{4}g_2\otimes h_1\\
g_2\otimes h_2\\
g_2\otimes h_3
\end{pmatrix},
\eneq
where $X$ is the matrix
\begin{align*}X&=[(M_1M_2^{-1})\otimes (N_1N_2^{-1})]^{-1}E_1^{-1}\\
&=
\left(
\begin{array}{cccccc}
 54 & 36 (\gamma +11 i \pi )&*&*&*&* \\
 -54 & -36 (\gamma +i \pi )&*&*&*&* \\
 54 & 36 (\gamma +3 i \pi ) &*&*&*&*\\
 54 & 36 (\gamma +9 i \pi ) &*&*&*&*\\
 -54 & -36 (\gamma -i \pi ) &*&*&*&*\\
 54 & 36 (\gamma +i \pi ) &*&*&*&*\\
\end{array}
\right).
\end{align*}
Multiply on the left each sides of \eqref{gh1} by the matrix
\[E_2=\left(
\begin{array}{cccccc}
 1 & 0 & 0 & 0 & -5 & -6 \\
 0 & 1 & 0 & 0 & 0 & 1 \\
 0 & 0 & 1 & 0 & -1 & -2 \\
 0 & 0 & 0 & 1 & -4 & -5 \\
 0 & 0 & 0 & 0 & 1 & 1 \\
 0 & 0 & 0 & 0 & 0 & 1 \\
\end{array}
\right).
\]
We obtain
\[
\begin{pmatrix}
s_{11}-5s_{22}-6s_{23}\\ s_{12}+s_{23}\\ s_{13}-s_{22}-2s_{23}\\ s_{21}-4s_{22}-5s_{23}\\s_{22}+s_{23}\\s_{23}
\end{pmatrix}
=E_2X\begin{pmatrix}
g_1\otimes h_1\\
g_1\otimes h_2\\
g_1\otimes h_3\\
\frac{1}{3}g_1\otimes h_2+\frac{1}{4}g_2\otimes h_1\\
g_2\otimes h_2\\
g_2\otimes h_3
\end{pmatrix},
\]
and we have
\beq\label{e2x} E_2X=\left(
\begin{array}{cc|cccc}
 0 & 0 &&&&\\
 0 & 0 &&C_1&&\\
 0 & 0 &&&&\\
 0 & 0 &&&&\\
 \hline
 0 & 72 i \pi  &*&*&*&*\\
 54 & 36 (\gamma +i \pi ) &*&*&*&*\\
\end{array}
\right).
\eneq This proves the claim.
\endproof

\begin{rem}\label{c1}
The matrix $C_1$ in equation \eqref{e2x} is
\[
C_1=\left(
\begin{array}{cccc}
 24 (-3 i \gamma -2 \pi ) \pi  & -216 i \pi  & 36 \pi  (-5 i \gamma +9 \pi ) & 3 \pi  \left(-42 i \gamma ^2+92 \gamma  \pi +17 i \pi ^2\right) \\
 72 i \gamma  \pi  & 216 i \pi  & 36 \pi  (5 i \gamma +\pi ) & 3 \pi  \left(42 i \gamma ^2+12 \gamma  \pi -i \pi ^2\right) \\
 -72 i \gamma  \pi  & -216 i \pi  & 36 \pi  (-5 i \gamma +\pi ) & 3 \pi  \left(-42 i \gamma ^2+12 \gamma  \pi +i \pi ^2\right) \\
 -48 \pi ^2 & 0 & 0 & -48 \gamma  \pi ^2 \\
\end{array}
\right).
\]
\end{rem}

\begin{cor}\label{corsig}
The functions 
\begin{align}
\Sig_1:=&\mathscr P(s_{11}-5s_{22}-6s_{23}),\\
\Sig_2:=&\mathscr P(s_{12}+s_{23}),\\
\Sig_3:=&\mathscr P(s_{13}-s_{22}-2s_{23}),\\
\Sig_4:=&\mathscr P(s_{21}-4s_{22}-5s_{23})
\end{align}
define a basis of solutions of the $qDE$ of $\Fb_1$.\qed
\end{cor}

\begin{prop}\label{siglbd}
The Stokes basis $\Xi_R$ of $\mc H^{\rm od}_0$ on the sector $\Pi_R(\eps)$ can be reconstructed, using formulae \eqref{rec1},\eqref{rec2},\eqref{rec3},\eqref{rec4}, from a basis $\bm{\Sig_\lambda}$ of solutions of the $qDE$ of $\Fb_1$ of the form
\[
\lambda_1\Sig_2,\quad \lambda_2\Sig_3+\lambda_3\Sig_2,\quad \lambda_4\Sig_4+\lambda_5\Sig_3+\lambda_6\Sig_2,\quad \lambda_7\Sig_1+\lambda_8\Sig_4+\lambda_9\Sig_3+\lambda_{10}\Sig_2,
\]
for a suitable choice of the coefficients $\lambda_j\in\C$, with $j=1,\dots,10$.
\end{prop}
\proof
The canonical coordinates $x_1,x_2,x_3,x_4$ are in lexicographical order w.r.t. a line of slope $\eps>0$ sufficiently small. The functions above have the expected expnential growth $\exp(x_i z)$ in the sector $\Pi_R(\eps)$ defined by an oriented line of slope $\eps$. This follows from the data in Tables \ref{table1} and \ref{table2}, and from the configuration of the Stokes rays $R_{ij}:=\{-r\sqrt{-1}(\overline{x_i}-\overline{x_j})\colon r\in\R_+\}$: these are given by
\begin{align*}
R_{12}&=\{\arg z=\pi \},& R_{13}&=\{\arg z=2.36573 \},\\
R_{14}&=\{\arg z=1.88197 \},& R_{23}&=\{\arg z=0.775863 \},\\
R_{24}&=\{\arg z=1.25962 \},& R_{34}&=\left\{\arg z=\frac{\pi}{2} \right\},
\end{align*} see Figure \ref{stokes}.
\endproof

\begin{figure}[ht!]
\centering
\def\svgscale{.95}
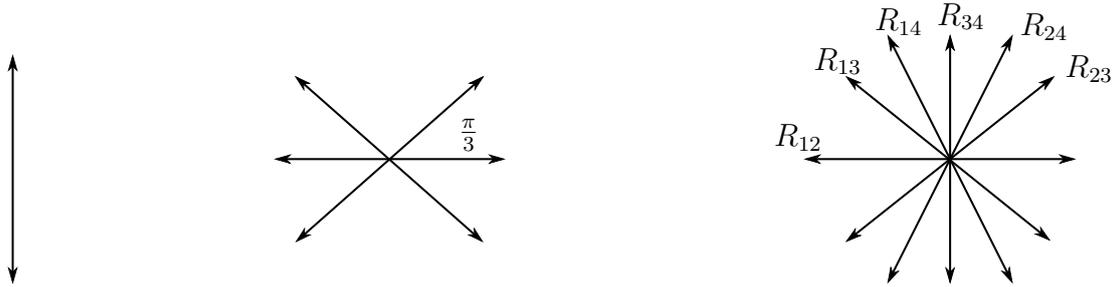
\caption{From the left to the right: Stokes rays corresponding to the origin of the quantum cohomology of $\P^1$, $\P^2$, and $\Fb_1$ respectively. }
\label{stokes}
\end{figure}

\begin{rem}
Notice that, according to Proposition \ref{A'B'asy}, the function $\Sig_3$ has the expected exponential growth $\exp(z x_2)$ in the sector in which this is minimal w.r.t. the dominance relation, i.e. in which it is dominated by any other exponential $\exp(z x_1), \exp(z x_3), \exp(z x_4)$. Hence, we expect that $\lambda_3=0$.
\end{rem}
\subsection{Computation of the central connection and Stokes matrices}
Denote by $\mc H''_0$ the system of differential equations $\mc H^{\rm od}_0$ specialized at $0\in QH^\bullet(\Fb_1)$.
Consider the fundamental system of solutions of $\mc H''_0$
\beq
\Xi_\Ups(z):=\begin{pmatrix}
z\Ups_1(z)&z\Ups_2(z)&z\Ups_3(z)&z\Ups_4(z)\\
\vdots&\vdots&\vdots&\vdots
\end{pmatrix},
\eneq
reconstructed from the basis $(\Ups_1,\Ups_2,\Ups_3,\Ups_4)$ of the $qDE$ of $\Fb_1$ (see Corollary \ref{corUps}) by formulae \eqref{rec1}, \eqref{rec2}, \eqref{rec3}, \eqref{rec4}.

\begin{prop}
We have
\beq
\Xi_\Ups(z)=\Xi_{\rm top}(z)\cdot C_0,
\eneq
where
\beq\label{c0}
C_0:=\begin{pmatrix}
\frac{1}{18}&0&0&0\\
-\frac{\gamma}{18}&\frac{1}{2}&0&0\\
-\frac{\gamma}{18}&-\frac{1}{2}&\frac{1}{3}&0\\
\frac{6\gm^2+\pi^2}{72}&-\frac{\gamma}{2}&-\frac{\gamma}{3}&1
\end{pmatrix}.
\eneq
\end{prop}
\proof
From Lemmata \ref{gi} and \ref{hi}, we can compute the asymptotic expansions of $\Ups_i(z)$ for $z\to 0$. We have
\begin{align*}
\Ups_1(z)=&\frac{1}{72} \left(16 \log ^2(z)-20 \gamma  \log (z)+6 \gamma ^2+\pi ^2\right)+\frac{1}{18} z (\log (z)-\gamma -2)\\
&+\frac{1}{72} z^2 \left(16 \log ^2(z)-20 \gamma  \log (z)-17 \log (z)+6 \gamma ^2+\pi ^2+13 \gamma +2\right)\\
&+\frac{z^3 \left(432 \log ^2(z)-540 \gamma  \log (z)-750 \log (z)+162 \gamma ^2+27 \pi ^2+426 \gamma +311\right)}{1944}\\
&+\dots,\\
\Ups_2(z)=&\frac{1}{2} (\log (z)-\gamma )-\frac{z}{2}+\frac{1}{8} z^2 (4 \log (z)-4 \gamma +5)+\frac{1}{36} z^3 (18 \log (z)-18 \gamma -37)\\
&+\frac{1}{192} z^4 (24 \log (z)-24 \gamma +13)+\dots,\\
\Ups_3(z)=&-\frac{\gamma }{3}+\frac{2 \log (z)}{3}+\frac{z}{3}+\frac{1}{12} z^2 (8 \log (z)-4 \gamma -9)+\frac{1}{54} z^3 (36 \log (z)-18 \gamma -17)\\
&+\frac{1}{288} z^4 (48 \log (z)-24 \gamma -49)+\dots,\\
\Ups_4(z)=&1+z^2+z^3+\frac{z^4}{4}+\dots.
\end{align*}
After some computations, one finds the first terms of the asymptotic expansion of $\Xi_\Ups(z)$ for $z\to 0$:
\begin{align*}
&\Xi_\Ups(z)=\\
& \left(
\begin{array}{cccc}
 \frac{1}{72} z \left(16 \log ^2(z)-20 \gamma  \log (z)+6 \gamma ^2+\pi ^2\right) & \frac{1}{2} z (\log (z)-\gamma ) & z \left(\frac{2 \log (z)}{3}-\frac{\gamma }{3}\right) & z \\
 \frac{\log (z)}{6}-\frac{\gamma }{9} & 0 & \frac{1}{3} & 0 \\
 \frac{\log (z)}{9}+\frac{1}{18} z (\log (z)-\gamma -1)-\frac{\gamma }{18} & \frac{1}{2}-\frac{z}{2} & \frac{z}{3} & 0 \\
 \frac{1}{18} z (2 \log (z)-\gamma -1)+\frac{1}{18 z} & \frac{z}{2} & 0 & 0 \\
\end{array}
\right)\\
\\
&\quad\quad+\text{h.o.t.}.
\end{align*}
The leading term of the asymptotic expansion of $\Xi_{\rm top}(z)$ for $z\to 0$ is 
\begin{align*}
\Xi_{\rm top}(z)&=\eta z^\mu z^R+\text{h.o.t.}\\
&=\left(
\begin{array}{cccc}
 4 z \log ^2(z) & 3 z \log (z) & 2 z \log (z) & z \\
 3 \log (z) & 1 & 1 & 0 \\
 2 \log (z) & 1 & 0 & 0 \\
 \frac{1}{z} & 0 & 0 & 0 \\
\end{array}
\right)+\text{h.o.t.},
\end{align*}
where $\mu={\rm  diag}(-1,0,0,1)$ and $R$ is the operator of $\cup$-multiplication by $c_1(\Fb_1)$ on $H^\bullet(\Fb_1,\C)$, that is
\beq
R=\left(
\begin{array}{cccc}
 0 & 0 & 0 & 0 \\
 2 & 0 & 0 & 0 \\
 1 & 0 & 0 & 0 \\
 0 & 3 & 2 & 0 \\
\end{array}
\right).
\eneq
By comparison of the leading terms of the asymptotic expansions of $\Xi_\Ups$ and $\Xi_{\rm top}$, one obtains the matrix $C_0$ in formula \eqref{c0}.
\endproof

\begin{thm}\label{CSF1}
The central connection and Stokes matrices at $0\in QH^\bullet(\Fb_1)$, computed w.r.t. an admissible oriented line of slope  $\eps>0$ sufficiently small, equal
\beq\label{cf1}
C=\left(
\begin{array}{cccc}
 \frac{1}{2 \pi } & -\frac{1}{2 \pi } & \frac{1}{2 \pi } & -\frac{1}{2 \pi } \\
 \frac{\gamma }{\pi } & -\frac{\gamma }{\pi } & i+\frac{\gamma }{\pi } & -i-\frac{\gamma }{\pi } \\
 \frac{1}{2} \left(-i+\frac{\gamma }{\pi }\right) & -\frac{\gamma +i \pi }{2 \pi } & \frac{1}{2} \left(-i+\frac{\gamma }{\pi }\right) & -\frac{\gamma +i \pi }{2 \pi } \\
 \gamma  \left(-i+\frac{2 \gamma }{\pi }\right) & \gamma  \left(-i-\frac{2 \gamma }{\pi }\right) & \frac{2 \gamma  (\gamma +i \pi )}{\pi } & -\frac{2 (\gamma +i \pi )^2}{\pi } \\
\end{array}
\right),
\eneq
\beq \label{sf1}
S=\left(
\begin{array}{cccc}
 1 & 2 & -1 & -3 \\
 0 & 1 & 1 & -1 \\
 0 & 0 & 1 & 2 \\
 0 & 0 & 0 & 1 \\
\end{array}
\right).
\eneq
\end{thm}
\proof
Denote
\begin{itemize} 
\item by $\Xi_{\bm \lambda}$ the fundamental system of solutions of $\mc H''_0$ constructed from the basis $\bm{\Sig_\lambda}$ of Proposition \ref{siglbd},
\item by $\Xi_{\Sig}$ the fundamental system of solutions of $\mc H''_0$ constructed from the basis $\bm{\Sig}$ of Corollary \ref{corsig}.
\end{itemize} 
We have
\begin{align*}
\Xi_{\bm\lambda}=&\Xi_\Sig\cdot \left(
\begin{array}{cccc}
 0 & 0 & 0 & \lambda _7 \\
 \lambda _1 & \lambda _3 & \lambda _6 & \lambda _{10} \\
 0 & \lambda _2 & \lambda _5 & \lambda _9 \\
 0 & 0 & \lambda _4 & \lambda _8 \\
\end{array}\right)
=\Xi_\Ups\Pi^TC_1^T\left(
\begin{array}{cccc}
 0 & 0 & 0 & \lambda _7 \\
 \lambda _1 & \lambda _3 & \lambda _6 & \lambda _{10} \\
 0 & \lambda _2 & \lambda _5 & \lambda _9 \\
 0 & 0 & \lambda _4 & \lambda _8 \\
\end{array}\right),\\
\end{align*}
where $C_1$ is as in Remark \ref{c1} and 
\[\Pi:=\left(
\begin{array}{cccc}
 0 & 1 & 0 & 0 \\
 1 & 0 & 0 & 0 \\
 0 & 0 & 1 & 0 \\
 0 & 0 & 0 & 1 \\
\end{array}
\right).
\]
Thus, we obtain
\[
\Xi_{\bm \lambda}=\Xi_{\rm top}C_{\bm \lambda},\quad C_{\bm \lambda}:=C_0\Pi^TC_1^T\left(
\begin{array}{cccc}
 0 & 0 & 0 & \lambda _7 \\
 \lambda _1 & \lambda _3 & \lambda _6 & \lambda _{10} \\
 0 & \lambda _2 & \lambda _5 & \lambda _9 \\
 0 & 0 & \lambda _4 & \lambda _8 \\
\end{array}\right),
\]where $C_0$ is given by \eqref{c0}. In order to determine the values of $\bm \lambda$ for which $\Xi_{\bm\lambda}$ is the Stokes basis, let us compute the product
\beq\label{invS} C_{\bm \lambda}^T\eta e^{\pi i\mu}e^{\pi i R}C_{\bm \lambda}.
\eneq If $\Xi_{\bm \lambda}$ is the Stokes basis, then the matrix above is the inverse of the Stokes matrix $S$,  by equation \eqref{const2}: in particular, it is an upper triangular matrix with 1's along the main diagonal. An explicit computation gives the following result: the columns of \eqref{invS} are 
\beq
\left(
\begin{array}{c}
 -576 \pi ^4 \lambda _1^2 \\
 -576 \pi ^4 \lambda _1 \lambda _3 \\
 -576 \pi ^4 \lambda _1 \lambda _6 \\
 -576 \pi ^4 \lambda _1 \left(3 \lambda _7+\lambda _{10}\right) \\
\end{array}
\right),
\eneq
\beq
\left(
\begin{array}{c}
 576 \pi ^4 \lambda _1 \left(2 \lambda _2-\lambda _3\right) \\
 -576 \pi ^4 \left(\lambda _2-\lambda _3\right){}^2 \\
 -576 \pi ^4 \left(\lambda _3 \lambda _6+\lambda _2 \left(\lambda _4+\lambda _5-2 \lambda _6\right)\right) \\
 576 \pi ^4 \left(\lambda _2 \left(\lambda _7-\lambda _8-\lambda _9+2 \lambda _{10}\right)-\lambda _3 \left(3 \lambda _7+\lambda _{10}\right)\right) \\
\end{array}
\right),
\eneq
\beq
\left(
\begin{array}{c}
 -576 \pi ^4 \lambda _1 \left(\lambda _4-2 \lambda _5+\lambda _6\right) \\
 -576 \pi ^4 \left(\lambda _2 \lambda _5+\lambda _3 \left(\lambda _4-2 \lambda _5+\lambda _6\right)\right) \\
 -576 \pi ^4 \left(\lambda _4^2+\left(\lambda _5+\lambda _6\right) \lambda _4+\left(\lambda _5-\lambda _6\right){}^2\right) \\
 -576 \pi ^4 \left(\lambda _6 \left(3 \lambda _7+\lambda _{10}\right)+\lambda _4 \left(5 \lambda _7+\lambda _8+\lambda _{10}\right)+\lambda _5 \left(-\lambda _7+\lambda _8+\lambda _9-2 \lambda _{10}\right)\right) \\
\end{array}
\right),
\eneq
\beq
\left(
\begin{array}{c}
 576 \pi ^4 \lambda _1 \left(6 \lambda _7-\lambda _8+2 \lambda _9-\lambda _{10}\right) \\
 -576 \pi ^4 \left(\lambda _2 \left(6 \lambda _7+\lambda _9\right)+\lambda _3 \left(-6 \lambda _7+\lambda _8-2 \lambda _9+\lambda _{10}\right)\right) \\
 -576 \pi ^4 \left(\lambda _5 \left(6 \lambda _7+\lambda _9\right)+\lambda _4 \left(6 \lambda _7+\lambda _8+\lambda _9\right)+\lambda _6 \left(-6 \lambda _7+\lambda _8-2 \lambda _9+\lambda _{10}\right)\right) \\
 -576 \pi ^4 \left(13 \lambda _7^2+\left(11 \lambda _8+5 \lambda _9-3 \lambda _{10}\right) \lambda _7+\lambda _8^2+\left(\lambda _9-\lambda _{10}\right){}^2+\lambda _8 \left(\lambda _9+\lambda _{10}\right)\right) \\
\end{array}
\right).
\eneq
The matrix \eqref{invS} is upper triangular with 1's along the diagonal if and only if
\begin{align}
\lambda_1^2=&-\frac{1}{576\pi^4}, &\lambda_2^2=&-\frac{1}{576\pi^4},\\
\lambda_3=&\ 0,& \lambda_4^2=&-\frac{1}{576\pi^4},\\
\lambda_5=&-\lambda_4,& \lambda_6=&\ 0,\\
\lambda_7^2=&-\frac{1}{576\pi^4},& \lambda_{8}=&-2\lambda_7,\\
\lambda_9=&-3\lambda_7,& \lambda_{10}=&-3\lambda_7.
\end{align}
For the choice $\lambda_1=\lambda_2=\lambda_4=\lambda_7=-\frac{i}{24\pi^2}$, we obtain the central connection and Stokes matrices \eqref{cf1} and \eqref{sf1}.
\endproof

\begin{thm}\label{DubrconjF2k1}
The central connection matrix of $QH^\bullet(\Fb_{2k+1})$, computed w.r.t. an oriented line of slope $\eps>0$ sufficiently small, and a suitable choice of the branch of the $\Psi$-matrix, equals
\beq\label{matconckf2k1}
C_k=\left(
\begin{array}{cccc}
 \frac{1}{2 \pi } & -\frac{1}{2 \pi } & \frac{1}{2 \pi } & -\frac{1}{2 \pi } \\
 \frac{\gamma }{\pi } & -\frac{\gamma }{\pi } & i+\frac{\gamma }{\pi } & -i-\frac{\gamma }{\pi } \\
 \frac{\gamma -2 \gamma  k-i \pi }{2 \pi } & -\frac{\gamma -2 \gamma  k+i \pi }{2 \pi } & \frac{-2 \gamma  k-i (2 \pi  k+\pi )+\gamma }{2 \pi } & \frac{(2 k-1) (\gamma +i \pi )}{2 \pi } \\
 \gamma  \left(-i+\frac{2 \gamma }{\pi }\right) & \gamma  \left(-i-\frac{2 \gamma }{\pi }\right) & \frac{2 \gamma  (\gamma +i \pi )}{\pi } & -\frac{2 (\gamma +i \pi )^2}{\pi } \\
\end{array}
\right).
\eneq
This is the matrix associated with the morphism 
\[\textnormal{\textcyr{D}}^-_{\mathbb F_{2k+1}}\colon K_0(\mathbb F_{2k+1})_\mathbb C\to H^\bullet(\mathbb F_{2k+1},\mathbb C)\colon [\mathscr F]\mapsto \frac{1}{2\pi}\widehat{\Gamma}^-_{\mathbb F_{2k+1}}\cup e^{-\pi i c_1(\mathbb F_{2k+1})}\cup {\rm Ch}(\mathscr F),
\]w.r.t. an exceptional basis $\frak E:=(E_i)_{i=1}^4$ of $K_0(\mathbb F_{2k+1})_{\mathbb C}$ and the basis $(T_{i,2k+1})_{i=0}^3$ of $H^\bullet(\mathbb F_{2k+1},\mathbb C)$. The exceptional basis $\frak E$  mutates to the exceptional basis
\beq\label{excbasis}
\left([\mathcal O],[\mathcal O(\Sigma_2^{2k+1})],[\mathcal O(\Sigma_4^{2k+1})],[\mathcal O(\Sigma_2^{2k+1}+\Sigma_4^{2k+1})]\right),
\eneq
by application of the following natural transformations:
\begin{enumerate}
\item  action of the braid $\bt_{3}\bt_{2}\bt_{1}\bt_{3}\bt_{2}$;
\item action of the element $\tilde J_k\in(\Z/2\Z)^4$
\begin{empheq}[left={\tilde J_k:=} \empheqlbrace]{align*} 
&(-1,-1,(-1)^p,(-1)^{p+1}),\quad\text{if }k=2p,\\
\\
&(-1,-1,(-1)^{p+1},(-1)^{p+1}),\quad\text{if }k=2p+1;
\end{empheq}
\item action of the element $\bt_{3}^k$.
\end{enumerate}
\end{thm}

\proof
Equations \eqref{chcordubro1} and Proposition \ref{flatcord} imply equation \eqref{matconckf2k1}.
The matrix associated to $\textnormal{\textcyr{D}}^-_{\mathbb F_{2k+1}}$ w.r.t. the basis \eqref{excbasis} is
\[
E_k:=\left(
\begin{array}{cccc}
 \frac{1}{2 \pi } & \frac{1}{2 \pi } & \frac{1}{2 \pi } & \frac{1}{2 \pi } \\
 -i+\frac{\gamma }{\pi } & -i+\frac{\gamma }{\pi } & \frac{\gamma }{\pi } & \frac{\gamma }{\pi } \\
 \frac{(1-2 k) (\gamma -i \pi )}{2 \pi } & \frac{-2 \gamma  k+i (2 \pi  k+\pi )+\gamma }{2 \pi } & \frac{(1-2 k) (\gamma -i \pi )}{2 \pi } & \frac{-2 \gamma  k+i (2 \pi  k+\pi )+\gamma }{2 \pi } \\
 \frac{2 (\gamma -i \pi )^2}{\pi } & \frac{2 \gamma  (\gamma -i \pi )}{\pi } & \gamma  \left(-i+2 i k+\frac{2 \gamma }{\pi }\right) & \gamma  \left(i+2 i k+\frac{2 \gamma }{\pi }\right) \\
\end{array}
\right).
\]
Set $C_k':=C_k^{\bt_{3}\bt_{2}\bt_{1}\bt_{3}\bt_{2}}$. We have
\[(C'_k)^{-1}E_k=\left(
\begin{array}{cccc}
 -1 & 0 & 0 & 0 \\
 0 & -1 & 0 & 0 \\
 0 & 0 & 1-k & -k \\
 0 & 0 & -k & -k-1 \\
\end{array}
\right).
\]
It is easy to see that this is the matrix representing the action of the element $(\tilde J_k,\bt_3^{k})\in (\Z/2\Z)^4\rtimes\mc B_4$: the argument is the same as in Step 3 of the proof of Theorem \ref{DubconjF2k}.
\endproof

\newpage
\appendix
\section{Proof of Theorem \ref{topJ}}\label{aproof}
We need some preliminary results.

\begin{lem}\label{lem1}
For $n\geq 0$, and $\delta\in H^2(X,\C)$, we have
\[\llangle \tau_n T_\alpha,1\rrangle_0(\delta)=\frac{1}{(n+1)!}\left(\int_XT_\al\cup\delta^{n+1}\right)+\sum_{\bt\neq 0}\sum_{\nu\geq 0}\frac{{\bf Q}^\bt e^{\int_\bt\delta}}{\nu!}\langle\tau_{n-\nu}T_{\al}\cup\delta^\nu,1\rangle^X_{0,2,\bt}.
\]
\end{lem}
\proof
We have
\begin{align*}
\llangle \tau_n T_\alpha,1\rrangle_0(\delta)&=\left.\frac{\partial}{\partial t^{\al,n}}\frac{\partial}{\partial t^{0,0}}\mc F^X_0\right|_\delta=\sum_{k=0}^\infty\sum_\bt\frac{{\bf Q}^\bt}{k!}\langle\tau_nT_{\al},1,\delta\,\dots,\delta \rangle^X_{0,k+2,\bt}.
\end{align*}
We have two cases:
\begin{itemize}
\item if $\bt\neq 0$, then for $k\geq 0$  we have
\[
\langle\tau_nT_{\al},1, \delta\,\dots, \delta\rangle^X_{0,k+2,\bt}=\sum_{\mu+\nu=k}\frac{k!}{\mu!\nu!}\left(\int_\bt\delta\right)^\mu\langle\tau_{n-\nu}T_{\al}\cup\delta^\nu,1\rangle^X_{0,2,\bt},
\]
by the Divisor Axiom of Gromov-Witten invariants. Here any invariant with $\tau_{-r}$ with $r>0$ is vanishing.
\item If $\bt =0$, then for $k>0$ by Divisor Axiom we have\footnote{Here, we use the fact that $\mc L_1$ is trivial on $\overline{\mc M}_{0,3}(X,0)$ and hence has zero Chern class. This follows from the fact that $\overline{\mc M}_{0,3}(X,0)\cong X$, and the frogetful morphism $\overline{\mc M}_{0,4}(X,0)\to\overline{\mc M}_{0,3}(X,0)$ is the projection $X\times \overline{\mc M}_{0,4}\to X$.}
\[\langle\tau_nT_{\al},1, \delta\,\dots, \delta\rangle^X_{0,k+2,0}=\langle\tau_{n-k+1}T_\al\cup \delta^k,1,\delta\rangle_{0,3,0}=\left(\int_XT_\al\cup\delta^k\right)\delta_{k,n+1}.
\]
\end{itemize}
So, we have
\begin{align*}
\llangle \tau_n T_\alpha,1\rrangle_0(\delta)=&\frac{1}{(n+1)!}\left(\int_XT_\al\cup\delta^{n+1}\right)\\
&+\sum_{\bt\neq 0}\sum_{k\geq 0}\frac{{\bf Q^\bt}}{k!}\sum_{\mu+\nu=k}\frac{k!}{\mu!\nu!}\left(\int_\bt\delta\right)^\mu\langle\tau_{n-\nu}T_{\al}\cup\delta^\nu,1\rangle^X_{0,2,\bt}\\
=&\frac{1}{(n+1)!}\left(\int_XT_\al\cup\delta^{n+1}\right)+\sum_{\bt\neq 0}\sum_{\nu\geq 0}\frac{{\bf Q}^\bt e^{\int_\bt\delta}}{\nu!}\langle\tau_{n-\nu}T_{\al}\cup\delta^\nu,1\rangle^X_{0,2,\bt}. \quad\quad
\puqed
\qedhere
\poqed
\end{align*}

\begin{lem}\label{restJ}
Let $\delta\in H^2(X,\C)$. We have
\[
J_X(\delta)=e^\frac{\delta}{\hbar}+\sum_\al\sum_{\bt\neq0}\sum_{n=0}^\infty\sum_{k+p=n}\hbar^{-(n+1)}\frac{{\bf Q}^\bt e^{\int_\bt\delta}}{p!}\langle\tau_{k}T_{\al}\cup\delta^p,1\rangle^X_{0,2,\bt}T^\al.
\]
\end{lem}
\proof 
By Lemma \ref{lem1}, we have
\begin{align*}
J_X(\delta)=&1+\sum_\al\sum_{n=0}^\infty\frac{\hbar^{-(n+1)}}{(n+1)!}\left(\int_XT_\al\cup\delta^{n+1}\right)T^\al\\
&+\sum_\al\sum_{\bt\neq0}\sum_{n=0}^\infty\sum_{k+p=n}\hbar^{-(n+1)}\frac{{\bf Q}^\bt e^{\int_\bt\delta}}{p!}\langle\tau_{k}T_{\al}\cup\delta^p,1\rangle^X_{0,2,\bt}T^\al\\
\puqed
=&e^\frac{\delta}{\hbar}+\sum_\al\sum_{\bt\neq0}\sum_{n=0}^\infty\sum_{k+p=n}\hbar^{-(n+1)}\frac{{\bf Q}^\bt e^{\int_\bt\delta}}{p!}\langle\tau_{k}T_{\al}\cup\delta^p,1\rangle^X_{0,2,\bt}T^\al.\qedhere
\poqed
\end{align*}

\begin{lem}\label{lem3}
For $\delta\in H^2(X,\C)$, we have
\beq
\label{ztopr}
Z_{\rm top}(\delta, z)T_\al=e^{z\delta}\cup z^\mu z^{c_1(X)}T_\alpha+\sum_{\beta\neq 0}\sum_{\lambda}e^{\int_\beta\delta}\bigg\langle\frac{ze^{z\delta}}{1-z\psi}\cup z^\mu z^{c_1(X)}T_\alpha,T_\lambda\bigg\rangle^X_{0,2,\beta}T^\lambda.
\eneq
\end{lem}
\proof
For $\bm \tau\in H^\bullet(X,\C)$, we have
\begin{align*}
\Theta(\bm \tau,z)T_\al=&\sum_\eps\Theta(\bm \tau,z)_\al^\eps T_\eps=\sum_\la\left.\frac{\der\theta_\al}{\der t^\la}\right|_{(\bm \tau,z)}T^\la\\
=&\sum_\la\sum_{p=0}^\infty z^p\llangle\tau_pT_\al,1, T_\la\rrangle_0(\bm \tau)|_{\bf Q=1}T^\la\\
=&\sum_\la\sum_{p=0}^\infty\sum_{k=0}^\infty\sum_\bt\frac{z^p}{k!}\langle\tau_pT_\al,1,T_\la,\bm \tau,\dots,\bm \tau\rangle^X_{0,3+k,\bt}T^\la.
\end{align*}
Consider the contribution coming from $(k,\bt)=(0,0)$: by the Mapping to point Axiom of Gromov-Witten invariants, we have\footnote{Also here, we use the fact that $\mc L_1$ is trivial on $\overline{\mc M}_{0,3}(X,0)$.}
\begin{align*}
\sum_\la\sum_{p=0}^\infty{z^p}\langle\tau_pT_\al,1,T_\la,\rangle^X_{0,3,0}T^\la=\sum_\la\sum_{p=0}^\infty z^p\left(\int_XT_\al\cup T_\la\right)\delta_{0,p}T^\la=T_\al.
\end{align*}
By the Fundamental class Axiom, instead, the contribution from $(k,\bt)\neq (0,0)$ can be re-written as
\begin{align*}
\sum_\la\sum_{p=0}^\infty\sum_{k=1}^\infty\sum_{\bt\neq 0}\frac{z^p}{k!}\langle\tau_{p-1}T_\al,T_\la,\bm \tau,\dots,\bm \tau\rangle^X_{0,2+k,\bt}T^\la.
\end{align*}
Thus, we have recovered the formula 
\[\Theta(\bm \tau,z)={\rm Id}+\sum_\la\sum_{p=0}^\infty z^{p+1}\llangle\tau_p(-),T_\la\rrangle_0(\bm \tau)|_{\bf Q=1}T^\la,
\]which was used in \cite[Proposition 7.1]{CDG} to define $\Theta$. At this point the proof is known, and can be found in \cite[Proposition 10.2.3]{cox}: the parameter $\hbar$ of \emph{loc. cit.} has to be replaced by our $z$, and pre-composition with $z^\mu z^{c_1(X)}$ has to be taken into account in order to obtain formula \eqref{ztopr}.
\endproof

We are now ready for the proof of Theorem \ref{topJ}.
\proof[Proof of Theorem \ref{topJ}]
Let us compute the entries of the first row of the matrix
\[\eta\Theta(\delta,z)z^\mu z^{c_1(X)}.
\]By Lemma \ref{lem3}, we have
\begin{align*}
&\left[\eta\Theta(\delta,z)z^\mu z^{c_1(X)}\right]^1_\alpha=\eta\left(1,\Theta(\delta,z)z^\mu z^{c_1(X)}T_\alpha\right)\\
&=\eta\left(1, e^{z\delta}\cup z^\mu z^{c_1(X)}T_\alpha+\sum_{\beta\neq 0}\sum_{\lambda}e^{\int_\beta\delta}\bigg\langle\frac{ze^{z\delta}}{1-z\psi}\cup z^\mu z^{c_1(X)}T_\alpha,T_\lambda\bigg\rangle^X_{0,2,\beta}T^\lambda\right)\\
&=\eta(1,e^{z\delta}\cup z^\mu z^{c_1(X)}T_\alpha)\\
&+\eta\left(1,\sum_{\beta\neq 0}\sum_{\lambda}e^{\int_\beta\delta}\bigg\langle\frac{ze^{z\delta}}{1-z\psi}\cup z^\mu z^{c_1(X)}T_\alpha,T_\lambda\bigg\rangle^X_{0,2,\beta}T^\lambda\right).
\end{align*}
Using the identity of endomorphisms of $H^\bullet(X,\mathbb C)$ 
\[z^{-\mu}\circ (h^k\cup)\circ z^\mu = z^{-k} (h^k\cup),\quad h\in H^2(X,\mathbb C),\quad k\in\mathbb N,
\]
and the $\eta$-skew-symmetry of $\mu$, we can rewrite the first summand as 
\begin{align*}
\eta(1,e^{z\delta}\cup z^\mu z^{c_1(X)}T_\alpha)&=\eta(1, z^\mu e^\delta z^{c_1(X)}T_\alpha)\\
&=\eta(z^{-\mu}(1),e^\delta z^{c_1(X)}T_\alpha)\\
&=z^{\frac{\dim_\mathbb CX}{2}}\int_Xe^\delta z^{c_1(X)}T_\alpha.
\end{align*}
For the second summand, notice that  
\begin{enumerate}
\item the only nonzero contribution comes from $\lambda =0$, 
\item for any $\phi\in H^\bullet(X,\mathbb C)$ we have that
\[\frac{ze^{z\delta}}{1-z\psi}\cup\phi=\sum_{n=0}^\infty\sum_{k=0}^n\frac{z^{n+1}}{(n-k)!}\psi^k\delta^{n-k}\phi,
\]
\item and that
\[z^\mu z^{c_1(X)}T_\alpha=\sum_{\ell=0}^\infty\frac{(\log z)^\ell}{\ell!}z^{\frac{2\ell+\deg T_\alpha-\dim X}{2}}c_1(X)^\ell T_\alpha.
\]
\item the Gromov-Witten invariant 
\[\langle\tau_k\delta^{n-k}c_1(X)^\ell T_\alpha,1\rangle^X_{0,2,\beta}
\]is nonzero 
only if
\[2k+2(n-k)+2\ell+\deg T_\alpha=2\dim_{\mathbb C}X+2\int_\beta c_1(X)-2.
\]
\end{enumerate}
So, we obtain that
\begin{align*}&\bigg\langle\frac{ze^{z\delta}}{1-z\psi}\cup z^\mu z^{c_1(X)}T_\alpha,1\bigg\rangle^X_{0,2,\beta}&\\
&=\sum_{n=0}^\infty\sum_{k=0}^n\sum_{\ell=0}^\infty\frac{(\log z)^\ell}{\ell!(n-k)!}z^{n+1+\frac{2\ell+\deg T_\alpha-\dim X}{2}}\langle \tau_k\delta^{n-k}c_1(X)^\ell T_\alpha, 1\rangle^X_{0,2,\beta}\\
&=z^{\frac{\dim X}{2}}z^{\int_\beta c_1(X)}\sum_{h=0}^\infty\sum_{m+\ell+k=h}\frac{(\log z)^\ell}{\ell!m!}\langle \tau_k\delta^mc_1(X)^\ell T_\alpha, 1\rangle^X_{0,2,\beta}\\
&=z^{\frac{\dim X}{2}}z^{\int_\beta c_1(X)}\sum_{h=0}^\infty\sum_{k+p=h}\frac{1}{p!}\langle \tau_k(\delta+\log z\cdot c_1(X))^p T_\alpha, 1\rangle^X_{0,2,\beta}.
\end{align*}
Putting this all together, we obtain that
\begin{align*}&\left[\eta\Theta(\delta,z)z^\mu z^{c_1(X)}\right]^1_\alpha\\
&=z^\frac{\dim X}{2}\left(\int_X e^\delta z^{c_1(X)}T_\alpha+\sum_{\beta\neq 0}e^{\int_\beta\delta}z^{\int_\beta c_1(X)}\sum_{h=0}^\infty\sum_{k+p=h}\frac{1}{p!}\langle \tau_k(\delta+\log z\cdot c_1(X))^p T_\alpha, 1\rangle^X_{0,2,\beta}\right)
\\&=z^\frac{\dim X}{2}\int_XT_\alpha\cup J_X(\delta+\log z\cdot c_1(X))\Big|_{\substack{{\bf Q}=1,\\ \hbar=1}}.
\end{align*}
The last equality follows by Lemma \ref{restJ}. This completes the proof. \endproof

\newpage
\section{Coefficients $\mathcal A_j^{(i)},\mathcal B_j^{(i)}$}\label{appa}
The coefficients $\mathcal A_j^{(i)},\mathcal B_j^{(i)}$, introduced in equations \eqref{theta} and \eqref{lambda}, are 
\begin{align*}
\mathcal A_1^{(1)}(m,n)=&-\frac{8}{9} m n^2 A_{0,0} B_{0,0},\\
\mathcal A_2^{(1)}(m,n)=&\ \frac{8}{9} n^2 A_{0,0} B_{0,0},\\
\mathcal A_3^{(1)}(m,n)=&\ \frac{8}{9} m A_{0,0} B_{0,0},
\end{align*}
\begin{align*}
\mathcal A_1^{(2)}(m,n)=&\ \frac{4}{9} n \left(4 m n A_{0,0} B_{0,0} H_m+6 m n A_{0,0} B_{0,0} H_n\right.\\
&-4 m n A_{0,1} B_{0,0}-3 m n A_{0,0} B_{0,1}-4 m n A_{0,0} B_{0,0} \psi ^{(0)}(m+n+1)\\
&\left.-4 m A_{0,0} B_{0,0}-2 n A_{0,0} B_{0,0}\right),\\
\mathcal A_2^{(2)}(m,n)=&-\frac{4}{9} n \left(4 n A_{0,0} B_{0,0} H_m+6 n A_{0,0} B_{0,0} H_n\right.\\
&\left.-4 n A_{0,0} B_{0,0} \psi ^{(0)}(m+n+1)-4 n A_{0,1} B_{0,0}-3 n A_{0,0} B_{0,1}-4 A_{0,0} B_{0,0}\right),\\
\mathcal A_3^{(2)}(m,n)=&- \frac{4}{9}  \left(6 m A_{0,0} B_{0,0} H_n+4 m A_{0,0} B_{0,0} H_m\right.\\
&\left.-4 m A_{0,0} B_{0,0} \psi ^{(0)}(m+n+1)-4 m A_{0,1} B_{0,0}-3 m A_{0,0} B_{0,1}-2 A_{0,0} B_{0,0}\right),
\end{align*}
\begin{align*}
\mathcal A_1^{(3)}(m,n)=&-\frac{2}{9} \left(24 m n^2 A_{0,0} B_{0,0} H_m H_n-24 m n^2 A_{0,1} B_{0,0} H_n-12 m n^2 A_{0,0} B_{0,1} H_m\right.\\
&-8 m n^2 A_{0,0} B_{0,0} H_m \psi ^{(0)}(m+n+1)-18 m n^2 A_{0,0} B_{0,0} H_n \psi ^{(0)}(m+n+1)\\
&-16 m n A_{0,0} B_{0,0} H_m-12 m n A_{0,0} B_{0,0} H_n-12 n^2 A_{0,0} B_{0,0} H_n\\
&+12 m n^2 A_{0,1} B_{0,1}+9 m n^2 A_{0,0} B_{n,2}+5 m n^2 A_{0,0} B_{0,0} \psi ^{(0)}(m+n+1)^2\\
&+4 n^2 A_{0,0} B_{0,0} \psi ^{(0)}(m+n+1)+5 m n^2 A_{0,0} B_{0,0} \psi ^{(1)}(m+n+1)\\
&+8 m n^2 A_{0,1} B_{0,0} \psi ^{(0)}(m+n+1)+9 m n^2 A_{0,0} B_{0,1} \psi ^{(0)}(m+n+1)\\
&+16 m n A_{0,1} B_{0,0}+6 m n A_{0,0} B_{0,1}+12 m n A_{0,0} B_{0,0} \psi ^{(0)}(m+n+1)\\
&+2 m A_{0,0} B_{0,0}+6 n^2 A_{0,0} B_{0,1}+8 n A_{0,0} B_{0,0}\left.\right),
\end{align*}
\begin{align*}
\mathcal A_2^{(3)}(m,n)=&\ \frac{2}{9} \left(24 n^2 A_{0,0} B_{0,0} H_m H_n-12 n^2 A_{0,0} B_{0,1} H_m\right.\\
&-8 n^2 A_{0,0} B_{0,0} H_m \psi ^{(0)}(m+n+1)-18 n^2 A_{0,0} B_{0,0} H_n \psi ^{(0)}(m+n+1)\\
&-16 n A_{0,0} B_{0,0} H_m-24 n^2 A_{0,1} B_{0,0} H_n-12 n A_{0,0} B_{0,0} H_n\\
&+5 n^2 A_{0,0} B_{0,0} \psi ^{(0)}(m+n+1)^2+5 n^2 A_{0,0} B_{0,0} \psi ^{(1)}(m+n+1)\\
&+8 n^2 A_{0,1} B_{0,0} \psi ^{(0)}(m+n+1)+9 n^2 A_{0,0} B_{0,1} \psi ^{(0)}(m+n+1)\\
&+12 n A_{0,0} B_{0,0} \psi ^{(0)}(m+n+1)+12 n^2 A_{0,1} B_{0,1}+9 n^2 A_{0,0} B_{n,2}\\
&+16 n A_{0,1} B_{0,0}+6 n A_{0,0} B_{0,1}+2 A_{0,0} B_{0,0}\left.\right),
\end{align*}
\begin{align*}
\mathcal A_3^{(3)}(m,n)=&\ \frac{2}{9} \left(24 m A_{0,0} B_{0,0} H_m H_n-24 m A_{0,1} B_{0,0} H_n\right.\\
&-8 m A_{0,0} B_{0,0} H_m \psi ^{(0)}(m+n+1)-18 m A_{0,0} B_{0,0} H_n \psi ^{(0)}(m+n+1)\\
&-12 m A_{0,0} B_{0,1} H_m-12 A_{0,0} B_{0,0} H_n+9 m A_{0,0} B_{n,2}\\
&+5 m A_{0,0} B_{0,0} \psi ^{(0)}(m+n+1)^2+4 A_{0,0} B_{0,0} \psi ^{(0)}(m+n+1)\\
&+8 m A_{0,1} B_{0,0} \psi ^{(0)}(m+n+1)+9 m A_{0,0} B_{0,1} \psi ^{(0)}(m+n+1)\\
&+5 m A_{0,0} B_{0,0} \psi ^{(1)}(m+n+1)+12 m A_{0,1} B_{0,1}+6 A_{0,0} B_{0,1}\left.\right),
\end{align*}
\begin{align*}
\mathcal A_1^{(4)}(m,n)=&-\frac{2}{9} \left(-18 m n^2 A_{0,0} H_m B_{n,2}-2 m n^2 A_{0,0} B_{0,0} H_m \psi ^{(0)}(m+n+1)^2\right.\\
&-6 m n^2 A_{0,0} B_{0,0} H_n \psi ^{(0)}(m+n+1)^2-6 n^2 A_{0,0} B_{0,0} H_n \psi ^{(0)}(m+n+1)\\
&+12 m n^2 A_{0,0} B_{0,0} H_m H_n \psi ^{(0)}(m+n+1)\\
&-12 m n^2 A_{0,1} B_{0,0} H_n \psi ^{(0)}(m+n+1)-6 m n^2 A_{0,0} B_{0,1} H_m \psi ^{(0)}(m+n+1)\\
&-2 m n^2 A_{0,0} B_{0,0} H_m \psi ^{(1)}(m+n+1)-6 m n^2 A_{0,0} B_{0,0} H_n \psi ^{(1)}(m+n+1)\\
&+24 m n A_{0,0} B_{0,0} H_m H_n-24 m n A_{0,1} B_{0,0} H_n-12 m n A_{0,0} B_{0,1} H_m\\
&-8 m n A_{0,0} B_{0,0} H_m \psi ^{(0)}(m+n+1)-12 m n A_{0,0} B_{0,0} H_n \psi ^{(0)}(m+n+1)\\
&-4 m A_{0,0} B_{0,0} H_m-12 n A_{0,0} B_{0,0} H_n+18 m n^2 A_{0,1} B_{n,2}\\
&+m n^2 A_{0,0} B_{0,0} \psi ^{(0)}(m+n+1)^3+n^2 A_{0,0} B_{0,0} \psi ^{(0)}(m+n+1)^2\\
&+2 m n^2 A_{0,1} B_{0,0} \psi ^{(0)}(m+n+1)^2+3 m n^2 A_{0,0} B_{0,1} \psi ^{(0)}(m+n+1)^2
\end{align*}
\begin{align*}
&+3 m n^2 A_{0,0} B_{0,0} \psi ^{(1)}(m+n+1) \psi ^{(0)}(m+n+1)\\
&+3 n^2 A_{0,0} B_{0,1} \psi ^{(0)}(m+n+1)+6 m n^2 A_{0,1} B_{0,1} \psi ^{(0)}(m+n+1)\\
&+9 m n^2 A_{0,0} B_{n,2} \psi ^{(0)}(m+n+1)+n^2 A_{0,0} B_{0,0} \psi ^{(1)}(m+n+1)\\
&+m n^2 A_{0,0} B_{0,0} \psi ^{(2)}(m+n+1)+2 m n^2 A_{0,1} B_{0,0} \psi ^{(1)}(m+n+1)\\
&+3 m n^2 A_{0,0} B_{0,1} \psi ^{(1)}(m+n+1)+12 m n A_{0,1} B_{0,1}\\
&+4 m n A_{0,0} B_{0,0} \psi ^{(0)}(m+n+1)^2+2 m A_{0,0} B_{0,0} \psi ^{(0)}(m+n+1)\\
&+4 n A_{0,0} B_{0,0} \psi ^{(0)}(m+n+1)+8 m n A_{0,1} B_{0,0} \psi ^{(0)}(m+n+1)\\
&+6 m n A_{0,0} B_{0,1} \psi ^{(0)}(m+n+1)+4 m n A_{0,0} B_{0,0} \psi ^{(1)}(m+n+1)\\
&+4 m A_{0,1} B_{0,0}+9 n^2 A_{0,0} B_{n,2}+6 n A_{0,0} B_{0,1}+2 A_{0,0} B_{0,0}\left.\right),
\end{align*}
\begin{align*}
\mathcal A_2^{(4)}(m,n)=&\ \frac{2}{9} \left(-18 n^2 A_{0,0} H_m B_{n,2}\right.\\
&-2 n^2 A_{0,0} B_{0,0} H_m \psi ^{(0)}(m+n+1)^2-6 n^2 A_{0,0} B_{0,0} H_n \psi ^{(0)}(m+n+1)^2\\
&+12 n^2 A_{0,0} B_{0,0} H_m H_n \psi ^{(0)}(m+n+1)-12 n^2 A_{0,1} B_{0,0} H_n \psi ^{(0)}(m+n+1)\\
&-6 n^2 A_{0,0} B_{0,1} H_m \psi ^{(0)}(m+n+1)-2 n^2 A_{0,0} B_{0,0} H_m \psi ^{(1)}(m+n+1)\\
&-6 n^2 A_{0,0} B_{0,0} H_n \psi ^{(1)}(m+n+1)+24 n A_{0,0} B_{0,0} H_m H_n-12 n A_{0,0} B_{0,1} H_m\\
&-8 n A_{0,0} B_{0,0} H_m \psi ^{(0)}(m+n+1)-12 n A_{0,0} B_{0,0} H_n \psi ^{(0)}(m+n+1)\\
&-4 A_{0,0} B_{0,0} H_m-24 n A_{0,1} B_{0,0} H_n+n^2 A_{0,0} B_{0,0} \psi ^{(0)}(m+n+1)^3\\
&+2 n^2 A_{0,1} B_{0,0} \psi ^{(0)}(m+n+1)^2+3 n^2 A_{0,0} B_{0,1} \psi ^{(0)}(m+n+1)^2\\
&+3 n^2 A_{0,0} B_{0,0} \psi ^{(1)}(m+n+1) \psi ^{(0)}(m+n+1)\\
&+6 n^2 A_{0,1} B_{0,1} \psi ^{(0)}(m+n+1)+9 n^2 A_{0,0} B_{n,2} \psi ^{(0)}(m+n+1)\\
&+n^2 A_{0,0} B_{0,0} \psi ^{(2)}(m+n+1)+2 n^2 A_{0,1} B_{0,0} \psi ^{(1)}(m+n+1)\\
&+3 n^2 A_{0,0} B_{0,1} \psi ^{(1)}(m+n+1)\\
&+4 n A_{0,0} B_{0,0} \psi ^{(0)}(m+n+1)^2+2 A_{0,0} B_{0,0} \psi ^{(0)}(m+n+1)\\
&+8 n A_{0,1} B_{0,0} \psi ^{(0)}(m+n+1)+6 n A_{0,0} B_{0,1} \psi ^{(0)}(m+n+1)\\
&+4 n A_{0,0} B_{0,0} \psi ^{(1)}(m+n+1)+18 n^2 A_{0,1} B_{n,2}+12 n A_{0,1} B_{0,1}+4 A_{0,1} B_{0,0}\left.\right),
\end{align*}
\begin{align*}
\mathcal A_3^{(4)}(m,n)=&-\frac{2}{9} \left(18 m A_{0,0} H_m B_{n,2}\right.\\
&+2 m A_{0,0} B_{0,0} H_m \psi ^{(0)}(m+n+1)^2+6 m A_{0,0} B_{0,0} H_n \psi ^{(0)}(m+n+1)^2\\
&-12 m A_{0,0} B_{0,0} H_m H_n \psi ^{(0)}(m+n+1)+6 A_{0,0} B_{0,0} H_n \psi ^{(0)}(m+n+1)\\
&+12 m A_{0,1} B_{0,0} H_n \psi ^{(0)}(m+n+1)+6 m A_{0,0} B_{0,1} H_m \psi ^{(0)}(m+n+1)\\
&+2 m A_{0,0} B_{0,0} H_m \psi ^{(1)}(m+n+1)+6 m A_{0,0} B_{0,0} H_n \psi ^{(1)}(m+n+1)\\
&-18 m A_{0,1} B_{n,2}-m A_{0,0} B_{0,0} \psi ^{(0)}(m+n+1)^3-A_{0,0} B_{0,0} \psi ^{(0)}(m+n+1)^2\\
&-2 m A_{0,1} B_{0,0} \psi ^{(0)}(m+n+1)^2-3 m A_{0,0} B_{0,1} \psi ^{(0)}(m+n+1)^2\\
&-3 m A_{0,0} B_{0,0} \psi ^{(1)}(m+n+1) \psi ^{(0)}(m+n+1)-3 A_{0,0} B_{0,1} \psi ^{(0)}(m+n+1)\\
&-6 m A_{0,1} B_{0,1} \psi ^{(0)}(m+n+1)-9 m A_{0,0} B_{n,2} \psi ^{(0)}(m+n+1)\\
&-A_{0,0} B_{0,0} \psi ^{(1)}(m+n+1)-m A_{0,0} B_{0,0} \psi ^{(2)}(m+n+1)\\
&-2 m A_{0,1} B_{0,0} \psi ^{(1)}(m+n+1)-3 m A_{0,0} B_{0,1} \psi ^{(1)}(m+n+1)-9 A_{0,0} B_{n,2}\left.\right),
\end{align*}
\begin{align*}
\mathcal B_1^{(1)}(m,n)=&\ \frac{2}{9} n A_{0,0} B_{0,0} (m-n),\\
\mathcal B_2^{(1)}(m,n)=&- \frac{2}{9} n A_{0,0} B_{0,0},\\
\mathcal B_3^{(1)}(m,n)=&\ \frac{2}{9} A_{0,0} B_{0,0},
\end{align*}
\begin{align*}
\mathcal B_1^{(2)}(m,n)=&-\frac{1}{9} (m-n) \left(4 n A_{0,0} B_{0,0} H_m+6 n A_{0,0} B_{0,0} H_n\right.\\
&\left.-4 n A_{0,0} B_{0,0} \psi ^{(0)}(m+n+1)-4 n A_{0,1} B_{0,0}-3 n A_{0,0} B_{0,1}-2 A_{0,0} B_{0,0}\right),\\
\mathcal B_2^{(2)}(m,n)=&\ \frac{1}{9} \left(4 n A_{0,0} B_{0,0} H_m+6 n A_{0,0} B_{0,0} H_n\right.\\
&\left.-4 n A_{0,0} B_{0,0} \psi ^{(0)}(m+n+1)-4 n A_{0,1} B_{0,0}-3 n A_{0,0} B_{0,1}-2 A_{0,0} B_{0,0}\right),\\
\mathcal B_3^{(2)}(m,n)=& \frac{1}{9} \left(-4 A_{0,0} B_{0,0} H_m-6 A_{0,0} B_{0,0} H_n\right.\\
&\left. +4 A_{0,0} B_{0,0} \psi ^{(0)}(m+n+1)+4 A_{0,1} B_{0,0}+3 A_{0,0} B_{0,1}\right),
\end{align*}
\begin{align*}
\mathcal B_1^{(3)}(m,n)=&\frac{1}{18} \left(-24 n^2 A_{0,0} B_{0,0} H_m H_n\right.\\
&+12 n^2 A_{0,0} B_{0,1} H_m+8 n^2 A_{0,0} B_{0,0} H_m \psi ^{(0)}(m+n+1)\\
&+18 n^2 A_{0,0} B_{0,0} H_n \psi ^{(0)}(m+n+1)+16 n A_{0,0} B_{0,0} H_m\\
&+24 m n A_{0,0} B_{0,0} H_m H_n-24 m n A_{0,1} B_{0,0} H_n\\
&-12 m n A_{0,0} B_{0,1} H_m-6 m A_{0,0} B_{0,0} H_n-8 m n A_{0,0} B_{0,0} H_m \psi ^{(0)}(m+n+1)\\
&-18 m n A_{0,0} B_{0,0} H_n \psi ^{(0)}(m+n+1)-8 m A_{0,0} B_{0,0} H_m+24 n^2 A_{0,1} B_{0,0} H_n\\
&-5 n^2 A_{0,0} B_{0,0} \psi ^{(0)}(m+n+1)^2-5 n^2 A_{0,0} B_{0,0} \psi ^{(1)}(m+n+1)\\
&-8 n^2 A_{0,1} B_{0,0} \psi ^{(0)}(m+n+1)-9 n^2 A_{0,0} B_{0,1} \psi ^{(0)}(m+n+1)\\
&+12 m n A_{0,1} B_{0,1}+9 m n A_{0,0} B_{n,2}+5 m n A_{0,0} B_{0,0} \psi ^{(0)}(m+n+1)^2\\
&-8 n A_{0,0} B_{0,0} \psi ^{(0)}(m+n+1)+5 m n A_{0,0} B_{0,0} \psi ^{(1)}(m+n+1)\\
&+8 m n A_{0,1} B_{0,0} \psi ^{(0)}(m+n+1)+9 m n A_{0,0} B_{0,1} \psi ^{(0)}(m+n+1)\\
&+6 m A_{0,0} B_{0,0} \psi ^{(0)}(m+n+1)+8 m A_{0,1} B_{0,0}+3 m A_{0,0} B_{0,1}\\
&-12 n^2 A_{0,1} B_{0,1}-9 n^2 A_{0,0} B_{n,2}-16 n A_{0,1} B_{0,0}+2 A_{0,0} B_{0,0}\left.\right),
\end{align*}
\begin{align*}
\mathcal B_2^{(3)}(m,n)=&\frac{1}{18} \left(-24 n A_{0,0} B_{0,0} H_m H_n+12 n A_{0,0} B_{0,1} H_m\right.\\
&+8 n A_{0,0} B_{0,0} H_m \psi ^{(0)}(m+n+1)+18 n A_{0,0} B_{0,0} H_n \psi ^{(0)}(m+n+1)\\
&+8 A_{0,0} B_{0,0} H_m+6 A_{0,0} B_{0,0} H_n+24 n A_{0,1} B_{0,0} H_n\\
&-5 n A_{0,0} B_{0,0} \psi ^{(0)}(m+n+1)^2-6 A_{0,0} B_{0,0} \psi ^{(0)}(m+n+1)\\
&-8 n A_{0,1} B_{0,0} \psi ^{(0)}(m+n+1)-9 n A_{0,0} B_{0,1} \psi ^{(0)}(m+n+1)\\
&-5 n A_{0,0} B_{0,0} \psi ^{(1)}(m+n+1)-12 n A_{0,1} B_{0,1}-9 n A_{0,0} B_{n,2}\\
&-8 A_{0,1} B_{0,0}-3 A_{0,0} B_{0,1}\left.\right),
\end{align*}
\begin{align*}
\mathcal B_3^{(3)}(m,n)=&\frac{1}{18} \left(24 A_{0,0} B_{0,0} H_m H_n-8 A_{0,0} B_{0,0} H_m \psi ^{(0)}(m+n+1)\right.\\
&-18 A_{0,0} B_{0,0} H_n \psi ^{(0)}(m+n+1)-12 A_{0,0} B_{0,1} H_m\\
&-24 A_{0,1} B_{0,0} H_n+5 A_{0,0} B_{0,0} \psi ^{(0)}(m+n+1)^2+8 A_{0,1} B_{0,0} \psi ^{(0)}(m+n+1)\\
&+9 A_{0,0} B_{0,1} \psi ^{(0)}(m+n+1)+5 A_{0,0} B_{0,0} \psi ^{(1)}(m+n+1)\\
&+9 A_{0,0} B_{n,2}+12 A_{0,1} B_{0,1}\left.\right),
\end{align*}
\begin{align*}
\mathcal B_1^{(4)}(m,n)=&\frac{1}{18} \left(-n^2 A_{0,0} B_{0,0} \psi ^{(0)}(m+n+1)^3+m n A_{0,0} B_{0,0} \psi ^{(0)}(m+n+1)^3\right.\\
&+2 m A_{0,0} B_{0,0} \psi ^{(0)}(m+n+1)^2-3 n A_{0,0} B_{0,0} \psi ^{(0)}(m+n+1)^2\\
&+2 n^2 H_m A_{0,0} B_{0,0} \psi ^{(0)}(m+n+1)^2-2 m n H_m A_{0,0} B_{0,0} \psi ^{(0)}(m+n+1)^2\\
&+6 n^2 H_n A_{0,0} B_{0,0} \psi ^{(0)}(m+n+1)^2-6 m n H_n A_{0,0} B_{0,0} \psi ^{(0)}(m+n+1)^2\\
&-2 n^2 A_{0,1} B_{0,0} \psi ^{(0)}(m+n+1)^2+2 m n A_{0,1} B_{0,0} \psi ^{(0)}(m+n+1)^2\\
&-3 n^2 A_{0,0} B_{0,1} \psi ^{(0)}(m+n+1)^2+3 m n A_{0,0} B_{0,1} \psi ^{(0)}(m+n+1)^2\\
&-4 m H_m A_{0,0} B_{0,0} \psi ^{(0)}(m+n+1)+8 n H_m A_{0,0} B_{0,0} \psi ^{(0)}(m+n+1)\\
&-6 m H_n A_{0,0} B_{0,0} \psi ^{(0)}(m+n+1)+6 n H_n A_{0,0} B_{0,0} \psi ^{(0)}(m+n+1)\\
&-12 n^2 H_m H_n A_{0,0} B_{0,0} \psi ^{(0)}(m+n+1)\\
&+12 m n H_m H_n A_{0,0} B_{0,0} \psi ^{(0)}(m+n+1)\\
&-3 n^2 \psi ^{(1)}(m+n+1) A_{0,0} B_{0,0} \psi ^{(0)}(m+n+1)\\
&+3 m n \psi ^{(1)}(m+n+1) A_{0,0} B_{0,0} \psi ^{(0)}(m+n+1)\\
&+4 m A_{0,1} B_{0,0} \psi ^{(0)}(m+n+1)-8 n A_{0,1} B_{0,0} \psi ^{(0)}(m+n+1)\\
&+12 n^2 H_n A_{0,1} B_{0,0} \psi ^{(0)}(m+n+1)-12 m n H_n A_{0,1} B_{0,0} \psi ^{(0)}(m+n+1)\\
&+3 m A_{0,0} B_{0,1} \psi ^{(0)}(m+n+1)-3 n A_{0,0} B_{0,1} \psi ^{(0)}(m+n+1)\\
&+6 n^2 H_m A_{0,0} B_{0,1} \psi ^{(0)}(m+n+1)-6 m n H_m A_{0,0} B_{0,1} \psi ^{(0)}(m+n+1)\\
&-6 n^2 A_{0,1} B_{0,1} \psi ^{(0)}(m+n+1)+6 m n A_{0,1} B_{0,1} \psi ^{(0)}(m+n+1)\\
&-9 n^2 A_{0,0} B_{n,2} \psi ^{(0)}(m+n+1)+9 m n A_{0,0} B_{n,2} \psi ^{(0)}(m+n+1)\\
&+4 H_m A_{0,0} B_{0,0}+12 m H_m H_n A_{0,0} B_{0,0}-24 n H_m H_n A_{0,0} B_{0,0}\\
&-6 H_n A_{0,0} B_{0,0}+2 m \psi ^{(1)}(m+n+1) A_{0,0} B_{0,0}\\
&-3 n \psi ^{(1)}(m+n+1) A_{0,0} B_{0,0}+2 n^2 H_m \psi ^{(1)}(m+n+1) A_{0,0} B_{0,0}\\
&-2 m n H_m \psi ^{(1)}(m+n+1) A_{0,0} B_{0,0}\\
&+6 n^2 H_n \psi ^{(1)}(m+n+1) A_{0,0} B_{0,0}-6 m n H_n \psi ^{(1)}(m+n+1) A_{0,0} B_{0,0}\\
&-n^2 \psi ^{(2)}(m+n+1) A_{0,0} B_{0,0}+m n \psi ^{(2)}(m+n+1) A_{0,0} B_{0,0}\\
&-12 m H_n A_{0,1} B_{0,0}+24 n H_n A_{0,1} B_{0,0}-2 n^2 \psi ^{(1)}(m+n+1) A_{0,1} B_{0,0}\\
&+2 m n \psi ^{(1)}(m+n+1) A_{0,1} B_{0,0}-4 A_{0,1} B_{0,0}-6 m H_m A_{0,0} B_{0,1}\\
&+12 n H_m A_{0,0} B_{0,1}-3 n^2 \psi ^{(1)}(m+n+1) A_{0,0} B_{0,1}\\
&+3 m n \psi ^{(1)}(m+n+1) A_{0,0} B_{0,1}+3 A_{0,0} B_{0,1}+6 m A_{0,1} B_{0,1}\\
&-12 n A_{0,1} B_{0,1}+9 n A_{0,0} B_{n,2}+18 n^2 H_m A_{0,0} B_{n,2}\\
&-18 m n H_m A_{0,0} B_{n,2}-18 n^2 A_{0,1} B_{n,2}+18 m n A_{0,1} B_{n,2}\left.\right),
\end{align*}
\begin{align*}
\mathcal B_2^{(4)}(m,n)=&\ \frac{1}{18} \left(-12 A_{0,0} B_{0,0} H_m H_n\right.\\
&+18 n A_{0,0} H_m B_{n,2}+2 n A_{0,0} B_{0,0} H_m \psi ^{(0)}(m+n+1)^2\\
&+6 n A_{0,0} B_{0,0} H_n \psi ^{(0)}(m+n+1)^2+4 A_{0,0} B_{0,0} H_m \psi ^{(0)}(m+n+1)\\
&-12 n A_{0,0} B_{0,0} H_m H_n \psi ^{(0)}(m+n+1)+6 A_{0,0} B_{0,0} H_n \psi ^{(0)}(m+n+1)\\
&+12 n A_{0,1} B_{0,0} H_n \psi ^{(0)}(m+n+1)+6 n A_{0,0} B_{0,1} H_m \psi ^{(0)}(m+n+1)\\
&+2 n A_{0,0} B_{0,0} H_m \psi ^{(1)}(m+n+1)+6 n A_{0,0} B_{0,0} H_n \psi ^{(1)}(m+n+1)\\
&+6 A_{0,0} B_{0,1} H_m+12 A_{0,1} B_{0,0} H_n-n A_{0,0} B_{0,0} \psi ^{(0)}(m+n+1)^3\\
&-2 A_{0,0} B_{0,0} \psi ^{(0)}(m+n+1)^2-2 n A_{0,1} B_{0,0} \psi ^{(0)}(m+n+1)^2\\
&-3 n A_{0,0} B_{0,1} \psi ^{(0)}(m+n+1)^2-3 n A_{0,0} B_{0,0} \psi ^{(1)}(m+n+1) \psi ^{(0)}(m+n+1)\\
&-4 A_{0,1} B_{0,0} \psi ^{(0)}(m+n+1)-3 A_{0,0} B_{0,1} \psi ^{(0)}(m+n+1)\\
&-6 n A_{0,1} B_{0,1} \psi ^{(0)}(m+n+1)-9 n A_{0,0} B_{n,2} \psi ^{(0)}(m+n+1)\\
&-2 A_{0,0} B_{0,0} \psi ^{(1)}(m+n+1)-n A_{0,0} B_{0,0} \psi ^{(2)}(m+n+1)\\
&-2 n A_{0,1} B_{0,0} \psi ^{(1)}(m+n+1)-3 n A_{0,0} B_{0,1} \psi ^{(1)}(m+n+1)\\
&-18 n A_{0,1} B_{n,2}-6 A_{0,1} B_{0,1}\left.\right),
\end{align*}
\begin{align*}
\mathcal B_3^{(4)}(m,n)=&\ \frac{1}{18} \left(-18 A_{0,0} H_m B_{n,2}\right.\\
&-2 A_{0,0} B_{0,0} H_m \psi ^{(0)}(m+n+1)^2-6 A_{0,0} B_{0,0} H_n \psi ^{(0)}(m+n+1)^2\\
&+12 A_{0,0} B_{0,0} H_m H_n \psi ^{(0)}(m+n+1)-12 A_{0,1} B_{0,0} H_n \psi ^{(0)}(m+n+1)\\
&-6 A_{0,0} B_{0,1} H_m \psi ^{(0)}(m+n+1)-2 A_{0,0} B_{0,0} H_m \psi ^{(1)}(m+n+1)\\
&-6 A_{0,0} B_{0,0} H_n \psi ^{(1)}(m+n+1)+A_{0,0} B_{0,0} \psi ^{(0)}(m+n+1)^3\\
&+2 A_{0,1} B_{0,0} \psi ^{(0)}(m+n+1)^2+3 A_{0,0} B_{0,1} \psi ^{(0)}(m+n+1)^2\\
&+3 A_{0,0} B_{0,0} \psi ^{(1)}(m+n+1) \psi ^{(0)}(m+n+1)+6 A_{0,1} B_{0,1} \psi ^{(0)}(m+n+1)\\
&+9 A_{0,0} B_{n,2} \psi ^{(0)}(m+n+1)+A_{0,0} B_{0,0} \psi ^{(2)}(m+n+1)\\
&+2 A_{0,1} B_{0,0} \psi ^{(1)}(m+n+1)+3 A_{0,0} B_{0,1} \psi ^{(1)}(m+n+1)\\
&+18 A_{0,1} B_{n,2}\left.\right).
\end{align*}

\newpage

\bibliographystyle{alpha}
\bibliography{biblio}

\end{document}

%% file: cyclic.pdf_tex
\begingroup%
  \makeatletter%
  \providecommand\color[2][]{%
    \errmessage{(Inkscape) Color is used for the text in Inkscape, but the package 'color.sty' is not loaded}%
    \renewcommand\color[2][]{}%
  }%
  \providecommand\transparent[1]{%
    \errmessage{(Inkscape) Transparency is used (non-zero) for the text in Inkscape, but the package 'transparent.sty' is not loaded}%
    \renewcommand\transparent[1]{}%
  }%
  \providecommand\rotatebox[2]{#2}%
  \ifx\svgwidth\undefined%
    \setlength{\unitlength}{416.20472699bp}%
    \ifx\svgscale\undefined%
      \relax%
    \else%
      \setlength{\unitlength}{\unitlength * \real{\svgscale}}%
    \fi%
  \else%
    \setlength{\unitlength}{\svgwidth}%
  \fi%
  \global\let\svgwidth\undefined%
  \global\let\svgscale\undefined%
  \makeatother%
  \begin{picture}(1,0.36447281)%
    \put(0,0){\includegraphics[width=\unitlength,page=1]{cyclic.pdf}}%
    \put(0.46997204,0.01034196){\color[rgb]{0,0,0}\makebox(0,0)[lb]{\smash{$\P^1\times M$}}}%
    \put(0.16243108,0.03917392){\color[rgb]{0,0,0}\makebox(0,0)[lb]{\smash{$M_0$}}}%
    \put(0.86263592,0.06251408){\color[rgb]{0,0,0}\makebox(0,0)[lb]{\smash{$M_\infty$}}}%
    \put(0.45075073,0.1929444){\color[rgb]{0,0,0}\makebox(0,0)[lb]{\smash{$\mc P_\La$}}}%
    \put(0,0){\includegraphics[width=\unitlength,page=2]{cyclic.pdf}}%
    \put(0.91918777,0.19431731){\color[rgb]{0,0,0}\makebox(0,0)[lb]{\smash{$\{\infty\}\times \mc I^\infty_\La$}}}%
    \put(0.00603616,0.19569028){\color[rgb]{0,0,0}\makebox(0,0)[lb]{\smash{$\{0\}\times \mc I^0_\La$}}}%
  \end{picture}%
\endgroup%

%% file: hankel2.pdf_tex
\begingroup%
  \makeatletter%
  \providecommand\color[2][]{%
    \errmessage{(Inkscape) Color is used for the text in Inkscape, but the package 'color.sty' is not loaded}%
    \renewcommand\color[2][]{}%
  }%
  \providecommand\transparent[1]{%
    \errmessage{(Inkscape) Transparency is used (non-zero) for the text in Inkscape, but the package 'transparent.sty' is not loaded}%
    \renewcommand\transparent[1]{}%
  }%
  \providecommand\rotatebox[2]{#2}%
  \ifx\svgwidth\undefined%
    \setlength{\unitlength}{199.49865455bp}%
    \ifx\svgscale\undefined%
      \relax%
    \else%
      \setlength{\unitlength}{\unitlength * \real{\svgscale}}%
    \fi%
  \else%
    \setlength{\unitlength}{\svgwidth}%
  \fi%
  \global\let\svgwidth\undefined%
  \global\let\svgscale\undefined%
  \makeatother%
  \begin{picture}(1,0.44658303)%
    \put(0,0){\includegraphics[width=\unitlength,page=1]{hankel2.pdf}}%
    \put(0.71820143,0.33986267){\color[rgb]{0,0,0}\makebox(0,0)[lb]{\smash{$\gamma$}}}%
  \end{picture}%
\endgroup%

%% file: stokes.pdf_tex
\begingroup%
  \makeatletter%
  \providecommand\color[2][]{%
    \errmessage{(Inkscape) Color is used for the text in Inkscape, but the package 'color.sty' is not loaded}%
    \renewcommand\color[2][]{}%
  }%
  \providecommand\transparent[1]{%
    \errmessage{(Inkscape) Transparency is used (non-zero) for the text in Inkscape, but the package 'transparent.sty' is not loaded}%
    \renewcommand\transparent[1]{}%
  }%
  \providecommand\rotatebox[2]{#2}%
  \ifx\svgwidth\undefined%
    \setlength{\unitlength}{434.67016646bp}%
    \ifx\svgscale\undefined%
      \relax%
    \else%
      \setlength{\unitlength}{\unitlength * \real{\svgscale}}%
    \fi%
  \else%
    \setlength{\unitlength}{\svgwidth}%
  \fi%
  \global\let\svgwidth\undefined%
  \global\let\svgscale\undefined%
  \makeatother%
  \begin{picture}(1,0.26672593)%
    \put(0,0){\includegraphics[width=\unitlength,page=1]{stokes.pdf}}%
    \put(0.41380044,0.13755538){\color[rgb]{0,0,0}\makebox(0,0)[lb]{\smash{$\frac{\pi}{3} $}}}%
    \put(0.69905641,0.1322969){\color[rgb]{0,0,0}\makebox(0,0)[lb]{\smash{$R_{12}$}}}%
    \put(0.7345962,0.20197206){\color[rgb]{0,0,0}\makebox(0,0)[lb]{\smash{$R_{13}$}}}%
    \put(0.78969794,0.23746696){\color[rgb]{0,0,0}\makebox(0,0)[lb]{\smash{$R_{14}$}}}%
    \put(0.9616662,0.19539891){\color[rgb]{0,0,0}\makebox(0,0)[lb]{\smash{$R_{23}$}}}%
    \put(0.92102539,0.23352309){\color[rgb]{0,0,0}\makebox(0,0)[lb]{\smash{$R_{24}$}}}%
    \put(0.84609173,0.24272543){\color[rgb]{0,0,0}\makebox(0,0)[lb]{\smash{$R_{34}$}}}%
  \end{picture}%
\endgroup%